\documentclass[11pt,a4paper]{book}

\usepackage[margin=1.0in]{geometry}
\usepackage[utf8]{inputenc}
\usepackage[T1]{fontenc}
\usepackage{graphicx,physics,amsmath,amsfonts,amssymb,amsthm,amsthm,tikz-cd,mathtools,tikz-cd,microtype,booktabs,siunitx,tikz,cancel,titlesec,xcolor,stmaryrd,hyperref,longtable}
\usetikzlibrary{positioning}
\usepackage[all]{xy}
\usepackage{cmbright,newtxmath}

\newcommand{\eq}{\begin{equation}}
\newcommand{\qe}{\end{equation}}
\renewcommand{\eq}{\begin{equation}\noindent}
\newcommand{\tonde}[1]{\left( #1 \right) }
\newcommand{\quadre}[1]{\left[ #1 \right] }
\newcommand{\graffe}[1]{\left\{   #1  \right\} }
\renewcommand{\qe}{\end{equation}}
\newcommand{\eqa}[1]{\begin{equation}\noindent\begin{aligned} #1 \end{aligned}\end{equation}}
\newcommand{\M}{\mathcal}

\newcommand{\B}{\mathbb} 

\renewcommand{\span}{\mathrm{span}}
\newcommand{\st}{ \ : \ }

\newcommand{\Hom}{\mathrm{Hom}}

\newcommand{\id}{\operatorname{id}}

\newcommand{\ver}{\mathfrak{ver}}
\newcommand{\hor}{\mathfrak{hor}}
\newcommand{\ad}{\mathrm{ad}}
\newcommand{\flip}{\mathrm{flip}}

\newtheoremstyle{note}
{12pt}
{12pt}
{}
{}
{\bfseries}
{.}
{.5em}
{}
\theoremstyle{note} 

\theoremstyle{plain}
\newtheorem{thm}{Theorem}[section] 

\newtheorem{prop}[thm]{Proposition}
\newtheorem{lemma}[thm]{Lemma} 
\theoremstyle{definition}
\newtheorem{defn}[thm]{Definition} 
\newtheorem{exmp}[thm]{Example} 

\newtheorem{notation}[thm]{Notation}
\theoremstyle{remark}
\newtheorem{rmk}[thm]{Remark} 

\titleformat{\chapter}[display]
{\normalfont\huge\bfseries\centering}{\chaptertitlename\ \thechapter}{20pt}{\Huge}[{\titlerule[0.8pt]}]

\author{Antonio Del Donno}
\date{\today}
\title{Thesis notes}

\begin{document}
	\begin{titlepage}
		\begin{center}
			{
				{\Large
					{\textsc{
						Alma Mater Studiorum $\cdot$ Università di Bologna
					}}
				}
			}
			\rule[0.1cm]{15.8cm}{0.1mm}
			\rule[0.5cm]{15.8cm}{0.6mm}
			\medskip
			{\small
				{ 
					SCUOLA DI SCIENZE\\
					Corso di Laurea Magistrale in Matematica 
				}
			}
			\end{center}
			\vspace{15mm}
			\begin{center}
				{\huge{ \bf Differential calculi on }}\\
				\vspace{3mm}
				{\huge{ \bf quantum principal bundles }}\\
				\vspace{3mm}
				{\huge{\bf  in the \DJ ur\dj evi\'c approach}}\\
				\vspace{19mm}
				{\large
					{
						M.Sc. Thesis in Geometry
					}
				}
			\end{center}
			\vspace{50mm}
			\par
			\noindent
			\begin{minipage}[t]{0.47\textwidth}
				{\large
					{
						Supervisor:\\
							\textbf{Prof.
					Emanuele Latini} \\ 
						Co-Supervisor: \\ 
						\textbf{Dott.
						Thomas Weber}
						
					}
				}
			\end{minipage}
			\hfill
			\begin{minipage}[t]{0.5\textwidth}\raggedleft
    \large
    {
          \begin{flushright}
          	 Author:\\ 
      \hfill  \textbf{Antonio Del Donno}  \end{flushright}
    }
\end{minipage}
			\vspace{50mm}
			\begin{center}
			{\large
				{
					A.A. 2022/2023 
				}
			}
		\end{center}
	\end{titlepage}
	\thispagestyle{empty}
	\cleardoublepage

		\thispagestyle{empty}
	\thispagestyle{empty}
	\section*{\centering \textsc{Abstract}}

In this thesis we study the \DJ ur\dj evi\'c theory of differential calculi on quantum principal bundles within the domain of noncommutative geometry. Throughout the exposition, an algebraic approach based on Hopf algebras is employed. We begin by briefly recalling the foundational concepts of Hopf algebras and comodule algebras. Hopf-Galois extensions are introduced, along with an important correspondence result which relates crossed product algebras and cleft extensions. Differential calculus over algebras is reviewed, with a particular focus on the covariant case. The \DJ ur\dj evi\'c theory is presented in a modern language. We compare this theory to existing literature on the topic. We extend the \DJ ur\dj evi\'c theory and provide explicit realisations of quantum principal bundles and complete differential calculi on the noncommutative algebraic 2-torus, the quantum Hopf fibration and on crossed product algebras, resulting in an original contribution within this thesis.  
	\section*{\centering \textsc{Sommario}}

In questa tesi presentiamo la teoria di \DJ ur\dj evi\'c dei calcoli differenziali su fibrati principali quantistici nel dominio della geometria noncommutativa. Durante l'esposizione, adottiamo un approccio algebraico basato sulle algebre di Hopf. Iniziamo brevemente ricordando i concetti fondamentali delle algebre di Hopf e delle comodule algebras. Introduciamo le estensioni di Hopf-Galois, insieme ad un importante risultato di corrispondenza che collega le crossed product algebras e le estensioni cleft. Esaminiamo il calcolo differenziale sulle algebre, con particolare attenzione al caso covariante. La teoria di \DJ ur\dj evi\'c viene presentata adottando un linguaggio moderno. Confrontiamo questa teoria con la letteratura esistente sull'argomento. Estendiamo la teoria di \DJ ur\dj evi\'c, e provvediamo alcune realizzazioni esplicite di fibrati principali quantistici e calcoli differenziali completi sul 2-toro noncommutativo, sulla fibrazione di Hopf quantistica e sulle crossed product algebras, risultando in un contributo originale all'interno di questa tesi. 
 			\thispagestyle{empty}
	\clearpage{\thispagestyle{empty}\cleardoublepage}\thispagestyle{empty} \newpage\thispagestyle{empty}  \begin{center}
    \thispagestyle{empty}
    \vspace*{\fill}
		\thispagestyle{empty}
   A Lea. 
    \vspace*{\fill}
  \clearpage{\thispagestyle{empty}}
\end{center} 

		\frontmatter\tableofcontents\thispagestyle{empty}	 \newpage \thispagestyle{empty}	
		
 \chapter{Introduction}    		 \setcounter{page}{1}

To do geometry, instead of working with points of a manifold we may  consider the commutative algebra of $C^{\infty}$ functions on such. Noncommutative geometry extends this idea considering such algebras to be possibly noncommutative \cite{Connes}. The corresponding underlying "space" is called a "quantum space", borrowing the terminology from physics, where the process of quantisation consists in turning the commutative algebra of observables on the phase space into a noncommutative one. Here we explore aspects of noncommutative geometry in which a fundamental role is played by "quantum groups", i.e. the deformed function algebras on Lie groups \cite{woronowicz,Drinfeld}, and more in general Hopf algebras, commonly accepted as the natural analogue of a group in this noncommutative setting. These structures originally appeared as generalisation of symmetry groups for certain physical integrable systems \cite{Faddev}. Moreover they appear in the noncommutative generalisation of principal bundles of classical differential geometry, known as Hopf-Galois extensions, or principal comodule algebras, originally in \cite{BrzMjd,Hajac} and more recently in \cite{fioresiaschierilatiniweber,reamonn}. The total and base space are understood as algebras on which the quantum group acts in a precise sense. Our main motivation to explore this subject is to extend the differential geometry of homogeneous spaces to the noncommutative case. 
Within this noncommutative setting the mathematical theory of quantum groups its undoubtedly interesting on its own, but since geometry originally emerged as a practical subject we would also like to provide some motivation coming from physics. The quantum theory of fields gives a description of fundamental particles and interactions in terms of quantum fields defined on a suitable Lorentzian manifold. This theory culminates in the "Standard Model", a gauge theory which gives a quantum description of almost all the known fundamental particles. Interactions are given by gauge fields (or gauge mediators) that  mathematically are understood as connections on certain principal bundles. A quantum description of gravity is ruled out in this framework owing to non-renormalisability of the latter. Efforts to reconcile gravity with quantum field theory have led to the search for a quantum theory of gravity, aiming to unify general relativity with quantum field theory. One notable approach of investigation lies in the concept of gravity as a gauge theory within the Cartan moving frame setting, where connections on certain principal bundles play a central role \cite{wise}. 
The commonly accepted general indication is that the small-scale structure of space-time could not be reasonably modelled according to usual continuum geometry. Geometry at the Planck scale may be modified in order to account for the presence of quantum effects, and noncommutative geometry could provide a possible realisation of the latter. In fact, within this approach, geometry can be specialised to deal with discrete spaces or finite dimensional algebras, with the possibility of recovering key insights on physics beyond the Planck scale. The classical \textit{emergent} geometry of smooth (possibly Riemannian) manifolds is recovered as a \textit{low energy} approximation via a suitable limit process \cite{majid2000}. 

    \smallskip

In this thesis we study \DJ ur\dj evi\'c's theory of quantum principal bundles, originally presented in \cite{Durdevic-1,Durdevic-2}. \DJ ur\dj evi\'c introduces a comprehensive and natural theory of principal bundles within the framework of noncommutative geometry. In this theory, quantum groups act as the structure groups, while general quantum spaces serve as total space and base manifolds. The study focuses on developing a differential calculus specifically tailored for quantum principal bundles. This includes the introduction and examination of algebras of horizontal and vertical differential forms on the bundle. A natural "braiding" emerges, from which the formalism of connections is elaborated upon, with attention given to operators such as bimodule connections, metrics, Levi-Civita connections, covariant derivative, and curvature. A quantum counterparts of infinitesimal gauge transformations is given. Our aim is to present this theory in a modern light, offering three significant contributions. First, we translate the original papers into a more modern mathematical language, making the approach accessible to a wider audience. Second, we establish connections between \DJ ur\dj evi\'c's work and existing literature on quantum principal bundles, particularly engaging with the influential research of \cite{BrzMjd}. This highlights how the \DJ ur\dj evi\'c approach builds upon, rather than competes with, established frameworks, extending them to higher-order forms.
We introduce new examples of \DJ ur\dj evi\'c's quantum principal bundles, including bicovariant calculi, crossed product calculi, and the fundamental Hopf fibration on the quantum sphere. This integration with established methodologies, along with the examples provided, is indicative of the significance of \DJ ur\dj evi\'c's framework in noncommutative differential geometry.
Furthermore, we clarify and expand \DJ ur\dj evi\'c's findings, particularly those presented in the work \cite{Durdevic-2}. Doing so, we adopt an algebraic perspective, focusing on principal comodule algebras and Hopf algebras, without presuming the existence of compact quantum spaces or groups. While we omit certain technical details, such as considerations related to $*$-involutions, these aspects remain ripe for future investigation. Moreover, exploration into the braiding, connections, and metrics associated with \DJ ur\dj evi\'c's work has also been postponed to upcoming inquiries.
 \section*{Outline} In Chapter \ref{Hopf algebras} we briefly introduce the main algebraic structures we deal with. We consider Hopf algebras $H$, and $H-$comodule algebras together with the corresponding subalgebra of (co)invariant elements. These algebras respectively play the role of structure group, total and base space of a quantum principal bundle provided some auxiliary conditions hold. We provide some explicit realisations of Hopf algebras and comodule algebras. 
 
 In Chapter \ref{Hopf-Galois extensions} we define Hopf-Galois, cleft and trivial extensions. Such extensions are of the form $B\subseteq A$, where $A$ is a (right) $H-$comodule algebra and $B$ the subalgebra of $H-$coinvariant elements \cite{Montgomery}. We discuss some properties of convolution invertible morphisms (cleaving maps) from  Hopf algebras to comodule algebras. It is proven that trivial extensions are cleft, and that cleft extensions are Hopf-Galois. Some explicit realisations of trivial, cleft and Hopf-Galois extensions are discussed. We introduce the translation map as an "inverse" of the canonical Hopf-Galois  map of a Hopf-Galois extension. Some properties of the translation map are discussed. Crossed product algebras $B\sharp_{\sigma}H$ are defined starting from an algebra $B$, satisfying some additional assumptions, and a Hopf algebra $H$.  We provide a correspondence result between cross product algebras $B\sharp_{\sigma}H$ and cleft extensions $B\subseteq B\sharp_{\sigma}H$ \cite{Doi-Takeuchi}. 
 
  In Chapter \ref{Differential calculus over algebras} we begin discussing differential geometry in the noncommutative setting. We define the notion of first order differential calculus over an algebra and provide some explicit realisations. Such structures are non-unique, in the sense that many different first order differential calculi can be defined on the same algebra. We introduce the first order universal differential calculus and provide a theorem \cite{woronowicz} stating how every first order calculus over an algebra can be induced as a quotient of the latter. The notion of first order covariant calculi over algebras is introduced \cite{Connes}. We provide a classification result of those first order differential calculi which are covariant \cite{woronowicz}, i.e. compatible with a coaction on the algebra. We consider the extension of the first order theory to higher order differential calculi. The universal differential calculus is introduced. Higher order covariant differential calculi are discussed. We consider the maximal prolongation \cite{Durdevic-1}, being the biggest differential calculus that can be defined starting from a given first order differential calculus. We prove how the maximal prolongation of the universal first order differential calculus is its tensor algebra, i.e. the universal differential calculus. Moreover, every differential calculus  turns out to be induced as a quotient of the universal. 
   
    In Chapter \ref{Quantum principal bundles chap}  we delve into \DJ ur\dj evi\'c's theory of quantum principal bundles. We start with a bicovariant first order differential calculus $\Gamma$ over a Hopf algebra $H$. The corresponding maximal prolongation $\Gamma^{\wedge}$ is considered. The (differential graded) subalgebra of coinvariant elements $\Lambda^{\bullet}$ is introduced. In \DJ ur\dj evi\'c's framework the definition of quantum principal bundle is that of faithfully flat Hopf-Galois extensions $B\subseteq A$. Vertical forms $\ver^{\bullet}(A)$ over the total space $A$ of the bundle are introduced. We show that these form a differential calculus. The notion of complete differential calculus $\Omega^{\bullet}(A)$ over the total space is introduced. In such scenario we provide a surjective morphism $\pi_{\ver}:\Omega^{\bullet}(A)\twoheadrightarrow \ver^{\bullet}(A)$ projecting a total space form onto the space of vertical forms. We show that the differential calculus of vertical forms is complete if the calculus on the total space is. The noncommutative analogue of horizontal forms $\hor^{\bullet}(A)$ is introduced. It is shown that these form a graded subalgebra of $\Omega^{\bullet}(A)$. In particular $\hor^{\bullet}(A)$ is a right $H-$comodule algebra.  It is shown how first order horizontal, vertical and total space constitute a natural short exact sequence of $A-$modules $0\rightarrow \hor^{1}(A)\hookrightarrow\Omega^{1}(A) \twoheadrightarrow \ver^{1}(A) \rightarrow 0$. We observe that within other approaches exactness of this sequence is not guaranteed, see for example \cite{fioresiaschierilatiniweber}. It is discussed how the same sequence fails to be exact for higher order forms. The differential calculus over the base space $B:=A^{coH}$ is introduced. It is shown how differential forms $\Omega^{\bullet}(B)$ over the base space coincides with the intersection of right $H-$coinvariant and horizontal forms. We argue that $\Omega^{\bullet}(B)$ is not a differential calculus in general. However, in all explicit examples we consider the base calculus will be a differential calculus generated by $B$. We compare \DJ ur\dj evi\'c's theory of quantum principal bundles with the one in \cite{BrzMjd,beggs-majid}. Explicit realisations of complete differential calculi and quantum principal bundles are presented. It is shown how every Hopf algebra $H$ over a field $\B{k}$ naturally forms a quantum principal bundle $\B{k}\subseteq H$. A complete differential calculus is constructed along the "group algebra" $\B{C}(\B{Z})$.  We show that the noncommutative algebraic 2-torus forms a cleft extension, and so a quantum principal bundle, under the right $H-$coaction of the Hopf algebra $\M{O}(\text{U}(1))$. We construct a complete differential calculus over the total space of such bundle. We show how the space of right $H-$coinvariant forms $1-$forms of the bundle is generated by the base space itself, and so forms a differential calculus over the base.  We consider the quantum Hopf fibration with total space $A=\M{O}_{q}(\text{SU}(2))$ under the action of the group algebra $H=\M{O}_{q^{2}}(\text{U}(1))$. There is a complete differential calculus over the total space exists and that the corresponding forms over the base space are a differential calculus.  Finally, we define differential calculi over crossed product algebras $B\sharp_{\sigma}H$. According to Chapter \ref{Hopf-Galois extensions} it is clear that every crossed product algebra $B\sharp_{\sigma}H$ is a quantum principal bundle and that every cleft extension is of this form. We show that the differential calculus over the total space $B\sharp_{\sigma}H$ is complete given the differential calculus on $H$ is complete. Moreover the corresponding base forms are a differential calculus.
    
   \smallskip
    We point out that we make a new contribution by establishing completeness of the differential calculi on total spaces of the examples discussed. Moreover, within the same examples, we are able to recover differential calculi on the corresponding base spaces. The results of this thesis will be summarised in a forthcoming publication.
\section*{Acknowledgements}
My profound gratitude goes to my supervisors, Prof. Emanuele Latini and Dr. Thomas Martin Weber. If I have had such a wonderful time during these months of work it is also thanks to their patience, commitment and understanding. Their guidance has been invaluable and I am deeply thankful for all the seminars, fruitful conversations and overall time spent together.\clearpage

\section*{Notation}

\begin{table}[h]
\begin{tabular}{c p{8cm} c}

\textbf{Symbol} & \textbf{Description}  & \textbf{First time appearing}\\ 
\\ 
$A$             & An algebra/$H$-comodule algebra & Page 1  \\ 
$\Delta$        & The coproduct on a coalgebra/bialgebra/Hopf algebra & Page 1 \\ 
$\epsilon$      & The counit on a coalgebra/bialgebra/Hopf algebra  &  Page 1 \\ 
$\mu$           & The product on an algebra/bialgebra/Hopf algebra & Page 1 \\ 
$\eta$          & The unit on an algebra/bialgebra/Hopf algebra  &  Page 1 \\ 
$\mathbb{k}$    & A field   &  Page 1   \\ 
$H$             & A Hopf algebra  & Page 3\\ 
$S$             & The antipode on a Hopf algebra & Page 3 \\ 
$j$             & The cleaving map & Page 5 \\ 
$*$             & The convolution product & Page 5 \\  
$\triangleright$ & Right $H$-module action on a vector space/algebra & Page 7 \\ 
$\triangleleft$  & Left $H$-action on a vector space/algebra & Page 7 \\ 
${}_{V}\Delta$   & Left $H$-coaction on a vector space/algebra & Page 7 \\  
$\Delta_{V}$    & Right $H$-coaction on a vector space/algebra & Page 7 \\  
${}^{coH}A$     & The subalgebra of left $H$-coinvariant elements & Page 8 \\ 
$A^{coH}$       & The subalgebra of right $H$-coinvariant elements & Page 8 \\ 
$\otimes_{B}$   & The balanced tensor product over $B$ & Page 13 \\  
$B\subseteq A$  & A trivial/cleft/Hopf-Galois extension & Page 13 \\  
$\cdot$ (usually omitted) & Left/right $A$-module action & page 24 \\   	
\end{tabular}
\end{table}

\mainmatter \chapter{Hopf algebras}\label{Hopf algebras}

This chapter provides a brief introduction to the subject of Hopf algebras, being the fundamental "group" structures upon which we develop the theory of quantum principal bundles presented in this thesis.  We start by discussing concepts of algebra and coalgebras in \ref{algebras and coalgebras}. We define bialgebras and Hopf algebras in \ref{bialgebras and hopf algebras}, the latter being bialgebras equipped with an "antipode" which we show to be unique and satisfies the property of being an antibialgebra map. In section \ref{comodule and comodule algebras} we introduce the notion of $H-$comodule algebras, playing the role of total spaces in the theory of quantum principal bundles. The main references we will be using are \cite{Kassel,beggs-majid,Sweedler,Montgomery}. 

\section{Algebras and coalgebras} \label{algebras and coalgebras}
We provide the definition of algebra, expressing all structures as linear maps. 
\begin{defn}
	An algebra $A$ over a field $\B{k}$ is a $\B{k}-$vector space together with a product map $\mu:A\otimes A \rightarrow A$ and a unit element $1_{A}$ which can be equivalently written as a map $\eta:\B{k}\rightarrow A$ by $\eta(1)=1_{A}$ such that the following diagrams commute $$ \begin{tikzcd}[column sep=small]
& A\otimes A \otimes A \arrow[ld,"\mathrm{id}\otimes \mu"'] \arrow[rd,"\mu\otimes \mathrm{id}"] & \\
A\otimes A \arrow[dr,"\mu"'] & & A\otimes A \arrow[dl,"\mu"] \\
& A &
\end{tikzcd}
\quad
\begin{tikzcd}
A\otimes A \arrow[rrdd,"\mu"', shorten >=-2pt,swap] & \\ \\ 
\mathbb{k}\otimes A \arrow[uu,"\eta\otimes \mathrm{id}", shorten >=-2pt] \arrow[rr,"\cong"', shorten >=-2pt] & &  A
\end{tikzcd}
\quad
\begin{tikzcd}
A\otimes A \arrow[rrdd,"\mu"',swap] &  \\ \\
A\otimes \mathbb{k} \arrow[uu,"\mathrm{id}\otimes \eta"] \arrow[rr,"\cong"'] & & A
\end{tikzcd},$$ that is \eqa{ & \mu \circ (\id\otimes \mu)  = \mu \circ (\mu\otimes \id), \\ & \mu\circ (\eta  \otimes \id) = \id, \\ & \mu\circ (\id\otimes \eta) = \id.} 

\end{defn} Commutativity of these diagrams are representative of associativity and unitality, respectively. 
The notion of coalgebra is obtained via such diagrammatic approach by reversing arrows. \begin{defn}
	A coalgebra $C$ over a field $\B{k}$ is a $\B{k}-$vector space togehter with a coproduct map $\Delta:C\rightarrow C\otimes C$ and a counit map $\epsilon:C\rightarrow \B{k}$ such that the following diagrams commute $$ \begin{tikzcd}[column sep=small]
& C\otimes C \otimes C & \\
C\otimes C \arrow[ur,"\Delta\otimes \mathrm{id}"',swap] & & C\otimes C \arrow[lu,"\mathrm{id}\otimes \Delta"'] \\
& C \arrow[lu,"\Delta"',swap] \arrow[ru,"\Delta"'] &
\end{tikzcd}
\quad
\begin{tikzcd}
C\otimes C \arrow[dd,"\epsilon\otimes \mathrm{id}"']\\ \\  \mathbb{k}\otimes C & &  C \arrow[ll,"\cong"',swap] \arrow[uull,"\Delta"']
\end{tikzcd}
\quad
\begin{tikzcd}
C \otimes C \arrow[dd,"\mathrm{id}\otimes \epsilon"']\\  \\ \mathbb{k}\otimes C & & C \arrow[ll,"\cong"',swap] \arrow[uull,"\Delta"']
\end{tikzcd},$$ that is \eqa{& (\id\otimes \Delta)\circ \Delta= (\Delta\otimes \id)\circ \Delta, \\ & (\epsilon\otimes \id)\circ \Delta = \id, \\ &(\id\otimes \epsilon)\circ \Delta= \id.}

\end{defn}
Commutativity of these diagrams are representative of \textit{coassociativity} and \textit{counitality} of the coalgebra $C$. 
\begin{notation}\label{sweedler notation} When dealing with such algebraic structures it is customary to adopt the "Sweedler notation" We denote \eqa{\Delta(c) :=  c_{1}\otimes c_{2}, \quad \text{for all c}\in C,} which is a concise way of keeping track of the numbering of tensor factors. This is particularly useful to write down in a more explicit way the axiom of a coalgebra (an even more).\qed\end{notation} Coassociativity reads \eqa{ \nonumber (\Delta\otimes \id)\circ \Delta(c)& = (\Delta\otimes \id)\circ( c_{1}\otimes c_{2} )\\ & = (\Delta(c_{1})\otimes c_{2}) \\ & =  c_{11}\otimes c_{12} \otimes c_{2};\\ (\id\otimes \Delta)\circ \Delta(c) & = (\id\otimes \Delta)\circ (c_{1}\otimes c_{2}) \\ & = c_{1}\otimes \Delta(c_{2}) \\ & = c_{1}\otimes c_{21}\otimes c_{22}  ,} or \eqa{c_{11}\otimes c_{12}\otimes c_{22} =  c_{1}\otimes c_{2}\otimes c_{3} = c_{1}\otimes c_{21}\otimes c_{22}.} The convention is to relable indexes from the lowest to the highest\footnote{Here indexes $11$ and $12$ are of course relative to tensor product component preceeding the index $2$. Similarly for $21$ and $22$.}  keeping track of the order of factors.  Unitality reads \eqa{ \nonumber(\epsilon\otimes \id) \circ \Delta (c) & =  (\epsilon \otimes \id) \circ ( c_{1}\otimes c_{2}) \\ & = (\epsilon(c_{1})\otimes c_{2}) \\ & =  \epsilon (c_{1})c_{2} \\ & = c \\ & =  c_{1}\epsilon(c_{2}) \\ & = (\id\otimes \epsilon)\circ \Delta(c),} or \eqa{  \epsilon(c_{1})c_{2} = c=  c_{1}\epsilon(c_{2}).}  
\begin{exmp}
	We provide a short example to explain the Sweedler notation. Let $H$ be a Hopf algebra. We consider the map $(\Delta\otimes \Delta)\circ \Delta:H\rightarrow H\otimes H \otimes H \otimes H$. For $h\in H$ we find \eqa{\nonumber (\Delta\otimes \Delta)\circ \Delta(h) & = (\Delta\otimes \Delta) (h_{1}\otimes h_{2}) \\ & = \Delta(h_{1})\otimes \Delta(h_{2}) \\ & = h_{11}\otimes h_{12}\otimes h_{21}\otimes h_{22},} and relabel according to Notation \ref{sweedler notation} as $(\Delta\otimes \Delta)\circ \Delta(h)=h_{1}\otimes h_{2}\otimes h_{3}\otimes h_{4}$.  
\end{exmp}
\section{Bialgebras and Hopf algebras}\label{bialgebras and hopf algebras}
The next definition is the one of bialgebra,  featuring simultaneously the algebraic structure of an algebra and a coalgebra. \begin{defn}
	A bialgebra $H$ over $\B{k}$ is a $\B{k}-$vector space equipped with linear maps $\mu:H\otimes H \rightarrow H$, $\eta:\B{k}\rightarrow H$, and algebra maps $\Delta:H\rightarrow H\otimes H$, $\epsilon:H\rightarrow \B{k}$ such that $(H,\mu,\eta)$ is an algebra and $(H,\Delta,\epsilon)$ is a coalgebra. 
\end{defn} We now define a Hopf algebra starting with a  bialgebra with an additional requirement, that is a map playing the role of "group inversion". \begin{defn}A Hopf algebra is a bialgebra $(H,\Delta,\epsilon,\mu,\eta)$ together with an "antipode" map $S:H\rightarrow H$ such that the following diagram commutes $$ 
\begin{tikzcd}
                                   & H\otimes H  \arrow[rr,"S\otimes \id"] &                       & H \otimes H  \arrow[rd,"\mu"] &   \\
H \arrow[rd,"\Delta"] \arrow[ru,"\Delta"] \arrow[rr,"\epsilon"] &                        & \mathbb{k} \arrow[rr,"\eta"] &                         & H \\
                                   & H\otimes H \arrow[rr,"\id\otimes S"]  &                       & H \otimes H \arrow[ru,"\mu"]  &  
\end{tikzcd},
$$ that is to say \eqa{ \nonumber \mu\circ (S\otimes \id) \circ \Delta (h)& = \mu\circ (S(h_{1})\otimes h_{2}) \\ & =  S(h_{1})h_{2} \\ & = \eta\circ \epsilon(h) \\ & = \epsilon(h) \\ & =  h_{1}S(h_{2}) \\ & = \mu\circ (\id\otimes S)\circ \Delta(h) ,} or \eqa{  S(h_{1})h_{2}=\epsilon(h)= h_{1}S(h_{2}).}  
	
\end{defn}
\begin{defn} A Hopf algebra is commutative if the underlying algebra structure is. Moreover, defining $\text{flip}:H\otimes H \rightarrow H\otimes H$ as $\text{flip}(a\otimes b) = b\otimes a$ we say a Hopf algebra is cocommutative if $\text{flip}\circ \Delta=\Delta$, i.e. \eqa{ \text{flip}\circ \Delta(h) & = \text{flip}(h_{1}\otimes h_{2}) \\ & = h_{2}\otimes h_{1} \\ & = h_{1}\otimes h_{2},} or $h_{1}\otimes h_{2}=h_{2}\otimes h_{1}$ for every $h\in H$.  \end{defn} We provide a few explicit realisations of Hopf algebras by the following examples. \begin{exmp}\label{group algebra}
	Let $(G,\cdot,e)$ be a group and let $\B{k}$ be a field. We define the group algebra $H=\B{k}G$ as the $\B{k}-$vector space generated by elements of $G$. Accordingly every element in $H$ is of the form $\sum_{g\in G}k_{g}g$ with finitely many $k_{g}\in \B{k}$ being non-zero.  The algebra structure on $\B{k}G$ is given by the associative product \eqa{\sum_{g\in G} k_{g}g\cdot \sum_{h\in G} k'_{h}h = \sum_{g,h\in G} k_{g}k'_{h}g\cdot h.} The unit element is $e$. The coalgebra structure is extended linearly by the following maps acting on elements of $G$  \eqa{&\Delta: G \rightarrow G\otimes G,\quad && g\mapsto g\otimes g, \\ & \epsilon:G\rightarrow \B{k},\quad && g\mapsto 1} as \eqa{& \Delta\tonde{\sum_{g\in G}k_{g}g}= \sum_{g\in G}k_{g}g\otimes g, \quad\epsilon\tonde{\sum_{g\in G}k_{g}g} =\sum_{g\in G} k_{g}.}   
	Finally, the Hopf algebra structure follows defining a morphism $S:H\rightarrow H$ defined as $S(g)=g^{-1}$ on elements of $G$ and linearly extending to elements of $H$ as $S(\sum_{g\in G}k_{g}g) = \sum_{g\in G} k_{g} g^{-1}.$  It is an easy check that $(H,\Delta,\epsilon,\mu,\eta,S)$ is a Hopf algebra by those assignments. \qed  
\end{exmp}

\begin{exmp}
	Let $V$ be a $\B{k}-$vector space. Recall the tensor algebra $\M{T}V$ over $V$ is defined as \eqa{ \nonumber\M{T}(V):= \B{k}\oplus V \oplus (V\otimes V) \oplus\dots. }  We define \eqa{ & \Delta: \M{T}V \rightarrow \M{T}V\otimes \M{T}V, \quad && x\mapsto x\otimes 1 + 1 \otimes x \\& \epsilon:\M{T}V\rightarrow \B{k}, \quad && x \mapsto 0 \\ &S:\M{T}V\rightarrow \M{T}V,\quad && x\mapsto -x} for every $x\in V$. Extending$\Delta,\epsilon$ to algebra maps and $S$ as an antibialgebra map over $\M{T}V$ realises a Hopf algebra space.\qed 
\end{exmp}
\begin{exmp}
	Let $\B{k}=\B{C}$ and let $q\in \B{C}$ be non-zero. We define the "quantum group" $H=\text{SL}_{q}(2)$ as the free algebra generated as $H=\span_{\B{C}}\{\alpha,\beta,\gamma,\delta\}$ modulo relations \eqa{ & \alpha\beta=q \beta  \alpha, \quad \alpha\gamma=q\gamma\alpha, \quad \beta\delta =q\delta \beta, \quad \gamma\delta=q\delta\gamma,\\ & \beta\gamma=\gamma \delta,\quad \alpha\delta-\delta\alpha=(q-q^{-1})\beta\gamma, \quad \alpha\delta-q\beta\gamma=1.} Coproduct, counit and antipode are provided as \eqa{& \Delta:H\rightarrow H \otimes H, \quad && \begin{pmatrix}
		\alpha & \beta \\ \gamma & \delta 
	\end{pmatrix}\mapsto \begin{pmatrix}
		\alpha & \beta \\ \gamma & \delta 
	\end{pmatrix} \otimes \begin{pmatrix}
		\alpha & \beta \\ \gamma & \delta 
	\end{pmatrix} = \begin{pmatrix} \alpha\otimes \alpha + \beta\otimes \gamma & \alpha\otimes \beta + \beta \otimes \delta \\ \gamma\otimes \alpha + \delta\otimes \gamma & \gamma\otimes \beta + \delta\otimes \delta \end{pmatrix}, \\ &\epsilon:H\rightarrow \B{k},\quad && \begin{pmatrix}
		\alpha & \beta \\ \gamma & \delta 
	\end{pmatrix}\mapsto\begin{pmatrix}1 & 0 \\ 0 & 1\end{pmatrix}, \\ &S:H\rightarrow H, \quad && \begin{pmatrix}
		\alpha & \beta \\ \gamma & \delta 
	\end{pmatrix}\mapsto \begin{pmatrix}
		\delta & q^{-1}\beta \\ -q \gamma & \alpha
	\end{pmatrix}.}  \qed 
\end{exmp}

\begin{prop}
	Let $H$ be a Hopf algebra. Then its antipode $S:H\rightarrow H$ is unique. 
	\proof Let $h\in H$ and let $S_{1},S_{2}$ be two antipode maps. Then we find \eqa{\nonumber S_{1}(h) & = S_{1}(h_{1}\epsilon(h_{2})) \\ & = S_{1}(h_{1})\epsilon(h_{2}) \\ & = S_{1}(h_{1}) h_{21}S_{2}(h_{22}) \\ & = S_{1}(h_{1})h_{2}S_{2}(h_{3}) \\ & = \epsilon(h_{1})S_{2}(h_{2}) \\ & = S_{2}(\epsilon(h_{1})h_{2}) \\ & = S_{2}(h).}  \qed 
\end{prop}

Next we introduce the convolution algebra in order to prove that the antipode map $S:H\rightarrow H$ is a morphism of antibialgebras. The first definition is the one of convolution invertible morphism between algebras.  \begin{defn}
  	A $\B{k}-$linear map $j:H\longrightarrow A$ is said to be convolution invertible if we can find another $\B{k}-$linear map $i:H\longrightarrow A$ such that $j(h_{1})i(h_{2})=\epsilon(h)1_{A}=i(h_{1})j(h_{2})$, for every element $h\in H$.  In this case one usually writes $i=j^{-1}$ or vice-versa.  
  \end{defn} \begin{prop}[Convolution algebra]\label{convolution algebra}
 Given an associative unital algebra $(A,\mu,\eta)$ and a a coassociative counital algebra $(C,\Delta,\epsilon)$ consider the vector space $\text{Hom}_{\B{k}}(C,A)$ of $\B{k}-$linear maps $C\rightarrow A$. Define a $\B{k}-$linear map  \eqa{ *:\Hom_{\B{k}} (C,A)\otimes \Hom_{\B{k}}&(C,A)\longrightarrow \Hom_{\B{k}}(C,A)\\& (f,f')\mapsto f*f' := \mu\circ (f\otimes f')\circ \Delta .} We call $(\Hom_{\B{k}}(C,A),*)$ the \textit{convolution algebra}. It is s an associative unital algebra with unit $\eta\circ \epsilon$.
 \proof We already know $\Hom_{\B{k}}(C,A)$ is a vector space. Given an element $x\in C$ we find \begin{enumerate}
 \item $*$ is bilinear: \eqa{ (kf+hg)*(k'f'& +h'g')(x) \\ &  = \mu\circ\tonde{(kf+hg)\otimes(k'f'+h'g')}\otimes\Delta(x) \\ & = \mu\circ\tonde{ kf\otimes k'f' + kf\otimes h'g' + hg\otimes k'f + hg\otimes h'g'}\circ(x_{1}\otimes x_{2}) \\ & = \mu\circ (kf(x_{1})\otimes k'f'(x_{2}) + kf(x_{1})\otimes h'g'(x_{2})\\  & + hg(x_{1})\otimes k'f'(x_{2})+hg(x_{1})\otimes h'g'(x_{2})) \\ & = kk' f(x_{1})f'(x_{2}) + kh'f(x_{1})g'(x_{2})+ hk'g(x_{1})f'(x_{2})+ hh'g(x_{1})g'(x_{2}) \\ & = kk'(f*f')(x)+kh'(f*g')(x)+hk'(f*g)(x)+hh'(g*g')(x) \nonumber} for every $f,g,f',g' \in \Hom_{\B{k}}(C,A)$ and $k,k',h,h'\in\B{k}$.  
 \item $*$ is associative:  \eqa{ \nonumber \tonde{(f*g)*h}(x) & = \mu \circ \tonde{ (f*g)\otimes h)\Delta(x)} \\ & = \mu\circ \tonde{ (f*g)(x_{1})\otimes h(x_{2})} \\ & = \mu \circ \tonde{\mu\circ(f\otimes g)(x_{11}\otimes x_{12}) \otimes h(x_{2})} \\ & = \mu\circ\tonde{ f(x_{11})g(x_{12}) \otimes h(x_{2})} \\ &  = f(x_{1})g(x_{2})h(x_{3}) \\ & = \tonde{f*(g*h)}(x)}  
 \item $\eta\circ\epsilon$ is a unit: \eqa{ \nonumber\tonde{(\eta\circ \epsilon)*f}(x) & = \mu\circ\tonde{(\eta\circ \epsilon)\otimes f}(x_{1}\otimes x_{2}) \\ & = \mu\circ \tonde{\eta(\epsilon(x_{1}))\otimes f(x_{2})} \\ & = \eta(\epsilon(x_{1}))f(x_{2}) \\ & = f(x),} and similarly for the other component.  \qed   
 \end{enumerate}
  \end{prop}

  In view of Proposition \ref{convolution algebra} the notion of convolution invertibility amounts to invertibility in the convolution algebra.  
 \begin{prop}\label{convolution inverse in hopf algebras}
 	Let $H$ be a Hopf algebra. \begin{enumerate}
 		\item if $f:C\rightarrow H$ is a coalgebra map then $S\circ f$ is the convolution inverse of $f$ in $\Hom_{\B{k}}(C,H)$;
 		\item if $f:H\rightarrow A$ is an algebra map then $f\circ S$ is the convolution inverse of $f$ in $\Hom_{\B{k}}(H,A)$.  
 	\end{enumerate}
 	\proof It follows easily by a direct computation \eqa{ \nonumber(S\circ f * f)(x) & = \mu\circ\tonde{ S\circ f \otimes f)(x_{1}\otimes x_{2}} \\ & = \mu\circ \tonde{S(f(x_{1}))\otimes f(x_{2})} \\ & = S(f(x_{1}))f(x_{2}) \\ & = S(f(x)_{1})f(x)_{2}\\ & = \epsilon(f(x))1.}  Similarly the other point follows.  \qed 
 \end{prop}
 In view of the last two proposition we have the following.
\begin{prop}
	Let $H$ be a Hopf algebra. Then its antipode $S:H\rightarrow H$ satisfies \begin{enumerate}
		\item $S(hg)=S(g)S(h)$, for every $h,g\in H$;
		\item $S(1)=1$; 
		\item $(S\otimes S)\circ \Delta(h) = \text{flip}\circ \Delta \circ S(h)$, for every $h\in H$;
		\item $\epsilon\circ S(h)=\epsilon(h)$, for every $h\in H$. 
	\end{enumerate}
	\proof 
	In order: 
	\begin{enumerate}
		\item  We define $f:=S\circ \mu :H\otimes H \rightarrow H$ and $f':=\mu\circ\text{flip}\circ (S\otimes S)$. By Proposition \ref{convolution inverse in hopf algebras} we have that $f$ is the convolution inverse of $\mu$ in the convolution algebra $\Hom_{\B{k}}(H\otimes H, H)$. Moreover $f'$ is also a left inverse of $\mu:H\otimes H \rightarrow H$, indeed \eqa{\nonumber (f'* \mu)(h\otimes h') & = f'(h_{1}\otimes h'_{1})\mu(h_{2}\otimes h'_{2}) \\ & = [\mu\circ \text{flip}\circ (S\otimes S)(h_{1}\otimes h'_{1})] h_{2}h'_{2} \\ & = S(h'_{1})S(h_{1})h_{2}h'_{2} \\ & = \epsilon(h)\epsilon(h') \\ & = \epsilon(hh').} Therefore, by uniqueness, $f=f'$ and we have the first property. 
		\item Since $\eta:\B{k}\rightarrow H$ is a coalgebra map we know $S\circ \eta$ is its convolution inverse. Moreover $\eta$ is the unit of the convolution algebra $\Hom_{\B{k}}(\B{k},H)$ and so the convolution inverse of itself. Accordingly $S\circ \eta=\eta$, or $S(1)=1$. 
		\item Let $g:=\Delta\circ S:H\rightarrow H\otimes H$ and $g':=(S\otimes S)\circ \text{flip}\circ \Delta :H\rightarrow H\otimes H$. We have $g,g'$ are both left inverses of $\Delta:H\rightarrow H\otimes H$ in the convolution algebra $\Hom_{\B{k}}(H,H\otimes H)$, indeed, since $\Delta $ is an algebra map $g$ is the convolution inverse of $\Delta$ by Proposition \ref{convolution inverse in hopf algebras}, moreover for every $h\in H$ we get \eqa{ \nonumber (g'* \Delta)(h) & = g'(h_{1})\Delta(h_{2}) \\ & = (S(h_{1})_{2} \otimes S(h_{1})_{1})(h_{21} \otimes h_{22} ) \\ & = (S(h_{2})\otimes S(h_{1}))(h_{3}\otimes h_{4}) \\ & = S(h_{2})h_{3}\otimes S(h_{1})h_{4} \\ & = 1 \otimes S(h_{1})h_{2} \\ & = 1 \otimes \epsilon(h).}  Therefore $g=g'$ by uniqueness proving point 3.
		\item Finally, since $\epsilon:H\rightarrow \B{k}$ is an algebra map $\epsilon\circ S$ must be its convolution inverse. On the other hand $\epsilon$ is the unit of the convolution algebra $\Hom_{\B{k}}(H,\B{k})$ and so its the convolution inverse of itself. Therefore $\epsilon\circ S=\epsilon.$ \qed 
	\end{enumerate}\end{prop}
	\subsection{Comodule and comodule algebras}\label{comodules and comodule algebras}
 Let $H$ be a Hopf algebra. 
 \begin{defn}
 	A left $H-$module algebra is a vector space $V$ which is a left module over $H$, such that $h\triangleright 1_{A}=\epsilon(h)1_{A}$ and $h\triangleright (vw)=(h_{1}\triangleright v)(h_{2}\triangleright w)$ whenever $v,w\in V$ and $h\in H$.  Similarly, A right $H-$module algebra is a vector space $V$ which is a right module over $H$, such that $ 1_{A}\triangleleft h=1_{A} \epsilon(h)$ and $(vw)\triangleleft h=(v\triangleleft h_{1})(w\triangleleft h_{2})$, whenever $v,w\in V$ and $h\in H$.
 \end{defn}
\begin{defn}
	A left $H-$comodule is a vector space $V$ together with a map ${}_{V}\Delta:V\rightarrow H\otimes V$ such that the diagrams $$
\begin{tikzcd}
V \arrow[rr, "{}_{V}\Delta"] \arrow[dd, "{}_{V}\Delta"] &  & H\otimes V   \arrow[dd, "\id \otimes {}_{V}\Delta"] &  & V \arrow[rr, "{}_{V}\Delta"] \arrow[rrdd, "\id"] &  & H\otimes V \arrow[dd, "\epsilon\otimes \id "] \\
                                                    &  &                                                   &  &                                                &  &                                               \\
H\otimes V \arrow[rr, "\Delta\otimes \id "]         &  & H \otimes H \otimes V                             &  &                                                &  & V                                            
\end{tikzcd} 
	$$ commute, i.e. \eqa{&(\id \otimes {}_{V}\Delta)\circ {}_{V}\Delta=(\Delta \otimes \id) \circ {}_{V}\Delta,\\ & (\epsilon\otimes \id)\circ {}_{V}\Delta= \id.} The map ${}_{V}\Delta$ is called a left $H-$coaction on $V$. 
	
	A right $H-$comodule on $V$ is a vector space $V$ together with a map $\Delta_{V}:V\rightarrow V\otimes H$ such that the diagrams $$ \begin{tikzcd}
V \arrow[rr, "\Delta_{V}"] \arrow[dd, "\Delta_{V}"] &  & V\otimes H \arrow[dd, "\Delta_{V}\otimes \id"] &  & V \arrow[rr, "\Delta_{V}"] \arrow[rrdd, "\id"] &  & V\otimes H \arrow[dd, "\id\otimes \epsilon"] \\
                                                    &  &                                                &  &                                                &  &                                              \\
V\otimes H \arrow[rr, "\id\otimes \Delta"]          &  & V\otimes H \otimes H                           &  &                                                &  & V                                           
\end{tikzcd}$$ commute, i.e. \eqa{&(\id \otimes \Delta)\circ \Delta_{V} =(\Delta_{V} \otimes \id) \circ \Delta_{V},\\ & (\id\otimes \epsilon)\circ \Delta_{V}= \id.} The map $\Delta_{V}$ is called a right $H-$coaction on $V$.

\end{defn} 
\begin{defn}
Let $(V,{}_{V}\Delta)$ and $(W,{}_{W}\Delta)$ be left $H$-comodules. A morphism of left $H$-comodules, or a left $H-$colinear map, is a $\B{k}$-linear map $\phi:V\rightarrow W$ such that the diagram

$$ 
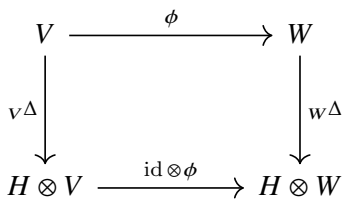
\begin{tikzcd}
V \arrow[rr,"\phi"] \arrow[dd,"{}_V\Delta"'] & & W \arrow[dd,"{}_{W}\Delta"] \\ \\
H\otimes V \arrow[rr,"\id\otimes \phi"] & &  H\otimes W
\end{tikzcd}$$ 
commutes, i.e. ${}_{W}\Delta \circ \phi = (\id\otimes \phi) \circ {}_{V}\Delta.$ Let $(V,\Delta_{V})$ and $(W,\Delta_{W})$ be right $H$-comodules. A morphism of right $H$-comodules, or a right $H-$colinear map, is a $\B{k}$-linear map $\phi:V\rightarrow W$ such that the diagram

$$ \begin{tikzcd}
V \arrow[rr,"\phi"] \arrow[dd,"\Delta_{V}"'] & &  W \arrow[dd,"\Delta_{W}"] \\ \\ 
V\otimes H \arrow[rr,"\phi\otimes \id"] & & W\otimes H
\end{tikzcd}$$ 
commutes, i.e.  i.e. $\Delta_{W} \circ \phi = (\phi\otimes \id) \circ \Delta_{V}.$
\end{defn}
\begin{notation}
	The Sweedler notation for right and left $H-$coactions on $V$ is the following: \eqa{{}_{V}\Delta v=:\sum v_{-1}\otimes v_{0}=:v_{-1}\otimes v_{0}, \quad \Delta_{V}v=:\sum v_{0}\otimes v_{1}=:v_{0}\otimes v_{1},} for all $v\in V$.  
\end{notation}

  \begin{defn}
  	A right $H-$comodule algebra is a right $H-$comodule $A$ such that $\mu:A\otimes A\rightarrow A$ and $\eta:\B{k}\rightarrow A$ are right $H-$colinear maps, namely \eqa{\Delta_{A}(ab)=\Delta_{A}(a)\Delta_{A}(b), \quad \Delta_{A}(1_{A})=1_{A}\otimes 1_{H}.} Similarly, a left $H-$comodule algebra is a left $H-$comodule $A$ such that $\mu:A\otimes A\rightarrow A$ and $\eta:\B{k}\rightarrow A$ are left $H-$colinear maps, namely \eqa{{}_{A}\Delta(ab)={}_{A}\Delta(a){}_{A}\Delta(b), \quad {}_{A}\Delta(1_{A})=1_{H}\otimes 1_{A}.} 
  \end{defn}
  \begin{defn}
  	Let $(V,{}_{V}\Delta)$ and $(W,{}_{W}\Delta)$ be left $H-$comodule algebras. A morphism of left $H-$comodule algebras is a morphism $\phi:V\rightarrow W$ of left $H-$comodules that is also an algebra map, namely $\phi(vv')=\phi(v)\phi(v')$ for every $v,v'\in V$.  Similarly for right $H-$comodule algebras.   
  \end{defn} \begin{defn}
  	Let $A$ be a left $H-$comodule algebra. We define the subalgebra ${}^{coH}A$ of left-coinvariant elements under the left $H-$coaction $\lambda_{A}:A\rightarrow H\otimes A$ as \eqa{ {}^{coH}A:=\graffe{ a\in A: \ {}_{A}\Delta(a)=1\otimes a}.} Similarly, if $A$ is a right $H-$comodule algebra, we define the subalgebra of  right-coinvariant elements under the right $H-$coaction $\Delta_{A}:A\rightarrow A\otimes H$ as \eqa{ A^{coH}:=\graffe{a\in A: \ \Delta_{A}(a)=a\otimes 1}.}  
  \end{defn}\begin{exmp} \label{adjoint H-coaction}
  	Let us define the adjoint right $H-$coaction as  \eqa{ \ad_{R}:H\rightarrow H\otimes H,\quad h\mapsto h_{2}\otimes S(h_{1})h_{3}.}  We have that $(H,\ad_{R})$ is a right $H-$comodule. Indeed \eqa{\nonumber (\id\otimes \Delta)\circ \ad_{R}(h) & = (\id\otimes \Delta) \circ (h_{2}\otimes S(h_{1})h_{3}) \\ & = h_{2}\otimes \Delta(S(h_{1})h_{3}) \\ & =  h_{2} \otimes \Delta(S(h_{1}))\Delta(h_{3}) \\ & = h_{2} \otimes (S(h_{1})_{1}\otimes S(h_{1})_{2})(h_{31}\otimes h_{32}) \\ & = h_{2} \otimes (S(h_{12})\otimes S(h_{11}))(h_{31}\otimes h_{32}) \\ & = h_{2} \otimes S(h_{12})h_{31}\otimes S(h_{11})h_{32} \\ & =  h_{3}\otimes S(h_{2})h_{4} \otimes S(h_{1})h_{5}.} 	On the other hand \eqa{\nonumber (\ad_{R}\otimes \id)\circ \ad_{R} (h) & = (\ad_{R}\otimes \id)\circ (h_{2}\otimes S(h_{1})h_{3}) \\ & = \ad_{R}(h_{2}) \otimes S(h_{1})h_{3} \\ & = h_{22}\otimes S(h_{21})h_{23}\otimes S(h_{1})h_{3} \\ & = h_{3}\otimes S(h_{2})h_{4}\otimes S(h_{1})h_{5}.}  \qed   
  \end{exmp}
  \begin{exmp}
  	Every Hopf algebra $H$ can be seen as a right $H-$comodule algebra via the coproduct $\Delta:H\rightarrow H\otimes H$. Indeed, let $h\in H$. We have that $(\Delta\otimes \id)\circ\Delta(h)= (\id\otimes \Delta)\circ \Delta(h)$, since it is the coassociativity axiom for a Hopf algebra.  Moreover $(\id\otimes \epsilon)\circ \Delta=\id $ is counitality. Finally $\Delta(ab)=\Delta(a)\Delta(b)$ and $\Delta(1_{H})=1_{H}\otimes 1_{H}$.  \qed 
  \end{exmp}
  \begin{defn}
  	Let $H$ be a Hopf algebra over a field $\B{k}$. An $H-$Hopf module is a vector space $V$ being both an $H-$module and an $H-$comodule and for which ${}_{V}\Delta:V\rightarrow H\otimes V$ is a left $H-$module map, i.e. ${}_{V}\Delta(h\cdot v)=(\Delta h).({}_{V}\Delta v)$, for all $h\in H$ and $v\in V$. The dot represents the action of $H$ on $V$.  
  \end{defn}

\begin{lemma}[\cite{beggs-majid},Lemma 2.17] Let $V$ be a left $H-$Hopf module. Then $V\simeq H \otimes ({}^{coH}V),$ where the right hand side is a Hopf module by the Hopf module structure of $H$. Conversely, every vector space defines a Hopf module by $H\otimes ({}^{coH}V)$, giving an equivalence between $H-$Hopf modules and vector spaces. 
	\qed  
\end{lemma}

  \chapter{Hopf-Galois extensions}\label{Hopf-Galois extensions}
  In this section $H$ is a Hopf algebra, $A$ is a right $H-$comodule algebra with right $H-$coaction $\Delta_{A}:A\rightarrow A\otimes H$.   
  In \ref{Hopf-Galois, cleft and trivial extensions} we give some basic results about convolution invertible morphisms, we provide the notions of trivial, cleft and Hopf-Galois extensions and we discuss some explicit realisations. 
  
  In \ref{Adjoint coactions} we introduce the adjoint coaction on the tensor product between a right $H-$comodule algebra and a Hopf algebra; we discuss some identities featuring such coaction. 
  
  In \ref{translation map} we introduce the translation map, which is the inverse of the canonical map of a Hopf-Galois extension, and discuss some properties.
  
   We introduce crossed product algebras in \ref{crossed product algebras} and discuss the proof of a theorem due to Doi-Takeuchi providing one-to-one correspondence between cross product algebras and cleft extensions, reducing to a correspondence with smash product algebras if the extension is trivial.   

Content of Sections \ref{Hopf-Galois, cleft and trivial extensions} and \ref{crossed product algebras} are essential for what we discuss in \ref{Quantum principal bundles}.
  
  A full understanding of the theory of principal bundle presented in \cite{Durdevic-1,Durdevic-2,durdevic1995quantum,urevi1998DifferentialSO} requires also sections \ref{Adjoint coactions} and \ref{translation map}.  
  
  The main works to which we refer are \cite{BRZ,BRZ2,Montgomery} 
  
   \section{Hopf-Galois, cleft and trivial extensions} \label{Hopf-Galois, cleft and trivial extensions}
Let $B:=A^{coH}$ be the subalgebra of coinvariant elements under $\Delta_{A}$. 
We start with a lemma.  
   \begin{lemma} \label{properties of cleaving map}
   	Let $j:H\rightarrow A$ be a right $H-$colinear convolution invertible map; denote by $j^{-1}:H\rightarrow A$ the convolution inverse of $j$. Then: \begin{enumerate}
   		\item  the convolution inverse satisfies $\Delta_{A}(j^{-1}(h))=j^{-1}(h_{2})\otimes S(h_{1})$ for all $h\in H$;
   		\item there is a map $A\rightarrow B$ assigning $a_{0}j^{-1}(a_{1})$ to every element $a\in A$; 
   		\item $j(1)\in A$ is invertible with inverse $j^{-1}(1)$;
   		\item there is a unital right $H-$colinear convolution invertible map $j':H\rightarrow A$ assigning $j(1)^{-1}j(h)$ to every element $h\in H$, with convolution inverse $j'^{-1}(h)=j^{-1}(h)j(1)$; 
   		\item if $j$ is an algebra morphism then $j^{-1}$ is anti-algebra morphism, that is $j^{-1}(hh')=j^{-1}(h')j^{-1}(h) $ and $j^{-1}(1)=1$ for every $h,h'\in H$.  
   	\end{enumerate}

  \proof In order: 
  \begin{enumerate}
  	\item Notice that by colinearity and properties of $\Delta_{A}$ the following equivalence holds: \eqa{\nonumber \epsilon(h)1_{A}\otimes 1_{H} = \Delta_{A}(\epsilon(h)1_{A}) & = \Delta_{A}(j(h_{1})j^{-1}(h_{2})) \\ & =\Delta_{A}(j(h_{1}))\Delta_{A}(j^{-1}(h_{2})) \\ & = (j(h_{1})_{0}\otimes j(h_{1})_{1})(j^{-1}(h_{2})_{0}\otimes j^{-1}(h_{2})_{1}) \\ & = j(h_{1})j^{-1}(h_{3})_{0}\otimes h_{2}j^{-1}(h_{3})_{1}, } for every $h\in H$.  If we define $f:=(j\otimes \id)\Delta:H\rightarrow A\otimes H$ we have \eqa{\nonumber \tonde{f*(\Delta_{A}\circ j^{-1})}(x) & = \mu \circ\tonde{ f\otimes (\Delta_{A}\circ j^{-1})}\circ \Delta(x) \\ & = \mu\circ\tonde{f(x_{1}) \otimes \Delta_{A}\circ j^{-1}(x_{2})} \\ & = \mu\circ\tonde{(j\otimes \id)\Delta(x_{1})\otimes \Delta_{A}\circ j^{-1}(x_{2})} \\ & = \mu\circ \tonde{(j(x_{11})\otimes x_{12}) \otimes \Delta_{A}(j^{-1}(x_{2}))} \\ & = \mu\circ\tonde{j(x_{1}) j^{-1}(x_{3})_{0} \otimes x_{2} j^{-1}(x_{3})_{1}} \\ & = \epsilon(x)1_{A} ,} so that $f$ is a convolution inverse of $\Delta_{A}\circ j^{-1}$ in the convolution algebra $\Hom_{\B{k}}((H,A\otimes H),*)$. Another convolution inverse of $f$ is provided by $j^{-1}(h_{2})\otimes S(h_{1})$. Therefore the claim follows from uniqueness of the inverse. 
  	\item By a direct computation we find \eqa{\nonumber  \Delta_{A}(a_{0}j^{-1}(a_{1})) & = \Delta_{A}(a_{0})\Delta_{A}(j^{-1}(a_{1})) \\  & = (a_{00}\otimes a_{01})(j^{-1}(a_{12})\otimes S(a_{11})) \\ & = a_{00}j^{-1}(a_{12}) \otimes a_{01}S(a_{11}) \\ & = a_{0}j^{-1}(a_{3})\otimes a_{1}S(a_{2}) \\ & = a_{0}j^{-1}(a_{1})\otimes 1,} and therefore we have a coinvariant element in $B:=A^{coH}$.   
  	\item As $j^{-1}$ is the convolution inverse of $j$ we find \eqa{\nonumber (j*j^{-1})(1) & = j^{-1}(1)j(1) \\ & = \epsilon(1)1 \\ & =1 \\ & =j(1)j^{-1}(1)\\ & =(j^{-1}*j)(1),} so $j^{-1}(1)$ is the inverse of $j(1)$.  
  	\item Follows from the previous point. 
  	\item By proposition \ref{convolution inverse in hopf algebras} we have $j^{-1}=j\circ S$. In particular \eqa{ \nonumber j^{-1}(hh')& =j(S(hh'))\\ & =j(S(h')S(h))\\ & = j^{-1}(h')j^{-1}(h),} for every elements in $h,h'\in H$. Accordingly $j^{-1}(1)=1.$ \qed
  \end{enumerate}
     \end{lemma}
  \bigskip
  This lemma contains several important features of the so called \textit{cleaving map}, which we define in the following. More specifically we define the notion of \textit{extension} in the setting of Hopf algebra coactions on comodule algebras. The following definition is of special interest for the development of the theory presented in \ref{Quantum principal bundles}. 
  \begin{defn}  We call $B\subseteq A$ \begin{enumerate}
  		\item a \textit{trivial extension} if there is a convolution invertible morphism $j:H\rightarrow A$ of right $H-$comodule algebras;
  		\item a \textit{cleft extension} if there is a convolution invertible morphism $j:H\rightarrow A$ of right $H$-comodules, to which we will refer to as \textit{cleaving map};
  		\item a \textit{Hopf-Galois extension} if the canonical map\footnote{ Notice the map $\chi:A\otimes_{B}A\rightarrow A\otimes H$ is well defined over $A\otimes_{B}A$, since $ab_0a'_0\otimes b_1a'_1=aba'_0\otimes a'_1$.} $$ \chi:A\otimes_{B} A \rightarrow A\otimes H, \quad a\otimes_{B}a'\mapsto aa'_{0}\otimes a'_{1} $$ is invertible.  
  	\end{enumerate}
  \end{defn}
   
  In virtue of Lemma \ref{properties of cleaving map} we assume $j,j^{-1}:H\rightarrow A$ to be unital maps. 
  
  The next result is presented to compare trivial, cleft and Hopf-Galois extensions. In particular We  show  trivial extensions to be cleft, and cleft extensions to be Hopf-Galois.  
   \begin{prop}
   	Every trivial extension is also cleft. Every cleft extension is also Hopf-Galois.
   	\proof 
   	
   	The first statement is obviously satisfied. Indeed, the difference is that in trivial extensions we require $j:H\rightarrow A$ to be an $H-$comodule algebra map, whereas in cleft extension $j:H\rightarrow A$ is an $H-$comodule map. 
   	
   	Let $B\subseteq A$ be a cleft extension and let $j:H\rightarrow A$ be the cleaving map with convolution inverse $j^{-1}:H\rightarrow A$. Define a map $\chi^{-1}:A\otimes H\rightarrow A\otimes_{B} A$ by $$\chi^{-1}(a\otimes h) := aj^{-1}(h_{1})\otimes_{B}j(h_{2}),$$ where $a\in A$ and $h\in H$.  From Lemma \ref{properties of cleaving map} we have \eqa{ \nonumber\chi(\chi^{-1}(a\otimes h)) & = \chi(aj^{-1}(h_{1})\otimes_{B}j(h_{2})) \\ & = aj^{-1}(h_{1})j(h_{2})_{0}\otimes j(h_{2})_{1} \\ & = aj^{-1}(h_{1})j(h_{2})\otimes h_{3} \\ & = a \otimes h. } On the other hand, given $a,a'\in A$, we have \eqa{ \nonumber\chi^{-1}(\chi(a\otimes_{B}a')) & = \chi^{-1}(aa'_{0}\otimes a'_{1}) \\ & = aa'_{0}j^{-1}(a'_{1})\otimes_{B}j(a'_{2}) \\ & = a\otimes_{B} a'_{0}j^{-1}(a'_{1})j(a'_{2}) \\ & = a\otimes_{B}a', } where we exploited the very definition of balanced tensor product over $B$ to move elements of $B$ along the tensor product.  \qed 
   \end{prop} 
   
   In virtue of this Proposition, and that $a'_0 j^{-1}(a'_1)\in B$ by Lemma \ref{properties of cleaving map}, we will generally refer to trivial, cleft and Hopf-Galois extensions as Hopf-Galois extension, specifying if the extension is trivial or cleft according to the special case. 

   \subsection{Examples}
    We discuss some explicit realisations of Hopf-Galois extension. 
  \begin{exmp}
  Right $H-$coinvariant elements of any bialgebra are simply scalars $\B{k}\cong H^{coH}$. Indeed $\Delta(h)=h\otimes 1$ for any $h\in H$ implies $h=\epsilon(h)1_{H}$, i.e. that $h$ must be a scalar multiple of the unit.  	On the other hand any such multiple is clearly right $H-$coinvariant. \begin{prop} Let $H$ be a bialgebra for the moment. The canonical map $\chi:H\otimes H \rightarrow H\otimes H$ sending $h\otimes h'\mapsto hh'_{1}\otimes h'_{2}$ is invertible if and only if $H$ is a Hopf algebra. 
  \proof 
If $S$ is an antipode for $H$ then $h\otimes h' \mapsto hS(h'_{1})\otimes h'_{2}$ is the inverse of the canonical map. If, on the other hand, $\chi$ is invertible, we have an antipode via $S:H\rightarrow H$, where $$ S(h):=(\id \otimes \epsilon)\circ\chi^{-1}\circ(1\otimes h).$$  The map $S$ satisfies the antipode axiom: \eqa{\nonumber h_{1}S(h_{2}) & = h_{1}(\id\otimes \epsilon)(\chi^{-1}(1\otimes h_{2}))\\ & = (\id\otimes \epsilon)(\chi^{-1}(h_{1}\otimes h_{2} ))\\ & = (\id\otimes \epsilon)(\chi^{-1}\chi(1\otimes h)) \\ & =1\otimes \epsilon(h). }  Similarly $S(h_{1})h_{2}=\epsilon(h)\otimes 1$. Therefore any Hopf algebra is in particular  a Hopf-Galois extension $\B{k}\subseteq H$.  Invertibility of the canonical map $\chi$ is equivalent to the existence of an antipode. \qed  \end{prop}
  \end{exmp}
  \begin{exmp}
  Let $F$ and $E$ be fields. We say that $E$ is a field extension of $F$ if $F\subseteq E$ is a subfield. This means $E$ can be considered as a vector space over $F$. We call the $F-$dimension of $E$ the \textit{degree of the extension}. Let now $G$ be a finite group acting by $\B{k}-$automorphisms $\Phi_{g}:E\rightarrow E$ on the field $E$; denote by $$F=E^{G}=\graffe{x\in E \st \Phi_{g}(x)=x,\  \text{for all g} \in G}$$ the set of \textit{fixed points} of $E$ under the action of $G$.  The field extension $F\subseteq E$ is called \textit{Galois} if and only if the action of $G$ is faithful (that is, $\Phi_{g}(x)=x$ for every $x\in E$ means $g=e$, the identity element of $G$) on $E$ if and only if the degree of the extension is equal to the cardinality of $G$.
  
   \begin{prop} Let $F\subseteq E$ be a Galois extension of degree $n=|G|$. Such field extension always corresponds to a Hopf-Galois extension with $F=B$, $E=A$ and $H=\Hom_{\B{k}}(\B{k}G,\B{k})$. 
   \proof 	
 
  Consider $G=\graffe{x_{1},\dots,x_{n}}$ and let $\graffe{b_{1},\dots,b_{n}}$ a basis of the quotient field $E/F$. Define a Hopf algebra $H$ as the dual of the group algebra $\B{k}G$ introduced as in Example \ref{group algebra}. $H$ has basis $\graffe{x^{1},\dots,x^{n}}$ dual to the basis of $G$, namely $x^{i}(x_{j})=\delta^{i}_{j}$.  This gives a right $H-$coaction  $\Delta_{E}:E\rightarrow E\otimes H$ explicitly presented as $\Delta_{E}(x)=\sum_{i=1}^{n}\Phi_{x_{i}}(x)\otimes x^{i}$. The canonical map reads $\chi(x\otimes_{F}x') =x\sum_{i=1}^{n}\Phi_{x_{i}}(x')\otimes x^{i}$ for all $x,x'\in E$.  
  
Faithfulness of the $G-$action follows since we are dealing with Galois field extensions. Accordingly the canonical map $\chi$ is a bijection.   \qed  \end{prop} 
  \end{exmp}
  \begin{exmp}
  	In differential geometry a fiber bundle is the datum of a smooth map $\pi:P\rightarrow M$ between smooth manifolds $P,M$ and another smooth manifold $F$, such that for every $x\in M$ there is an open neighbourhood $U\subseteq M$ of $x$ such that $\varphi_{U}:\pi^{-1}(U)\rightarrow U\cross F$ are diffeomorphic and the diagram $$ \begin{tikzcd}
  		\pi^{-1}(U) \arrow[dd ,"\pi"] \arrow[rr,"\varphi_{U}"] & &  \arrow[lldd,"pr_{1}"] U\cross F  \\ \\ U
  	\end{tikzcd},$$ where $pr_{1}:U\cross F\rightarrow U$ is to projection to the first factor, commutes. One calls $P$ the total space, $M$ the base space and $F$ the fiber. The commutative diagram above is called a \textit{local trivialisation} of the bundle. Here $\pi$ must be a surjective summersion. 
  	
 Let $G$ be a Lie group acting freely and trasitively from the right on $P$, such that the fibers are preserved, in the sense that $$\pi(p\cdot g)=\pi(p)\text{ for all} \ p\in P \text{ and }\ g\in G.$$ We call $\pi:P\rightarrow M$ a principal $G-$bundle. In particular the canonical map $$P\cross_{M}G\rightarrow P\cross_{M}P, \quad (p,g)\mapsto(p,p\cdot g)$$ into the fibered product $P\cross_{M}P=\graffe{(p,p')\in P\cross P':\pi(p)=\pi(p')}$ is a diffeomorphism. 
 
 Denote $A=\M{C}^{\infty}(P)$ and $B=\M{C}^{\infty}(M), H=\M{C}^{\infty}(G)$.  The pull-back of the right action $P\cross G\rightarrow P$ determines a right $H-$coaction $\Delta_{A}:A\rightarrow A\otimes H$ such that $B=A^{coH}$. Moreover the canonical map $\chi:A\otimes_{B} A\rightarrow A\otimes H$ is the pull-back of the canonical map of the principal bundle, and therefore it a bijection.  Therefore $B\subseteq A$ forms a Hopf-Galois extension.  \qed 
  \end{exmp}

  \subsection{Adjoint coactions} \label{Adjoint coactions}
  
  Recall the adjoint right $H-$coaction $\ad_{R}:H\rightarrow H$  of Example \ref{adjoint H-coaction}, with $\ad_{R}(h)=h_{2}\otimes S(h_{1})h_{3}$ on $H$. Consider the corresponding diagonal right $H-$coaction on the tensor product $A\otimes H$ defined as $$ \Delta_{A\otimes H}^{\ad} :A\otimes H \rightarrow A\otimes H \otimes H, \quad a\otimes h\mapsto a_{0}\otimes h_{2}\otimes a_{1}S(h_{1})h_{3}.$$  Define two auxiliary maps \eqa{\label{projected hopf galois map}\chi':A&\otimes A \longrightarrow  A \otimes H, \quad &&a\otimes a'\mapsto aa'_{0}\otimes a'_{1}, \\ \nu:A &\otimes H\longrightarrow H\otimes A \otimes H,\quad &&a\otimes h\mapsto a_{1}S(h_{1})\otimes a_{0}\otimes h_{2}.} We  discuss some properties of those maps when $B\subseteq A$ is a Hopf-Galois extension. Let $$\flip_{A,H}:A\otimes H \rightarrow H\otimes A, \quad a\otimes h \mapsto h\otimes a.$$
   \begin{lemma}\label{properties of chi}
  Let $A$ be a right $H-$comodule algebra. The following equations hold
  	\begin{enumerate}
  		\item  $(\id\otimes \chi')\circ ((\flip_{A,H}\circ\Delta_{A})\otimes \id)=\nu\circ \chi'$; 
  		\item $(\chi'\otimes \id)\circ(\id\otimes \Delta_{A})=(\id\otimes \Delta)\circ \chi'$; 
  	\item $(\chi'\otimes \id)\circ \Delta_{A\otimes A} = A^{\ad}_{A\otimes H}\circ \chi'$.
  	\end{enumerate}
  	\proof We proceed in both cases by a direct calculation. 
  	\begin{enumerate} \item We find\eqa{ \nonumber(\id \otimes \chi')\circ((\flip_{A,H}\circ \Delta_{A})\otimes \id )(a\otimes a') & = (\id \otimes \chi')\circ(\flip_{A,H}(a_{0}\otimes a_{1})\otimes a') \\ & = (\id\otimes \chi')\circ (a_{1}\otimes a_{0}\otimes a') \\ & = a_{1} \otimes \chi'(a_{0}\otimes a') \\ & = a_{1} \otimes a_{0}a'_{0}\otimes a'_{1}. } Similarly \eqa{\nonumber (\nu\circ \chi')(a\otimes a') & = \nu(\chi'(a\otimes a')) \\ & = \nu(aa'_{0}\otimes a'_{1}) \\ & = (aa'_{0})_{1}S(a'_{11})\otimes (aa'_{0})_{0}\otimes a'_{12} \\ & = a_{1}a'_{1}S(a'_{2})\otimes a_{0}a'_{0}\otimes a'_{3} \\ & = a_{1}\otimes a_{0}a'_{0}\otimes a'_{1}.  }

  \item	Next we have  \eqa{\nonumber (\chi'\otimes \id)\circ (a\otimes \Delta_{A}(a')) & = (\chi'\otimes \id)\circ (a\otimes (a'_{0}\otimes a'_{1})) \\ & = \chi'(a\otimes a'_{0})\otimes a'_{1} \\ & = aa'_{00}\otimes a'_{01} \otimes a'_{1} \\ & = aa'_{0}\otimes a'_{1}\otimes a'_{2};\\ \\ (\id\circ \Delta)\circ \chi'(a\otimes a') & = (\id \otimes \Delta)(aa'_{0}\otimes a'_{1}) \\ & = aa'_{0}\otimes a'_{1} \otimes a'_{2}.}

\item Finally \eqa{\nonumber(\chi'\otimes \id)\circ\Delta_{A\otimes A}(a\otimes a') & = (\chi'\otimes \id)(a_{0}\otimes a'_{0}\otimes a_{1}a'_{1} )\\ & = \chi'(a_{0}\otimes a'_{0})\otimes a_{1}a'_{1} \\ & = a_{0}a'_{00}\otimes a'_{01} \otimes a_{1}a'_{1} \\ & =  a_{0}a'_{0}\otimes a'_{1}\otimes a_{1}a'_{2}, \\ 
  	\\  \Delta^{\ad}_{A\otimes H} \circ \chi'(a\otimes a') & = A^{\ad}_{A\otimes H}(aa'_{0}\otimes a'_{1}) \\ & = (aa')_{0}\otimes a'_{12}\otimes (aa'_{0})_{1}S(a'_{11})a'_{13} \\ & = a_{0}a'_{0}\otimes a'_{3}\otimes a_{1}a'_{1} S(a'_{2})a'_{4} \\ & = a_{0}a'_{0}\otimes a'_{1}\otimes a_{1}a'_{2}.}\end{enumerate}
  	 \qed 
  \end{lemma}
  \subsection{The translation map}\label{translation map} 
    Let us consider the quotient map $\pi:A\otimes A \rightarrow A \otimes_{B}A'$. By the very definition of the Hopf-Galois canonical map $\chi:A\otimes_{B} A\rightarrow A\otimes H$ the composition $\chi\circ \pi:A\otimes A \rightarrow A\otimes H$ is exactly $\chi':A\otimes A\rightarrow A\otimes H$ defined in Equation \eqref{projected hopf galois map}. Consequently relations derived in the Lemma \ref{properties of chi} descend to the quotient.
    
Let $\chi$ be the Hopf-Galois canonical map of a Hopf-Galois extension then its inverse $$\chi^{-1}:A\otimes H \rightarrow A\otimes_{B}A$$ is well defined. 
\begin{defn} 
The \textit{translation map} is defined as $\kappa=\chi^{-1}|_{1_{A}\otimes H}:H\rightarrow A\otimes_{B}A$.    \end{defn}
\begin{notation} We use the Sweedler-like shorthand notation $\kappa(h)=h^{\langle 1\rangle } \otimes h^{\langle 2\rangle }$ to denote the action for the translation map, given any $h\in H$.
\end{notation} In the next proposition we exploit some important features of the translation map.  
  \begin{prop}
  	Let $B\subseteq A$ be a Hopf-Galois extension. For any $h,h' \in H$ and $a\in A$, we have: \begin{enumerate}
  		\item $h^{\langle 1\rangle }(h^{\langle 2\rangle})_{0}\otimes (h^{\langle 2\rangle})_{1}=1_{A}\otimes h$; 
  		\item $a_{0}(a_{1})^{\langle 1\rangle} \otimes_{B}(a_{1})^{\langle 2\rangle} = 1_{A} \otimes_{B}a $; 
  		\item $\kappa(hh')=h'^{\langle 1\rangle}h^{\langle 1 \rangle } \otimes_{B}h^{\langle 2\rangle}h'^{\langle 2\rangle}; $
  		\item $h^{\langle 1\rangle}h^{\langle 2 \rangle} =\epsilon(h)1_{A}$;
  		\item $h^{\langle 1\rangle} \otimes_{B}(h^{\langle 2\rangle})_{0}\otimes (h^{\langle 2\rangle})_{1}=(h_{1})^{\langle 1\rangle} \otimes_{B}(h_{1})^{\langle 2\rangle}\otimes h_{2};$ 
  	\end{enumerate}
  	\proof In order: \begin{enumerate}
  		\item  by a direct calculation \eqa{ \nonumber h^{\langle 1\rangle }(h^{\langle 2\rangle})_{0}\otimes (h^{\langle 2\rangle})_{1} &  = \chi(h^{\langle 1\rangle}\otimes h^{\langle 2\rangle}) \\ & = \chi(k(h)) \\ & = \chi(\chi^{-1}(1_{A}\otimes h)) =1_{A}\otimes h. }  
  		\item Using that $\chi:A\otimes_{B} A\rightarrow A\otimes H$ and  $\chi^{-1}:A\otimes H\rightarrow A\otimes_{B}A $ are left $A-$linear we have \eqa{\nonumber  a_{0}(a_{1})^{\langle 1\rangle} \otimes_{B}(a_{1})^{\langle 2\rangle} & = a_{0}k(a_{1})=a_{0}\chi^{-1}(1\otimes a_{1})\\ & =\chi^{-1}(a_{0}\otimes a_{1})\\ & =\chi(1_{A}\otimes_{B}a)\\ & =1_{A}\otimes_{B} a.} 
  		\item Applying two times the first property: \eqa{\nonumber \chi(h'^{\langle 1\rangle}h^{\langle 1 \rangle } \otimes_{B}h^{\langle 2\rangle}h'^{\langle 2\rangle}) & = h'^{\langle 1\rangle}h^{\langle 1\rangle} (h^{\langle 2\rangle})_{0}(h'^{\langle 2\rangle})_{0}\otimes (h^{\langle 2\rangle})_{1}(h'^{\langle 2\rangle})_{1}\\ & = h'^{\langle 1 \rangle }(h'^{\langle 2\rangle})_{0}\otimes h(h'^{\langle 2\rangle})_{1} \\ & =1_{A}\otimes hh' ,} so we find $k(hh')=\chi(h'^{\langle 1\rangle}h^{\langle 1 \rangle } \otimes_{B}h^{\langle 2\rangle}h'^{\langle 2\rangle})$.
  		\item We consider a concatenation between the relations in Equations \eqref{projected hopf galois map} descended to the quotient. \eqa{ \nonumber (\id \otimes \chi)\circ(\tau_{A,H}\circ \Delta_{A}\otimes_{B} \id)\circ \kappa(h) & = \nu\circ\chi\circ \kappa(h) \\ (\id\otimes \chi)\circ ((h^{\langle 1 \rangle})_{1} \otimes (h^{\langle 1\rangle})_{0}\otimes_{B}h^{\langle 2\rangle})& =S(h_{1})\otimes 1_{A}\otimes h_{2} \\ (h^{\langle 1\rangle})_{1}\otimes (h^{\langle 1\rangle})_{0} (h^{\langle 2\rangle})_{0}\otimes (h^{\langle 2\rangle})_{1} & = S(h_{1})\otimes 1_{A}\otimes h_{2} , } giving the claim after applying  $(\epsilon\otimes \id\otimes \epsilon)$ to both sides.  
  		\item We perform a calculation similar to the previous point. This time we concatenate the second equation of Equations \eqref{projected hopf galois map} with $(\chi^{-1}\otimes \id)$ on the left, and with $\kappa$ on the right. \eqa{ \nonumber(\chi^{-1}\otimes \id)(\chi\otimes \id)\circ(\id\otimes_{B}\Delta_{A})\kappa(h) & = (\chi^{-1}\otimes \id)(\id\otimes\Delta)\chi(\kappa(h)) \\ (\chi^{-1}\otimes \id)(\chi(h^{\langle 1\rangle}\otimes_{B}h^{\langle 2\rangle}_{0})\otimes (h^{\langle 2\rangle})_{1}) & = (\chi^{-1}\otimes \id)(1_{A}\otimes h_{1}\otimes h_{2}) \\ h^{\langle 1\rangle }\otimes_{B} (h^{\langle 2\rangle})_{0} \otimes (h^{\langle 2\rangle})_{1} & = \chi^{-1}(1_{A}\otimes h_{1})\otimes h_{2} \\ & = \kappa(h_{1})\otimes h_{2} \\ & = (h^{\langle 1\rangle})_{1}\otimes_{B}(h^{\langle 2\rangle})_{1}\otimes h_{2}.}  \qed 
  	\end{enumerate}
  \end{prop}

  \section{Crossed product algebras}\label{crossed product algebras}
  In this section we define \textit{crossed product algebras} and \textit{smashed product algebras}. The Doi-Takeuchi correspondence between cleft extension and crossed product algebras (trivial extensions and smash product algebras) of \cite{Doi-Takeuchi} is proven. 
  
   Let $B$ an algebra.
   \begin{defn} We say that $H$ measures $B$ if there is a $\B{k}-$linear map $$\cdot:H\otimes B\rightarrow B, \quad  h\otimes b\mapsto h\cdot b,$$ such that \eqa{\nonumber  h\cdot 1_{B} & =\epsilon(h)1_{B}, \\  h\cdot (bb')& =(h_{1}\cdot b)(h_{2}\cdot b'). } \end{defn} We stress that $\cdot:H\otimes B \rightarrow B$ is not assumed to be an action.   
   \begin{defn}
   	Assume $H$ measures $B$. Consider a map $\sigma:H\otimes H\rightarrow B$.  \begin{enumerate}
   		\item  We call $\sigma:H\otimes H \rightarrow B$ a 2-\textit{cocyle with values in} $B$ if it is a convolution invertible morphism such that \eqa{& \label{property 1} \sigma(h\otimes 1) =\epsilon(h)1_{B}=\sigma(1\otimes h), \\   & (h_{1}\cdot \sigma(h'_{1}\otimes h''_{1}))\sigma(h_{2}\otimes h'_{2}h''_{2}) =\sigma(h_{1}\otimes h'_{1})\sigma(h_{2}h'_{2}\otimes h''),} for all $h,h',h''\in H$ and $b\in B$.
   		\item We call $B$ a $\sigma-$\textit{twisted left} $H-$\textit{module} if there is a 2-cocycle $\sigma:H\otimes H\rightarrow B$ with values in $B$ such that \eqa{ \label{property 2} & 1\cdot b=b, \\ & h\cdot (h'\cdot b) =\sigma(h_{1}\otimes h'_{1})((h_{2}h'_{2})\cdot b)\sigma^{-1}(h_{3}\otimes h'_{3}),}  for all $h,h'\in H$ and $b\in B$.
   	\end{enumerate}
   \end{defn}
    \begin{lemma}\label{crossed product algebras are associative unital}
 Given $B$ a $\sigma-$twisted left $H-$module, we define a $\B{k}-$linear map  \eqa{\mu_{\sharp_{\sigma}}:(B\otimes H)\otimes (B\otimes H)\rightarrow B\otimes H,\quad (b\otimes h)\otimes (b'\otimes h')\mapsto b(h_{1}\cdot b')\sigma(h_{2}\otimes h'_{1})\otimes h_{3}h'_{2},} for every $h,h'\in H$ and $b\in B$,  providing an associative unital product $\mu_{\sharp_{\sigma}}$ on $B\otimes H$, with unit $1_{B}\otimes 1_{H}$. 
   	 \proof First notice the $\B{k}-$linear map $\mu_{\sharp_{\sigma}}$ is well defined. To prove associativity we proceed considering
   	 \eqa{\nonumber \mu_{\sharp_{\sigma}}(\mu_{\sharp_{\sigma}}((b\otimes h)& \otimes (b'\otimes h'))\otimes (b''\otimes h''))  = \mu_{\sharp_{\sigma}}((b(h_{1}\cdot b')\sigma(h_{2}\otimes h'_{1})\otimes h_{3}h'_{2} \otimes b'' \otimes h'')\\ & = b(h_{1}\cdot b')\sigma(h_{2}\otimes h'_{1})(h_{3}h'_{2})_{1}\cdot b''_{1}\cdot \sigma((h_{3}h'_{2})_{2}\otimes h''_{1})\otimes (h_{3}h'_{2})_{3}h''_{2} \\ \small{\textbf{relabeling}}\quad  & = b(h_{1}\cdot b')\sigma(h_{2}\otimes h'_{1})(h_{31}h'_{21}\cdot b''_{1})\cdot \sigma(h_{32}h'_{22}\otimes h''_{1})\otimes h_{33}h'_{23}h''_{2} \\ \small{\textbf{Equation \eqref{property 2}} }\quad & = b(h_{1}\cdot b')\sigma(h_{2}\otimes h'_{1})\sigma^{-1}(h_{3}\otimes h'_{2})(h_{4}\cdot h'_{3}\cdot b'')\sigma(h_{5}\otimes h'_{4})\sigma(h_{6}h'_{5}\otimes h''_{1})\otimes h_{7}h'_{6}h''_{2} \\ \small{\textbf{simplifying}}\quad & = b(h_{1}\cdot b')(h_{2}\cdot h'_{1}\cdot b'')\sigma(h_{3}\otimes h'_{2})\sigma(h_{4}h'_{3}\otimes h''_{1})\otimes h_{5}h'_{4}h''_{2} \\ \small{\textbf{measurability}} \quad & = b(h_{1}\cdot(b'(h'_{1}\cdot b'')))\sigma(h_{2}\otimes h'_{2})\sigma(h_{3}h'_{3}\otimes h''_{1})\otimes h_{4}h'_{4}h''_{2} \\ \small{\textbf{Equation \eqref{property 1}}} \quad & = b(h_{1}\cdot (b'(h'_{1}\cdot b'')\sigma (h_{2}\otimes h''_{1})))\sigma(h_{2}\otimes h'_{3}h''_{2})\otimes h_{3}h'_{4}h''_{3} \\ \small{\textbf{measurability}}\quad & =  b(h_{1}\cdot (b'(h'_{1}\cdot b'')\sigma(h_{2}\otimes h''_{1})))\sigma(h_{2}\otimes h'_{3}h''_{2})\otimes h_{3}h'_{4}h''_{3} \\ & =  \mu_{\sharp_{\sigma}}((b\otimes h)\otimes (b'(h'_{1}\cdot b'')\cdot \sigma (h'_{2}\otimes h''_{1})\otimes h'_{3}h''_{2}) \\ & = \mu_{\sharp_{\sigma}}(b\otimes h \otimes \mu_{\sharp_{\sigma}}((b'\otimes h')\otimes (b''\otimes h''))),  } 
   	 where we used property of $B$ being a $\sigma-$twisted left $H$-module, so that $$ h(h'\cdot b)  = \sigma^{-1}(h_{1}\otimes h'_{1})((h_{2}h'_{2}\cdot b)(\sigma(h_{3}\otimes h'_{3})), $$ and the property of 2-cocycle with values in $B$. The last part of the proof is to show that $1_{B}\otimes 1_{H}$ is a unit, but this follows immediately. \qed     \end{lemma} 
   	 
   	$B\sharp_{\sigma}H$ is called a \textit{cross product algebra}. We define on $B\sharp_{\sigma}H$ the structure of a right $H-$comodule algebra with respect to the right $H-$coaction $\Delta_{\sharp_{\sigma}}:=\id_{B}\otimes\Delta:B\sharp_{\sigma}H\rightarrow(B\sharp_{\sigma}H)\otimes H.$ Accordingly, the subalgebra of coinvariant elements in $B\sharp_{\sigma}H$ is  $$(B\sharp_{\sigma}H)^{coH}=(B\otimes 1) \cong B.$$
   \begin{rmk}
   	 Every left $H-$module algebra $B$ is in particular a $\sigma-$twisted left $H-$module with respect to the trivial $2-$cocycle $$ \sigma:H\otimes H \rightarrow B, \quad h\otimes h'\mapsto \epsilon(hh')1_{B}.$$ The corresponding product $$ \mu_{\sharp}:(B\otimes H)\otimes (B\otimes H) \rightarrow B\otimes H, \quad (b\otimes h)\otimes (b'\otimes h')\mapsto b(h_{1}\cdot b')\otimes h_{2}h'$$ makes $B\sharp H:=(B\otimes H,\mu_{\sharp})$ an associative unital algebra, the \textit{smash product algebra}. 
   \end{rmk}
   The next theorem was proven in \cite{Doi-Takeuchi}. It shows that crossed product algebras are in $1:1$ correspondence with cleft extensions. Moreover, trivial extensions are in $1:1$ correspondence with smash product algebras. 
   
   \begin{thm}\label{Doi-Takeuchi}
   	Any crossed product algebra $B\sharp_{\sigma}H$ is a cleft extension $B\subseteq B\sharp_{\sigma}H$ with cleaving map $j:H\rightarrow B\sharp_{\sigma}H$ assigning $h\mapsto 1\sharp_{\sigma}h$. Conversely, giving a cleft extension $B\subseteq A$ with cleaving map $j:H\rightarrow A$ we define a $\sigma-$twisted left $H-$module action $$\cdot :H\otimes B \rightarrow B,\quad \quad h\otimes b \mapsto h\cdot b:=j(h_{1})bj^{-1}(h_{2})$$ on B, and a 2-cocycle $$\sigma:H\otimes H \rightarrow B,\quad \quad h\otimes h'\mapsto\sigma(h\otimes h'):=j(h_{1})j(h'_{1})j^{-1}(h_{2}h'_{2})$$ with values in $B$.  Then $A\cong B\sharp_{\sigma}H$ are isomorphic as right $H-$comodule algebras. A cleft extension is a trivial extension if and only if the corresponding crossed product algebra is a smash product algebra.  
   	
   	\proof 
   	Given a crossed product algebra $B\sharp_{\sigma}H$ we consider a right $H-$colinear map \eqa{ j:H\rightarrow B\sharp_{\sigma}H,\quad h\in H\mapsto 1\sharp_{\sigma}h.} This map is convolution invertible with inverse \eqa{j^{-1}:H\rightarrow B\sharp_{\sigma}H, \quad h\mapsto \sigma^{-1}(S(h_{2})\otimes h_{3})\sharp_{\sigma} S(h_{1}).}  Indeed \eqa{\nonumber j^{-1}(h_{1})j(h_{2}) & = (\sigma^{-1}(S(h_{2})\otimes h_{3}) \sharp_{\sigma} S(h_{1}))(1\sharp_{\sigma} h_{4}) \\ & = \mu_{\sharp_{\sigma}}\quadre{(\sigma^{-1}(S(h_{2})\otimes h_{3}) \otimes S(h_{1}))\otimes (1\otimes h_{4})} \\ & = \sigma^{-1}(S(h_{2})\otimes h_{3})(S(h_{1})_{1}\cdot 1)\sigma(S(h_{1})_{2}\otimes h_{4})\otimes S(h_{1})_{3}h_{5} \\ & = \sigma^{-1}(S(h_{3})\otimes h_{4})\sigma(S(h_{2})\otimes h_{5})\otimes S(h_{1})h_{6} \\ & = \sigma^{-1}(S(h_{2})_{2}\otimes h_{3})\sigma(S(h_{2})_{2}\otimes h_{4})\otimes S(h_{1})h_{5} \\ & = 1\sharp_{\sigma} S(h_{1})h_{2} \\ &  = \epsilon(h)1\sharp_{\sigma} 1, } and similarly evaluating in the other order, proving $B\subseteq B\sharp_{\sigma}H$ is a cleft extension.  
   	
   	\smallskip
   	
   	For the other implication, given a cleft extension $B\subseteq A$ with cleaving map $j:H\rightarrow A$ we prove $\sigma:H\otimes H \rightarrow B$ and $\cdot:H\otimes B \rightarrow B$ provide a $\sigma-$twisted left $H-$module structure on $B$. We have that $\cdot:H\otimes B \rightarrow B$ is a $H-$measure on $B$, indeed \eqa{\nonumber(h\cdot 1_{B}) & =j(h_{1})j^{-1}(h_{2})\\ & =\epsilon(h)1_{B};}  moreover \eqa{\nonumber(h_{1}\cdot b)(h_{2}\cdot b') & = j(h_{1})bj^{-1}(h_{2})j(h_{3})b'j(h_{4})\\ & =j(h_{1})bb'j(h_{2}) \\ &  = h\cdot (bb'). }Moreover, $\sigma:H\otimes H\rightarrow B$ is a \textbf{2-cocycle} with values in $B$:\begin{enumerate} \item The image of $\sigma$ is contained in $B$, indeed \eqa{\nonumber \Delta_{A}(j(h_{1})j(h_{2})j^{-1}(h_{2}h'_{2})) & = \Delta_{A}(j(h_{1}))\Delta_{A}(j(h'_{1}))\Delta_{A}(j^{-1}(h_{2}h'_{2})) \\ & = (j(h_{1})_{0}\otimes j(h_{1})_{1})(j(h'_{1})_{0}\otimes j(h'_{1})_{1})(j^{-1}(h_{2}h'_{2})_{2} \otimes S(h_{2}h'_{2})_{1}) \\ & = j(h_{1})j(h'_{1})j(h_{4}h'_{4}) \otimes h_{2}h'_{2}S(h_{3}h'_{3}) \\ & = j(h_{1})j(h'_{1})j(h_{2}h'_{2})\otimes 1,} where we used Lemma \ref{Doi-Takeuchi}. The convolution inverse of $\sigma$ is $$\sigma^{-1}(h\otimes h')=j(h'_{1}h_{1})j^{-1}(h_{2})j^{-1}(h'_{2}).$$  \item  $\sigma(1\otimes h)=\sigma(h\otimes 1)=\epsilon(h)1$. 
   	\item Finally \eqa{\nonumber (h_{1} \cdot \sigma(h'_{1}\otimes h''_{1}))& (\sigma(h_{2}\otimes h'_{2}h''_{2})) \\ &  = j(h_{11})\sigma(h'_{1}\otimes h''_{1})j^{-1}(h_{12})j(h_{21})j(h'_{21}h''_{21})j^{-1}(h_{22}h'_{22}h''_{22}) \\ & = j(h_{11})j(h'_{11})j(h''_{11})j^{-1}(h'_{12}h''_{12})j^{-1}(h_{12})j(h_{21})j(h'_{21}h''_{21})j^{-1}(h_{22}h'_{22}h''_{22}) \\ & = j(h_{1})j(h'_{1})j(h''_{1})j^{-1}(h'_{2}h''_{2})j^{-1}(h_{2})j(h_{3})j(h'_{3}h''_{3})j^{-1}(h_{4}h'_{4}h''_{4}) \\ & = j(h_{1})j(h'_{1})j(h''_{1})j^{-1}(h_{2}h'_{2}h''_{2})\\ & = j(h_{1})j(h'_{1})j^{-1}(h_{2}h'_{2})j(h_{3}h'_{3})j(h''_{1})j^{-1}(h_{4}h'_{4}h''_{2})\\  & = \sigma(h_{1}\otimes h'_{1})\sigma(h_{2}h'_{2}\otimes h'').} \end{enumerate}

   	The measure $\cdot:H\otimes B\rightarrow B $ is a $\sigma-$twisted left $H-$module action, since \eqa{\nonumber 1\cdot b& =j(1)bj^{-1}(1)=b; \\ \sigma(h_{1}\otimes h'_{1})((h_{2}h'_{2})\cdot b)\sigma^{-1}(h_{3}\otimes h'_{3}) & = \sigma(h_{1}\otimes h'_{1})(j(h_{21}h'_{21})bj^{-1}(h_{22}h'_{22}))\sigma^{-1}(\sigma_{3}\otimes \sigma'_{3}) \\ & = j(h_{1})j(h'_{1})j^{-1}(h_{2}h'_{2})j(h_{3}h'_{3})bj^{-1}(h_{4}h'_{4})j(h_{5}h'_{5})j^{-1}(h'_{6})j(h_{6}) \\ & = j(h_{1})j(h'_{1})bj^{-1}(h'_{2})j^{-1}(h_{2}) \\ & = h\cdot (h'\cdot b),  } for every $h,h'\in H$ and $b\in B$.  	Consequently we can construct the crossed product algebra $B\sharp_{\sigma}H$.  
   	
   	The only thing left to show it that $B\sharp_{\sigma} H$ is isomorphic to $A$ as a right $H-$comodule algebra, the explicit isomorphism is provided by $$\theta:A\rightarrow B \sharp_{\sigma}H,\quad a\mapsto a_{0}j^{-1}(a_{1})\otimes a_{2}.$$ This map is a right $H-$comodule morphism and is well defined according to Lemma \ref{properties of cleaving map}. It is left to prove it is also an algebra morphism. We have \eqa{\nonumber \theta(a)\theta(a')& = (a_{0}j^{-1}(a_{1})\sharp_{\sigma}a_{2})(a'_{0}j^{-1}(a_{1})\sharp_{\sigma}a'_{2}) \\ & = \mu_{\sharp_{\sigma}}(a_{0}(j^{-1}(a_{1}))\otimes a_{2})\otimes (a'_{0}j^{-1}(a_{1})\otimes a'_{2}) \\ & = a_{0}j^{-1}(a_{1})(a_{21}\cdot (a'_{0}(j^{-1}(a'_{1})))\sigma(a_{22}\otimes a'_{21})\otimes a_{23}a'_{22} \\  & = a_{0}j^{-1}(a_{1})j(a_{21})a'_{0}j^{-1}(a'_{1})j^{-1}(a_{22})j(a_{31})j(a'_{21})j^{-1}(a_{32}a'_{22})\otimes a_{33}a'_{22}\\ & =  a_{0}j^{-1}(a_{1})j(a_{2})a'_{0}j^{-1}(a'_{1})j^{-1}(a_{3})j(a_{4})j(a'_{2})j^{-1}(a_{5}a'_{3})\otimes a_{6}a'_{4}\\ & = a_{0}a'_{0}j^{-1}(a_{1}a'_{1}) \otimes a_{2}a'_{2} \\ & = \theta(aa'). } 
   	
   	In the context of a trivial extension, then $j:H\rightarrow A$ is an algebra morphism. The induced $2-$cocycle becomes trivial \eqa{ \nonumber\sigma(h\otimes h')& =j(h_{1})j(h'_{1})j^{-1}(h_{2}h'_{2})\\ &  = j(h_{1})j(h'_{1})j^{-1}(h'_{2})j^{-1}(h_{2})\\ & =\epsilon (hh')1.}  Thus, $B\sharp H$ is a smash product algebra.
   	
   	 Conversely if $(B\sharp H,\mu_{\sharp})$ is a smash product algebra, the induced cleaving map is an algebra map.  \qed     
   \end{thm}
   \chapter{Differential calculus over algebras}\label{Differential calculus over algebras}
   In this chapter we explore the topic of differentials over algebras. In the context of (classical) differential geometry we usually considers a $n-$dimensional manifold $M$  with a differentiable structure, that is, every open set of $M$ is identified with a open set of $\B{R}^{n}$ in "some" smooth way (smooth atlas). These \textit{local patches} fit in such way that we can talk of a global structure on $M$. Tangent and cotangent bundles to $M$ are defined as dual structures. They are respectively obtained by gluing the spaces of tangent vectors and spaces of covectors for every point in the \textit{manifold} $M$. Sections of the tangent bundle act as derivations of the (commutative) algebra of $C^{\infty}$ functions on $M$. Dually, we have a map $\dd:C^{\infty}(M)\rightarrow \Gamma(T^{*}M)$ turning smooth functions into covectors. Differential forms are obtained by subsequent applications of the exterior derivative $\dd:\Gamma\tonde{\Lambda^{\bullet}(T^{*}M)} \rightarrow \Gamma\tonde{\Lambda^{\bullet}(T^{*}M)}$ on sections of the exterior powers of the cotangent bundle.
   
   This construction can be generalised to the setting of \textit{noncommutative} algebras. In this picture no actual topological space is required; we actually consider differential geometry over algebras in this sense.  
   
   In \ref{first order calculi} we introduce first order differential calculi over  $\B{k}-$algebras, we provide some examples and prove a theorem by Woronowicz stating that every first order differential calculus on an algebra is induced as a quotient of the so called \textit{universal} differential calculus. Then we introduce the notion of \textit{covariant} calculus. 
    
     In \ref{Woronowicz classification section} we prove another important classification result due to Woronowicz, providing a $1$ to $1$ correspondence between covariant calculi and certain ideals defined on the underlying algebra. 
   
   In \ref{higher order calculi} we generalise the first order theory to higher order differential forms. In \ref{higher order covariant calculi} we discuss higher order covariant differential calculi and provide an explicit formula for the $H-$coaction on higher order differential forms.

  In \ref{maximal prolongation} we introduce the maximal prolongation, an higher order differential calculus defined just by the datum of a first order calculus over a $\B{k}-$algebra. This is of particular interest for the theory developed in \ref{Quantum principal bundles}. 
   
    The main references we follow here are \cite{beggs-majid, woronowicz,SCHAUENBURG1996239}  
   \section{First order differential calculi}\label{first order calculi}
  \begin{defn} A first order differential calculus $(\Gamma,\dd)$ over a $\B{k}-$algebra $A$ is the datum of\begin{enumerate}
   	\item an $A-$bimodule $\Gamma$; 
   	\item a linear map $\dd:A\rightarrow \Gamma$ satisfying the Leibniz rule $ \dd(ab)=(\dd a)b+ a \dd b$ for every $a,b\in A$;
   	\item a surjectivity condition $\Gamma=A\dd A$, i.e. $\Gamma=\span\{a\dd b :  a,b\in A\}.$ 
   \end{enumerate}\end{defn}  In classical differential geometry $A=C^{\infty}(M)$, and left and right module structure on $\Gamma$ always coincide, meaning that for every $a,b\in A$ we have $a\dd b = \dd b a$. This is not true in the general noncommutative setting. 
   \begin{rmk}[Surjectivity from right] The surjectivity condition $\Gamma=A\dd A$ of a differential calculus $(\Gamma,\dd)$ on $A$ is equivalent to surjectivity from the right $\Gamma=\dd(A)A$. In fact, for every $\omega\in \Gamma$ there are $a^{i},b^{i}\in A$, with $1\leq i \leq n$ such that $\omega=a^{i}\dd b^{i}$. We define $e^{j},f^{j}\in A$ for $1\leq j \leq 2n $ by  $e^{j}=a^{j}b^{j}$, with $f^{j}=1$ for $1\leq j \leq n$ and $e^{j}=a^{j-n}$, $f^{j}=b^{j-n}$ for $n< j \leq 2n$. Then $$ \omega= a^{i} \dd b^{i} =  \dd(a^{i}b^{i})-\dd(a^{i})b^{i} = \dd (e^{j})f^{j},$$ by the Leibniz rule. Similarly the other implication is proven.
\end{rmk}
   \begin{rmk}[Connected calculus]
	From the Leibniz rule it immediately follows that every first order differential calculus satisfies $ \dd(1)=\dd(1\cdot 1)=2\dd(1)$, from which we find $\B{k}\cdot 1 \subseteq \ker \dd$. If the last is an equality we say the calculus is \textit{connected}.
\end{rmk}
   
We provide a few explicit realisations of first order differential calculi.  
   \begin{exmp}[q-differential calculus on the circle]\label{circle}
Let $q\in \B{C}$ not a root of unity. Consider the algebra $A=\B{S}_{q}^{1}=\B{C}[t,t^{-1}]$ of rational polynomials in one variable $t$. Define $\Gamma$ as the free left $A-$module generated by $\dd t$.  Every element of $\Gamma$ will be of the form $f(t)\cdot \dd t$. The free left $A-$module $\Gamma$ becomes an $A-$bimodule via $$f(t)\cdot \dd t \cdot g(t) = f(t)g(qt)\cdot \dd t $$ for $f,g\in A$. We define a $\B{C}-$linear map $\dd:A\rightarrow \Gamma$, the exterior derivative, by 
$$ \dd(f(t))= \frac{f(qt)-f(t)}{t(q-1)}\cdot \dd t, $$ for every $f\in A$. The exterior derivative satisfyies the Leibniz rule.  
\eqa{\nonumber \dd(f(t)) \cdot g(t) + f(t) \cdot \dd(g(t)) & = \frac{f(qt)-f(t)}{t(q-1)}\cdot \dd t \cdot g(t) + f(t) \frac{g(qt)-g(t)}{t(q-1)} \cdot \dd t \\ & = \frac{ f(qt)g(qt)-f(t)g(qt)+f(t)g(qt)-f(t)g(t)}{t(q-1)} \cdot \dd t \\ & = \frac{f(qt)g(qt)-f(t)g(t)}{t(q-1)}\cdot \dd t \\ & = \dd (f\cdot g)(t), } for all $f,g\in A$. 

  By construction $\Gamma =A\cdot \dd t = A \dd A$. Thus $(\Gamma,\dd)$ is a first order differential calculus over $A$.  
  \end{exmp}
   \begin{exmp}[Quotient differential calculus] \label{quotient dc}
  Given a surjective algebra map $\pi:A\rightarrow H$, and a first order differential calculus $(\Gamma,\dd)$ on $A$, we define $I:=\ker \pi$ and $N_{I}:=I\dd A + A \dd I$. Clearly $I\subseteq A$ is an algebra ideal such that $H\cong A/I$. 
  
  Using the Leibniz rule one shows that $N_{I}\subseteq \Gamma$ is an $A-$bimodule. The induced calculus $(\Gamma_{H},\dd_{H})$ is the quotient calculus on $H$. Namely, we set $\Gamma_{H}:=\Gamma/N_{I}$ and $\dd_{H}=\pi_{\Gamma}\circ \dd$, where $\pi_{\Gamma}:\Gamma\rightarrow\Gamma_{H}$ is the quotient map. 
  \end{exmp} 
  \begin{defn}
  	Let $(\Gamma,\dd)$ on $A$ and $(\Gamma',\dd')$ on $A'$ be first order differential calculi. A morphism of differential calculi is a tuple $(\Phi,\phi)$, where $\Phi:\Gamma\rightarrow \Gamma'$ is a $\B{k}-$linear map, and $\phi:A\rightarrow A'$ is an algebra morphism, such that \eqa{\Phi(a\cdot \omega\cdot b)=\phi(a)\cdot' \Phi\cdot'\phi(b),} for all $a,b\in A$ and $\omega\in \Gamma$, and such that the diagram $$ \begin{tikzcd}
\Gamma \arrow[rr, "\Phi"]              &  & \Gamma'               \\
                                       &  &                       \\
A \arrow[rr, "\phi"] \arrow[uu, "\dd"] &  & A' \arrow[uu, "\dd'",swap]
\end{tikzcd}$$ commutes. In this case we say $\phi$ is \textit{differentiable}. We write  $\dd \phi=\Phi$.
  \end{defn}
  \begin{prop}[Universal differential calculus]\label{UDC}
  	For every algebra $A$ there is a first order differential calculus $(\Gamma_{u},\dd_{u})$ defined by $$ \Gamma_{u}:=\ker \mu_{A}= \graffe{\sum_{i}a^{i}\otimes b^{i}\in A\otimes A : \sum_{i}a^{i}b^{i}=0},$$ with differential $\dd_{u}(a):=1\otimes a - a\otimes 1$ for all $a\in A$.
  	\proof 
  	  The $A-$bimodule structure on $\Gamma_{u}$ is induced from $A\otimes A$, i.e. the multiplication on the first tensor factor from the left and the second tensor factor from the right, respectively. 
  	  
  	  The map $\dd_{u}$ maps into $\Gamma_{u}$ and satisfies the Leibniz rule \eqa{ \nonumber \dd_{u}(a)b+a\dd_{u}(b)& =1\otimes ab - a\otimes b +a\otimes b -ab\otimes 1 \\ & =\dd_{u}(ab),} for all $a,b\in A$. 
  	  
  Let $a^{i}\otimes b^{i} \in \Gamma_{u}:=\ker \mu_{A}$ be arbitrary. Then \eqa{\nonumber \sum_{i}a^{i}\dd_{u}(b^{i}) & = \sum_{i}a^{i}\otimes b^{i} -\sum_{i}a^{i}b^{i}\otimes 1 \\ & =  \sum_{i}a^{i}\otimes b^{i},} which proves surjectivity.  Therefore $(\Gamma_{u},\dd_{u})$ is a first order differential calculus over $A$.

  \qed \end{prop}
The next theorem was first proven by Woronowicz in \cite{woronowicz}. It states that every first order differential calculus is induced as quotient of the \textit{universal} differential calculus in Proposition \ref{UDC}. 
\begin{thm}[Woronowicz] \label{Woronowicz classification}Let $A$ be a unital $\B{k}-$algebra. Every first order differential calculus $(\Gamma,\dd)$ over $A$ can be obtained as a quotient of the \textit{universal differential calculus} $(\Gamma_{u},\dd_{u})$, meaning there exists an $A-$subbimodule $N\subseteq \Gamma_{u}$ such that $\Gamma\cong \Gamma_{u}/N$ as $A-$bimodules, and $\dd=\pi \circ \dd_{u}$, where $\pi:\Gamma_{u}\rightarrow \Gamma_{u}/N\cong \Gamma$ is the quotient map. By construction

$$\begin{tikzcd}
\Gamma_{u} \arrow[rr, "\pi"] &                                             & \Gamma \\
                             &                                             &        \\
                             & A \arrow[luu, "\dd_{u}"] \arrow[ruu, "\dd",swap] &       
\end{tikzcd}. $$ \end{thm}The denomination "universal" is used since for every calculus $\Gamma$ on $A$ there is a unique surjective morphism $\pi$ such that the diagram shown commutes. 
\proof Let us consider any $\B{k}-$algebra $A$ and the first order differential calculus $(\Gamma_{u},\dd_{u})$ of Proposition \ref{UDC}.  Let $(\Gamma,\dd)$ be any other differential calculus on $A$. Consider the $A-$bilinear map $\pi :\Gamma_{u}\rightarrow \Gamma$ assigning $a^{i}\otimes b^{i}\mapsto a^{i}\dd b^{i}$. 

The map $\pi$ is surjective. Indeed, for any $\omega=a^{i}\dd b^{i}\in \Gamma$ we have \eqa{ \nonumber \pi (\sum_{i}a^{i}\otimes b^{i}-\sum_{i}a^{i}b^{i}\otimes 1) & = \sum_{i}a^{i}\dd b^{i} - \sum_{i}a^{i}b^{i} \dd (1) \\ & = \sum_{i}a^{i}\dd b^{i}.} Considering the quotient of $\Gamma_{u}$ by the subbimodule $N=\ker\pi $ we provide the isomorphism $\Gamma\cong \Gamma_{u}/ N.$ 

Moreover  \eqa{\nonumber \pi(\dd_{u}(a)) & = \pi (1\otimes a -a\otimes 1)\\ & =1\cdot \dd (a) \\ & = \dd (a),} for every $a\in A$.  \qed 

	\subsection{The classification of first order $H-$covariant calculi}\label{Woronowicz classification section}
Here we introduce the notion of $H-$covariant differential calculi over $A$, where $H$ is a Hopf algebra and $A$ is an $H-$comodule algebra. 
	\begin{defn}
 A first order differential calculus $(\Gamma),\dd)$ on a right $H-$comodule algebra $(A,\Delta_{A})$ is called \textit{right} $H-$\textit{covariant} if $\Gamma$ is a right $H-$covariant $A-$bimodule with coaction $$\Delta_{\Gamma}:\Gamma\rightarrow \Gamma\otimes H$$ such that the differential $\dd:A\rightarrow \Gamma$ is right $H-$colinear. 
 
 Similarly, it is called \textit{left} $H-$\textit{covariant} if $\Gamma$ is a left $H-$covariant $A-$bimodule with coaction $${}_{\Gamma}\Delta:\Gamma\rightarrow H\otimes \Gamma,$$ such that the differential $\dd:A\rightarrow \Gamma$ is left $H-$colinear. 
 
 It is called $H-$\textit{bicovariant} if its both left and right $H-$covariant. 
	\end{defn}

Let us fix $(\Gamma,\dd)$  a left covariant first order differential calculus over a Hopf algebra $H$. We denote the module of left coinvariant forms by \eqa{\Lambda^{1}:=\{\omega\in \Gamma:  {}_{\Gamma}\Delta(\omega)=1\otimes\omega\}.} 
\begin{defn} We define the \textit{quantum Maurer-Cartan form} as the canonical $\B{k}$-linear map from the kernel $H^{+}=\ker\epsilon$ of the counit $\epsilon:H\rightarrow \B{k}$  to the coinvariant forms as  \eqa{\varpi:H^{+}\rightarrow\Lambda^{1},\quad  h\mapsto S(h_{1})\dd h_{2}.} \end{defn} 
\begin{lemma} \label{maurer cartan}
	The quantum Maurer-Cartan form $\varpi:H^{+}\rightarrow\Lambda^{1}$ is a surjective morphism and moreover $I=\ker\varpi\subseteq H$ is a right ideal. \end{lemma} We will use this Lemma in the proof of the classification theorem due to Woronowicz, originally proven in \cite{woronowicz}. This result provides a correspondence between right ideals $I\subseteq H$ that are in $H^{+}$ and left covariant first order differential calculi $\Gamma$ over $H$. In particular explicit forms of the bimodule structure, the differential and left $H-$coaction are provided. 
	 \proof Let  $\triangleright:H\otimes \Gamma \rightarrow \Gamma$ be the left module action of $H$ on $\Gamma$. Considering a left-coinvariant form $\omega=a^{i}\dd b^{i}$ we find \eqa{\nonumber \triangleright\circ (S\otimes\id)\circ{}_{\Gamma}\Delta(\omega) &=\triangleright\circ (S\otimes \id)(1\otimes \omega) \\ & =\omega} by left $H-$coinvariance. 
	 
	 Moreover\eqa{ \nonumber \triangleright\circ (S\otimes\id)\circ{}_{\Gamma}\Delta(\omega) & = \triangleright\circ (S\otimes \id)(a_{0}^{i}b_{0}^{i}\otimes a_{1}^{i}\dd b_{1}^{i}) \\ & = \triangleright\circ (S(a_{0}^{i}b_{0}^{i})\otimes a_{1}^{i}\dd b_{1}^{i})\\ & = S(b_{0}^{i})S(a_{0}^{i})a_{1}^{i}\dd b_{1}^{i} \\ & = \epsilon(a^{i})S(b_{0}^{i})\dd b_{1}^{i} \\ & = \varpi(\epsilon(a^{i})b^{i}-\epsilon(a^{i}b^{i})1_{H}), } where we used that $\varpi(1_{H})=0$. Notice that every expression is well defined as $\epsilon(a^{i})b^{i}-\epsilon(a^{i}b^{i})1_{H}$ is in the kernel of $\epsilon$, since  $\epsilon$ is an algebra map. Accordingly $\varpi$ maps onto left coinvariant forms. To show that $I=\ker\varpi$ is a right ideal we simply consider $\eta\in I$ and $h\in H$, to obtain \eqa{ \nonumber \varpi(\eta h) & = S(\eta_{1}h_{1})\dd(\eta_{2}h_{2}) \\ & =S(h_{1})S(\eta_{1}) (\dd \eta_{2} h_{2} +\eta_{2}\dd h_{1}) \\ & =  S(h_{1})\varpi(\eta)h_{2} + S(h_{1})\epsilon(\eta)\dd h_{2} \\ & = 0, } therefore $\eta h$ is in $I$.  \qed 
	 \smallskip
	 
	Denote the quotient map by \eqa{\pi:H^{+}\rightarrow H^{+}/I,\quad h\mapsto \pi(h)=[h].}

\begin{thm}
	For any right ideal $I\subseteq H$ with $I\subseteq H^{+}$ we have a left covariant first order differential calculus $(\Gamma,\dd)$ on $H$, where $$\Gamma=H\otimes (H^{+}/I),\quad \dd h=(\id \otimes \pi)(\Delta(h)-h\otimes 1). $$$\Gamma$ is an $H-$bimodule via \eqa{ h\cdot (h'\otimes [g])=hh'\otimes [g], \quad (h\otimes [g])\cdot h'=hh'_{1}\otimes [gh'_{2}],} for $h,h'\in H$ and $g\in H^{+}$. The left $H-$coaction on $\Gamma$ is ${}_{\Gamma}\Delta=\Delta \otimes \id_{H/I}$. 
	
	If $\ad_{R}(I)\subseteq I\otimes H$, $(\Gamma,\dd)$ is bicovariant with right $H-$coaction \eqa{\Delta_{\Gamma}(h\otimes [g]) =(h_{1}\otimes [g_{2}])\otimes h_{2}S(g_{1})g_{3},} for $h\in H$ and $g\in H^{+}$. Moreover, every $H-$covariant first order differential calculus is of this form.  
	\proof The bimodule relations do not depend upon the choice of a representative for a class under the projection map. Indeed $$h\cdot (h'\otimes [g])=h\cdot (h'\otimes [g']),$$ so we define the same element. The same for the other relation. Therefore $\Gamma$ is an $H-$bimodule;in particular a right $H-$module. 
	
	Consider the given definition of the $\B{k}-$linear map $\dd:H\rightarrow \Gamma$. We have \eqa{\nonumber \dd(h)h'+h\dd(h') & = (\id\otimes \pi)((h_{1}\otimes h_{2} -h\otimes 1)\cdot h' + h\cdot (h'_{1}\otimes [h'_{2}]-h'\otimes 1) \\ & = hh'_{1}\otimes [h_{2}h'_{2}]- hh'_{1}\otimes [h'_{2}] + hh'_{1}\otimes [h'_{2}] - hh'\otimes [1] \\ & = hh'_{1}\otimes [h_{2}h'_{2}] - hh'\otimes [1] \\ & = \dd (hh')} for all $h,h'\in H$. 
	
	Let $h^{i}\otimes [g^{i}]\in \Gamma$ be arbitrary, with $h^{i}\in H$ and $g^{i}\in H^{+}$. Consider the combination $h^{i}S(g_{1}^{i}) \dd(g^{i}_{2})$, defining an element of the form $H\dd H$. We find \eqa{\nonumber h^{i}S(g_{1}^{i})\dd (g_{2}^{i})& = h^{i}S(g_{1}^{i})g_{2}^{i}\otimes [g_{3}^{i}] - h^{i}S(g_{1}^{i})g_{2}^{i}\otimes [1]\\ &  = h^{i}\otimes [g^{i}]-h^{i}\epsilon(g^{i})\otimes [1]\\ &  =h^{i}\otimes [g^{i}].} Surjectivity follows. 
	
	The map ${}_{\Gamma}\Delta:\Gamma\rightarrow H\otimes \Gamma$ is a left $H-$coaction, as it is defined by $\Delta$ and the identity; the diagram $$ \begin{tikzcd} H\otimes \Gamma \arrow[rr,"\Delta\otimes \id_{\Gamma}"]& & H\otimes H \otimes \Gamma \\ \\ \Gamma \arrow[uu,"{}_{\Gamma}\Delta"] \arrow[rr,"{}_{\Gamma}\Delta"] & & H\otimes \Gamma \arrow[uu,"\id_{H}\otimes {}_{\Gamma}\Delta",swap] 
 \end{tikzcd}$$  commutes as \eqa{\nonumber(\id_{H}\otimes {}_{\Gamma}\Delta)(\Delta \otimes \id_{\Gamma})(h\otimes [g]) & = (\id_{H}\otimes {}_{\Gamma}\Delta)(h_{1}\otimes h_{2} \otimes [g]) \\ & = (\id_{H}\otimes \Delta \otimes \id_{\Gamma})(h_{1}\otimes h_{2}\otimes [g]) \\ & = h_{1}\otimes h_{2}\otimes h_{3}\otimes [g] \\& = (\Delta \otimes \id_{\Gamma})(h_{1}\otimes h_{2}\otimes [g]) \\ & = (\Delta\otimes \id_{\Gamma}){}_{\Gamma}\Delta(h\otimes [g]). } 
 
Compatibility with the $H-$bimodule structure follows  \eqa{\nonumber{}_{\Gamma}\Delta(h\cdot (h'\otimes [g])\cdot h'') & = {}_{\Gamma}\Delta(hh'h''_{1}\otimes [g h''_{2}]) \\ & = h_{1}h'_{1}h''_{1}\otimes h_{2}h'_{2}h''_{2} \otimes [gh''_{3}] \\ & = \Delta(h)\cdot {}_{\Gamma}\Delta(h'\otimes [g])\cdot \Delta (h''),} for all $h,h',h''\in H$ and $g\in H^{+}$. 

Left $H-$colinearity of the differential follows as \eqa{\nonumber(\id\otimes \dd)\Delta(h)&  = h_{1}\otimes \dd(h_{2}) \\ & = h_{1}\otimes h_{2}\otimes [h_{3}]-h_{1}\otimes h_{2}\otimes [1] \\ & ={}_{\Gamma}\Delta(\dd (h)).} Accordingly, $(\Gamma,\dd)$ is a left covariant first order differential calculus. 
	 
	 If in addition $\ad_{R}(I)\subseteq I\otimes H$, the map $\Delta_{\Gamma}:\Gamma\rightarrow\Gamma\otimes H$ is well defined and moreover a right $H-$coaction: \eqa{\nonumber(\Delta_{\Gamma}\otimes \id)(\Delta_{\Gamma}(h\otimes [g]) & = \Delta_{\Gamma}(h_{1}\otimes [g_{2}])\otimes h_{2}S(g_{1})g_{3} \\ & = h_{1}\otimes [g_{3}]\otimes h_{2}S(g_{2})g_{4}\otimes h_{3}S(g_{1})g_{5} \\ & = h_{1}\otimes [g_{2}] \otimes (h_{2}S(g_{1})g_{3})_{1}\otimes (h_{2}S(g_{1})g_{3})_{2} \\ & = (\id \otimes \Delta)(\Delta_{\Gamma}(h\otimes [g])), } for all $h\in H$ and $g\in H^{+}$. Proofs of compatibility of $\Delta_{\Gamma}$ with the bimodule structure and right $H-$colinearity of the differential are similar. We have shown that given a right ideal $I\subseteq H$ with $I\subseteq  H^{+}$ we can construct $(\Gamma,\dd)$ a bicovariant first order calculus on $H$.    
	 
	 \smallskip
	 
	 On the other hand, given a first order differential calculus we construct a right ideal $I\subseteq H$ as the kernel of $\varpi$ according to Lemma \ref{maurer cartan}

For the $1:1-$correspondence we prove that, given a left covariant first order differential calculus $\Gamma$, we obtain $\Gamma\cong H\otimes H^+/I$ as left covariant first order differential calculi. 
	 
	  The first thing we prove is the isomorphism $\Gamma\cong H\otimes H^{+}/I$ as left covariant $H-$bimodules. Consider \eqa{\phi:\Gamma\rightarrow H\otimes H^{+}/I,\quad \omega\mapsto \omega_{-2}\otimes \varpi^{-1}(S(\omega_{-1})\omega_{0}).} Here $\varpi^{-1}$ is the inverse of the quantum Maurer-Cartan form restricted as $\varpi:H^{+}/\ker\varpi \rightarrow \Lambda^{1}$, which is an isomorphism by Lemma \ref{maurer cartan}. Notice that $S(\omega_{-1})\omega_{0}$ is a left coinvariant form. 
	  
	  We provide the inverse of $\phi$ as \eqa{\phi^{-1}:H\otimes H^{+}/I\rightarrow \Gamma,\quad h\otimes [g]\mapsto h\varpi(g).} This map is $H-$bilinear and left $H-$colinear, indeed \eqa{\nonumber \phi^{-1}(h\cdot (h'\otimes [g]) \cdot h'') & = \phi^{-1}(hh'h''_{1}\otimes [gh''_{2}]) \\ & = h h' h''_{1}S(g_{1}h''_{2})\dd (g_{2}h''_{3}) \\ & = hh'h''_{1}S(h''_{2})S(g_{1}) \dd (g_{2}h''_{3}) \\ & = h h' S(g_{1})(\dd g_{2} h'' + g_{2} \dd h'') \\ & = hh'\varpi(h)h''+hh'\dd h''\epsilon(g) \\ & = h\cdot \phi^{-1}(h'\otimes [g])\cdot h'', } for all $h,h',h''\in H$ and $g\in H^{+}$. 
	  
	  Furthermore \eqa{ \nonumber (h\otimes [g])_{-1}\otimes \phi^{-1}((h\otimes [g])_{0}) & = h_{1}\otimes \phi^{-1}(h_{2}\otimes [g]) \\ & = h_{1}\otimes h_{2}\varpi(g) \\ & = (h\varpi(g))_{-1} \otimes (h\varpi(g))_{0} \\ & = \phi^{-1}(h\otimes [g])_{-1}\otimes \phi^{-1}(h\otimes [g])_{0}. } The last thing to show is the relation between the differential on $\Gamma$ and the one on $H\otimes H^{+}/I$. We have \eqa{ \nonumber \phi^{-1}(\dd_{I}(h)) & =\phi^{-1}(h_{1}\otimes [h_{2}] - h \otimes [1]) \\ & = h_{1}\varpi(h_{2})-h\varpi(1) \\ & = h_{1}S(h_{2})\dd (h_{3})-h\dd (1) \\ & = \dd h,} proving that $\phi$ is an isomorphism of left covariant first order differential calculi.   
	 
	 If we furthermore assume the calculus $(H\otimes H^{+}/I,\dd_{I})$ to be bicovariant, the ideal $I$ satisfies $\Delta_{\Gamma}(H\otimes I)\subseteq H \otimes I \otimes H$, in particular $$\Delta_{\Gamma}(1\otimes g)=1\otimes g_{2}\otimes S(g_{3})g_{1}\in H\otimes I \otimes H,$$ for all $g\in I$, namely $\ad_{R}(I)\subseteq I\otimes H$. In this case the isomorphism is a right $H-$covariant map, since \eqa{\nonumber\phi^{-1}((h\otimes [g]_{0}\otimes (h\otimes [g])_{1}) & = \phi^{-1}(h_{1}\otimes [g_{2}])\otimes h_{2}S(g_{1})g_{3} \\ & = h_{1}\varpi(g_{2})\otimes h_{2}S(g_{1})g_{3} \\ & = h_{1}S(g_{2})\dd(g_{3})\otimes h_{2}S(g_{1})g_{4} \\ & = (hS(g_{1})\dd (g_{2}))_{0}\otimes (hS(g_{1})\dd (g_{2}))_{1} \\ & = \phi^{-1}(h\otimes [g])_{0}\otimes \phi^{-1}(h\otimes [g])_{1},} for every $h\in H$ and every $g\in H^{+}$ by the right $H-$covariance of $\dd$.      
	 \qed 

\end{thm}

   \section{Higher order differential calculi}\label{higher order calculi}
   In this section we introduce higher order differential calculi as differential graded algebras with a surjectivity requirement. 
   
   \begin{defn}[Differential graded algebra]
   	A differential $\B{N}_{0}-$graded algebra $(\Omega^{\bullet},\wedge,\dd )$ is a $\B{N}_{0}-$graded vector space $\Omega^{\bullet}=\bigoplus_{k\geq 0} \Omega^{k}$ endowed with \begin{enumerate}
   		\item a map $\wedge:\Omega^{\bullet}\otimes \Omega^{\bullet}\rightarrow\Omega^{\bullet}$ of degree zero, i.e. $\Omega^{k}\wedge \Omega^{\ell}\subseteq \Omega^{k+\ell},$ 
   		\item a map $\dd:\Omega^{\bullet}\rightarrow\Omega^{\bullet+1}$ of degree $1$, meaning $\dd(\Omega^{k})\subseteq \Omega^{k+1}$,
   	\end{enumerate}
   	such that the \textit{wedge product} $\wedge$ is associative and unital and the differential satisfies $\dd^{2}=0$ and the graded Leibniz rule $$\dd(\omega\wedge \eta)=\dd \omega \wedge \eta +(-1)^{|\omega|}\omega \wedge \dd \eta$$ for all $\eta \in \Omega^{\bullet}$ and homogeneous\footnote{With homogeneous we simply mean an element that belongs to a fixed vector space in the grading. $|\omega|$ is the degree of $\omega$.} elements $\omega$. 
   \end{defn}

   \begin{defn}
  Let $(\Omega^{\bullet}(A),\dd_{A},\wedge_{A})$ over $A$ and $(\Omega^{\bullet}(B),\dd_{B},\wedge_{B})$ over $B$ be differential graded algebras. A morphism of differential graded algebras is a morphism of graded algebras $\Phi:\Omega^{\bullet}(A)\rightarrow\Omega^{\bullet}(B)$ that respects the differential. 
   \end{defn}
      We say that a morphism of differential graded algebras $\Phi:\Omega^{\bullet}(A)\rightarrow\Omega^{\bullet}(B)$ is an extension of a morphism of algebras $\phi:A\rightarrow B$ if $\Phi|_{A}=\phi.$ 
   	\begin{prop}\label{surjective morphism of DGA}
	Let $(\Omega^{\bullet}(A),\dd_{A},\wedge_{A})$ and $(\Omega^{\bullet}(B),\dd_{B},\wedge_{B})$ be differential graded algebras. Let $\Phi:\Omega^{\bullet}(A)\rightarrow\Omega^{\bullet}(B)$ be a morphism of differential graded algebras such that $\Phi|_{A}=\phi:A\rightarrow B$. We have that $\Phi$ is unique as an extension of $\phi$. Moreover if $\phi$ is surjective then $\Phi$ is.  
	\proof 
	Since the involved operations are linear it is enough to consider a homogeneous element $\omega=a^{0}\dd a^{1}\wedge\dots \wedge \dd a^{k}\in \Omega^{k}(A)$, for $a^{0},\dots, a^{k}\in A$.  We find \eqa{ \label{morphism of DGA} \Phi(a^{0}\dd_{A} a^{1}\wedge_{A}\dots \wedge_{A} \dd_{A} a^{k}) & = \Phi(a^{0}) \Phi(\dd_{A} a^{1})\wedge_{B} \cdots \wedge_{B} \Phi(\dd_{A} a^{k}) \\ & = \phi(a^{0}) \dd_{B}(\phi(a^{1}))\wedge_{B}\cdots\wedge_{B} \dd_{B}(\phi(a^{k})), } where we only used that $\Phi$ is a morphism of differential graded algebra. Accordingly we find that $\Phi$ is the unique extension of $\phi$.  
	Moreover, let $\phi$ be surjective and let $b^{0}=\phi(a^{0}),\dots, b^{k}=\phi(a^{k})\in B$ for some $a^{0},\dots, a^{k}\in A$. Then by equation \eqref{morphism of DGA} we have that $\Phi$ is surjective, since \eqa{ \nonumber b^{0}\dd_{B} b^{1}\wedge_{B}\cdots\wedge_{B}\dd_{B} b^{k} & = \phi(a^{0})\dd_{B}(\phi(a^{1})) \wedge_{B}\cdots\wedge_{B}\dd_{B} \phi(a^{k})\\ & = \Phi(a^{0}) \Phi(\dd_{A} a^{1})\wedge_{B} \cdots \wedge_{B} \Phi(\dd_{A} a^{k}) \\ & = \Phi(a^{0}\dd_{A} a^{1}\wedge_{A}\dots \wedge_{A} \dd_{A} a^{k}).} \qed 
\end{prop}
	\begin{defn}[Differential calculus]
   A differential calculus on a $\B{k}-$algebra $A$ is a differential $\B{N}_{0}-$graded algebra $(\Omega^{\bullet},\wedge,\dd)	$ which is generated in degree zero and such that $\Omega^{0}=A$. The former means that $$\Omega^{k}=\text{span}_{\B{k}}\{a^{0}\dd a^{1}\wedge\dots\wedge \dd a^{k}: a^{0},\dots,a^{k}\in A\}.$$     \end{defn}
     We call elements of $\Omega^{k}$ differential $k-$forms.  

\begin{lemma}\label{morphism of DGAs on fixed degree forms}
Let $(\Omega^{\bullet}(A),\dd_{A},\wedge_{A})$ be a differential calculus over $A$ and let $(\Omega^{\bullet}(B),\dd_{B},\wedge_{B})$ be a differential graded algebra.
	A map $\Phi:\Omega^{\bullet}(A)\rightarrow \Omega^{\bullet}(B)$ is a morphism of differential graded algebras if and only if \eqa{\label{morphism of DGAs on fixed degree forms}\Phi(a^{0}\dd_{A} a^{1}\wedge_{A}\cdots \wedge_{A}\dd_{A} a^{k}) = \Phi(a^{0})\wedge_{B} \dd_{B}\circ \Phi(a^{1})\wedge_{B}\cdots \wedge_{B} \dd_{B}\circ \Phi(a^{k}),} for all $a^0,...,a^k\in A,$ with $k\geq 0$.     
	\proof If $\Phi$ is a morphism of differential graded algebras then equation \eqref{morphism of DGAs on fixed degree forms} is clear. On the other hand assume equation \eqref{morphism of DGAs on fixed degree forms} to be true.  Let $a^{0}\dd_{A}a^{1}\wedge_{A}\cdots \wedge_{A} \dd a^{k} \in\Omega^{k}(A)$ and $ b^{0}\dd_{A}b^{1}\wedge_{A}\cdots \wedge_{A} \dd_{A}b^{\ell}\in \Omega^{\ell}(A).$ Then \eqa{\nonumber\Phi[(a^{0}\dd_{A}a^{1}& \wedge_{A}\cdots \wedge_{A} \dd a^{k})  \wedge_{A}(b^{0}\dd_{A}b^{1}\wedge_{A}\cdots \wedge_{A} \dd_{A}b^{\ell})] \\ &  = [\Phi(a^{0})\dd_{B} \Phi(a^{1})\wedge_{B}\cdots \wedge_{B}\dd_{B}\Phi(a^{k}) ] \wedge_{B} [\Phi(b^{0})\dd_{B} \Phi(b^{1})\wedge_{B}\cdots \wedge_{B}\dd_{B}\Phi(b^{\ell}) ]\\ & = \Phi(a^{1}\dd_{A}a^{1}\wedge_{A}\cdots \wedge_{A}\dd_{A}a^{k})\wedge_{B}\Phi(b^{1}\dd_{A}b^{1}\wedge_{A}\cdots \wedge_{A}\dd_{A}b^{\ell}). }
	
	Moreover, $\Phi$ is compatible with the differential, as \eqa{\nonumber\Phi(\dd_{A}\omega)& = \Phi(\dd_{A}(a^{0}\dd_{A}a^{1}\wedge\cdots \dd_{A}a^{k}) \\ & = \Phi(\dd_{A}a^{0}\wedge\cdots \dd_{A}a^{k}) \\ & = \dd_{B}\phi(a^{0}) \wedge_{B}\cdots \wedge_{B} \dd_{B}\phi(a^{k}) \\ & = \dd_{B} b^{0} \wedge_{B}\cdots \wedge_{B} \dd_{B} b^{k}} equals \eqa{\nonumber\dd_{B}\Phi(\omega) & = \dd_{B} \Phi(a^{0}\dd_{A} a^{1}\wedge_{A}\cdots \wedge_{A} \dd_{A} a^{k}) \\ & = \dd_{B}(b^{0}\dd_{B}b^{1} \wedge_{B}\cdots \wedge_{B} \dd_{B}b^{k}) \\ & = \dd b^{0}\wedge_{B}\cdots \wedge_{B}\dd_{B} b^{k}.}    \qed 
\end{lemma}

    \begin{exmp}[Trivial DC] On every algebra $A$ there is a trivial DC $\Omega^{k}=0$ for all $k> 0$ with the zero differential. 	
 \end{exmp}

    \begin{exmp}\label{universal differential calculus} The universal first order differential calculus of Example \ref{UDC} can be extended to higher order forms, defining $\Omega^{n}_{u}:=\cap_{k=1}^{n}\ker\mu_{k},$ where $\mu_{k}:=\id_{A^{\otimes(k-1)}}\otimes \mu\otimes \id_{A^{\otimes(n-k)}}:A^{\otimes(n+1)}\rightarrow A^{\otimes n}$. The differential is \eqa{\dd_{u}:\Omega_{u}^{k}\rightarrow \Omega_{u}^{k+1},\quad a^{0}\otimes \dots \otimes a^{k}\mapsto \sum_{i=1}^{k+1}a^{0}\otimes \dots \otimes a^{i-1}\otimes 1 \otimes a^{i+1} \otimes \dots \otimes a^{k}.} The wedge product $\wedge_{u}:\Omega_{u}^{k}\otimes \Omega_{u}^{\ell}\rightarrow\Omega^{k+\ell}$ is defined on a $k-$form $\omega$ and a $\ell-$form $\eta$ as \eqa{ \omega\wedge_{u} \eta = a^{0}\otimes \dots \otimes a^{k-1} \otimes a^{k}b^{0}\otimes b^{2} \dots \otimes b^{\ell},} i.e. as the concatenation of forms. \begin{lemma}
		The universal differential calculus $(\Omega^{\bullet}_{u},\wedge_{u},\dd_{u})$ is a differential calculus on $A$. 
		\proof During the proof we will denote $\dd_{u}$ as $\dd$ and $\wedge_{u}$ as $\wedge$ for simplicity.
		The $A-$bimodule structure on $\Omega_{u}^{\bullet}$ is the one induced from $A\otimes A$, i.e. the multiplication on the first tensor factor from the left and the second tensor factor from the right, respectively. The squared differential vanishes, as on the order $0$ we find  \eqa{\nonumber \dd^{2}(a) & = \dd (1\otimes a-a\otimes 1) \\ & =\dd(1\otimes a)-\dd (a\otimes 1) \\& = 1\otimes 1 \otimes a - 1\otimes 1 \otimes a +1\otimes a \otimes 1 \\ & \ \ - (1\otimes a \otimes 1 - a\otimes 1 \otimes 1 + a\otimes 1 \otimes 1) \\ & = 0, } and moreover at higher orders we find \eqa{ \nonumber \dd^{2}(a^{0}\otimes \dots \otimes a^{k}) & = \dd (\sum_{i=0}^{k+1} a^{0}\otimes \dots \otimes a^{i-1}\otimes 1 \otimes a^{i}\otimes \dots \otimes a^{k}) \\ & = \dd(1\otimes a^{0} \otimes \dots \otimes a^{k}) - \dd (a^{0}\otimes 1 \otimes \dots \otimes a^{k}) +\dots + (-1)^{k}\dd(a^{0}\otimes \dots \otimes a^{k}\otimes 1) \\ & = 1\otimes 1 \otimes a^{0}\otimes \dots \otimes a^{k}-1\otimes 1 \otimes a^{0}\otimes \dots \otimes a^{k} + 1 \otimes a^{0}\otimes 1 \otimes \dots \otimes a^{k} +\dots \\ &- 1 \otimes a^{0}\otimes 1 \otimes \dots \otimes a^{k} + a^{0}\otimes 1 \otimes 1 \otimes \dots \otimes a^{k} +\dots \\ & = 0. } 
		We already know the Leibniz rule holds at order zero, namely on elements of $A$. For elements in $\Gamma_{u}$ we find \eqa{ \nonumber \dd(\omega\wedge \gamma) & =  \dd ((a^{0}\otimes a^{1}) \wedge (b^{0}\otimes b^{1})) \\ & = \dd ( a^{0}\otimes a^{1}b^{0}\otimes b^{1}) \\ & =  1\otimes a^{0}\otimes a^{1}b^{0}\otimes b^{1} - a^{0}\otimes 1 \otimes a^{1}b^{0} \otimes b^{1}  \\ & \ \ + a^{0}\otimes a^{1}b^{0}\otimes 1 \otimes b^{1} - a^{0}\otimes a^{1}b^{0}\otimes b^{1}\otimes 1 \\ & = (1\otimes a^{0}\otimes a^{1} - a^{0}\otimes 1 \otimes a^{1} + a^{0}\otimes a^{1}\otimes 1)\wedge(b^{0}\otimes b^{1})\\ &  -(a^{0}\otimes a^{1})\wedge (1\otimes b^{0}\otimes b^{1} - b^{0}\otimes 1 \otimes b^{1} + b^{0}\otimes b^{1}\otimes 1)\\ & =  \dd \omega \wedge \gamma - \omega \wedge\dd\gamma,     } and we can infer that the rule continues to hold at higher orders.  
		Finally surjectivity of the calculus follows as any element $a^{0}\otimes \dots \otimes a^{n}$ can be obtained as \eqa{\nonumber a^{0}\dd a^{1}\wedge\dots \wedge \dd a^{n}  & = a^{0}(1\otimes a^{1}-a^{1}\otimes 1 )\wedge \dots \wedge (1\otimes a^{n}-a^{n}\otimes 1) \\ & = (a^{0}\otimes a^{1})\wedge(1\otimes a^{2}-a^{2}\otimes 1)\wedge\dots \wedge(1\otimes a^{n}-a^{n}\otimes 1)\\ & = (a^{0}\otimes a^{1}\otimes a^{2})(1\otimes a^{3}-a^{3}\otimes 1)\wedge\dots \wedge (1\otimes a^{n}-a^{n}\otimes 1)\\  & = \cdots \\ & = a^{0}\otimes \dots \otimes a^{n}. } Therefore we have $(\Omega_{u}^{\bullet},\wedge_{u},\dd_{u})$ is a differential calculus on $A$.  \qed  
	\end{lemma} 
  \end{exmp}
\begin{defn}
	A morphism from a differential calulus $(\Omega^{\bullet},\wedge,\dd)$ over an algebra $A$ to another differential calculus $(\Omega'^{\bullet},\wedge',\dd')$ over an algebra $A'$ is a morphism of the underlying differential $\B{N}_{0}-$graded algebras. Explicitly, it is a graded map $\Phi:\Omega^{\bullet}\rightarrow \Omega'^{\bullet}$ of degree $0$, such that \eqa{\Phi(\omega\wedge \eta) =\Phi(\omega)\wedge \Phi(\eta),\quad \Phi\circ \dd = \dd'\circ \Phi,} for all $\omega,\eta\in \Omega^{\bullet}$. An isomorphism of differential calculi is a morphism of differential calculi which is also invertible.
\end{defn}

 \subsection{Higher order $H-$covariant calculi} \label{higher order covariant calculi}In the following we discuss a special class differential calculi. Those are differential calculi on $H-$comodule algebras such that the coaction \textit{lifts} to graded algebra of forms and the differential is colinear, and are called differential covariant calculi. A $\B{N}_{0}-$graded $\B{k}-$vector space $V^{\bullet}=\bigoplus_{k\geq 0}V^{k}$ is a graded right $H-$comodule if there is a right $H-$coaction $\Delta^{\bullet}:V^{\bullet}\rightarrow V^{\bullet}\otimes H$ which is graded of degree $0$\footnote{Meaning that $\Delta^{\bullet}(V^{k})\subseteq V^{k}\otimes H$ for all positive k's.}. It follows that for any $k\geq 0$ the vector space $V^{k}$ is a right $H-$comodule via  $$\Delta^{k}:=\Delta^{\bullet}|_{V^{k}}:V^{k}\rightarrow V^{k}\otimes H;$$ accordingly we write $\Delta^{\bullet}=\oplus_{k\geq 0} \Delta^{k}$. 
 
 A graded map $\phi:V^{\bullet}\rightarrow W^{\bullet+\ell}$ of degree $\ell$ between graded right $H-$comodules $(V^{\bullet},\Delta_{V}^{\bullet})$ and $(W^{\bullet},\Delta^{\bullet}_{W})$ is right $H-$colinear if the diagram $$\begin{tikzcd}
  V^{k} \arrow[dd,"\Delta_{V}^{k}"'] \arrow[rr,"\phi^{k}"] & & W^{k+\ell}\arrow[dd,"\Delta_{W}^{k+\ell}"] \\ \\
  V^{k}\otimes H \arrow[rr,"\phi^{k}\otimes \id_{H}"'] & & W^{k+\ell}\otimes H 	
\end{tikzcd}$$ 
 commutes for all $k\geq 0$, where here $\phi^{k}=\phi|_{V^{k}}:V^{k}\rightarrow W^{k+\ell}$. 
 
 In case $V^{\bullet}$ is a graded algebra we call it a \textit{graded right $H-$comodule algebra} if it is a graded right $H-$comodule and the algebra structure is compatible with the right $H-$coaction in the usual sense.  
 \begin{defn}
 	A differential calculus over a right $H-$comodule algebra $(A,\Delta_{A})$ is called \textit{right $H$-covariant} if $(\Omega^{\bullet},\wedge)$ is a graded right $H-$comodule algebra with graded coaction $\Delta^{\bullet}:\Omega^{\bullet}\rightarrow\Omega^{\bullet}\otimes H$ such that $\Delta^{0}=\Delta_{A}:A\rightarrow A\otimes H$ and the differential $\dd:\Omega^{\bullet}\rightarrow\Omega^{\bullet+1}$ is right $H-$colinear. Similarly, left covariant and bicovariant calculi are defined. 
 \end{defn}
 If we consider $H$ as an $H-$comodule algebra with corresponding coaction given by the coproduct, we call the corresponding $H-$covariant calculus \textit{covariant}.

 \begin{prop} \label{equivalent definitions of covariance}
 	A differential calculus $(\Omega^{\bullet},\wedge,\dd)$ on a right $H-$comodule algebra $A$ is right $H-$covariant if and only if one of the following conditions holds: \begin{enumerate}
 		\item  for all $k>0$ and $a^{i,0},\dots,a^{i,k}\in A$ such that $\sum_{i}a^{i,0}\dd a^{i,1} \wedge \dots \wedge \dd a^{i,k}=0$ we have $$ \sum_{i}a_{0}^{i,0}\dd a_{0}^{i,1} \wedge \dots \wedge \dd a_{0}^{i,k} \otimes a_{1}^{i,0}a_{1}^{i,1}\dots a_{1}^{i,k}=0; $$ 
 		\item the right $H-$coaction $\Delta_{A}:A\rightarrow A\otimes H$ is $k-$times differentiable for all $k>0$, namely $$\begin{tikzcd}
 	A\otimes H \arrow[rr,"\dd_{\otimes}"]  & & \Omega^{1}\otimes H \arrow[rr, "\dd_{\otimes}"] & & \dots \arrow[rr, "\dd_{\otimes}"] & &  \Omega^{k}\otimes H \\ \\
 	A \arrow[uu,"\Delta_{A}"] \arrow[rr,"\dd"] & & \Omega^{1} \arrow[uu," \Delta_{\Omega^{1}(A)}"] \arrow[rr,"\dd"] & & \dots \arrow[rr,"\dd"] & & \Omega^{k}\arrow[uu,"\Delta_{\Omega^{k}(A)}"] 
\end{tikzcd}
$$ commute, and for all $k>0$ the map $\Delta_{\Omega^{k}(A)}:\Omega^{k}\rightarrow\Omega^{k}\otimes H$ is left $A-$linear, in the sense that $\Delta_{\Omega^{k}(A)}(a\cdot \omega) = \Delta_{A}(a)\Delta_{\Omega^{k}(A)}(\omega)$ for all $a\in A$ and $\omega\in \Omega^{k}$.  

 	\end{enumerate}
 	In this case the right $H-$coactions $\Delta_{\Omega^{k}(A)}:\Omega^{k}\rightarrow\Omega^{k}\otimes H$, $k>0$ are explicitely given as $$ \Delta_{\Omega^{k}(A)} \tonde{ \sum_{i}a^{i,0} \dd a^{i,1} \wedge \dots \wedge \dd a^{i,k}} = \sum_{i}a_{0}^{i,0}\dd a_{0}^{i,1} \wedge\dots \wedge \dd a_{0}^{i,k} \otimes a_{1}^{i,0}a_{1}^{i,1}\dots a_{1}^{i,k} $$  for all $a^{i,0},\dots,a^{i,k} \in A$.  
 	\proof 
 	Let $(\Omega^{\bullet},\wedge,\dd)$ be a differential calculus on a right $H-$comodule algebra $(A,\Delta_{A})$. First we prove that $1.$ implies right $H-$covariance. 
 	
 	 For $k>0$ define a $\B{k}-$linear map \eqa{\nonumber\Delta_{\Omega^{k}(A)}:\Omega^{k}(A)\rightarrow\Omega^{k}(A)\otimes H, \quad \sum_{i}a^{i,0}\dd a^{i,1} \wedge \dots \wedge \dd a^{i,k} \mapsto \sum_{i}a_{0}^{i,0}\dd a_{0}^{i,1}\wedge\dots \wedge \dd a_{0}^{i,k} \otimes a_{1}^{i,0}a_{1}^{i,1}\dots a_{1}^{i,k}. \nonumber} The assumption of point $1$ ensures this expression is well-defined, since it is zero whenever the input object is zero.
 
 	 Moreover $\Delta_{\Omega^{k}(A)}:\Omega^{k}(A)\rightarrow\Omega^{k}(A)\otimes H$ is a right $H-$coaction: \eqa{\nonumber(\Delta_{\Omega^{k}(A)}\otimes \id)\Delta_{\Omega^{k}(A)}(\omega) & = \sum_{i}(\Delta_{\Omega^{k}(A)}\otimes\id)(a_{0}^{i,0}\dd a_{0}^{i,1} \wedge \dots \wedge \dd a_{0}^{i,k} \otimes a_{1}^{i,0}\dots a_{1}^{i,k}) \\ & = \sum_{i}a_{0}^{i,0}\dd a_{0}^{i,1}\wedge\dots \wedge \dd a_{0}^{i,k} \otimes a_{1}^{i,0}\dots a_{1}^{i,k} \otimes a_{2}^{i,0}\dots a_{2}^{i,k} \\ & = \sum_{k}(\id \otimes \Delta)(a_{0}^{i,0}\dd a_{0}^{i,1}\wedge \dots \wedge \dd a_{0}^{i,k} \otimes a_{1}^{i,0}\dots a_{1}^{i,k} ) \\ & = (\id \otimes \Delta)\Delta_{\Omega^{k}(A)}(\omega),} for every $\omega \in \Omega^{k}$.  
 	
 	  Consequently, $\Omega^{\bullet}$ is a graded right $H-$comodule with respect to $\Delta_{\Omega^{\bullet}(A)}=\Delta_{A}\oplus \bigoplus_{k>0}\Delta_{\Omega^{k}(A)}$.

The differential $\dd|_{\Omega^{k}}:\Omega^{k}(A)\rightarrow\Omega^{k+1}(A)$ is right $H-$colinear for every $k\geq 0$. Indeed, let $\omega \in \Omega^{k}$ . By the Leibniz rule and $\dd^{2}=0$ we have $$\dd(\omega)=\sum_{i}\dd a^{i,0}\wedge \dd a_{0}^{i,1} \wedge \dots \wedge \dd a_{0}^{i,k}.$$ Therefore \eqa{\nonumber\Delta_{\Omega^{k+1}(A)}\circ \dd (\omega )& =\sum_{i}\Delta^{k+1}(\dd a^{i,0}\wedge \dots \wedge \dd a^{i,k}) \\ & = \sum_{i}\dd a_{0}^{i,0}\wedge\dots \wedge \dd a_{0}^{i,k}\otimes a_{1}^{i,0}\dots a_{1}^{i,k} \\ & = \sum_{i}(\dd \otimes \id_{H})(a_{0}^{i,0}\dd a_{0}^{i,0}\wedge \dots \wedge \dd a_{0}^{i,k} \otimes a_{1}^{i,0}\dots a_{1}^{i,k}) \\ & = (\dd \otimes \id_{H})\circ\Delta_{\Omega^{k}(A)}(\omega). }

To prove $\Delta_{\Omega^{\bullet}(A)}$ is a graded algebra map we show $\Delta_{\Omega^{\bullet}(A)}(\omega\wedge \eta)=\Delta_{\Omega^{\bullet}(A)}(\omega)(\wedge\otimes \mu)\Delta_{\Omega^{\bullet}(A)}(\eta)$ for all $\omega\in \Omega^{k}(A)$ and $\eta\in \Omega^{\ell}(A)$. For $k=0$ this is automatic, since $\Delta_{A}$ is an algebra map. For $k=1$, using the Leibniz rule, we find \eqa{ \nonumber\Delta_{\Omega^{1+\ell}(A)} & \tonde{ (\sum_{i}a^{i,0}\dd a^{i,1})\wedge (\sum_{j}b^{j,0}\dd b^{j,1}\wedge\dots\wedge \dd b^{j,\ell})} \\ & =  \sum_{i,j} \Delta^{1+\ell} \tonde{a^{i,0}\dd (a^{i,1}b^{j,0})\wedge \dd b^{j,1} \wedge\dots\wedge \dd b^{j,\ell} } - \sum_{i,j}\Delta^{1+\ell}\tonde{a^{i,0}a^{i,1}\dd b^{j,0} \wedge \dd b^{j,1} \wedge \dots \wedge \dd b^{j,\ell}} \\ & = \sum_{i,j}\tonde{a_{0}^{i,0}\dd(a_{0}^{i,1}b^{j,0})\wedge\dd b_{0}^{j,1} \wedge\dots \wedge\dd b_{0}^{j,1}-a_{0}^{i,0}a_{0}^{i,1}\dd b^{j,0} \wedge \dots\wedge \dd b^{j,\ell}}\otimes a_{1}^{i,0}a_{1}^{i,1}b_{1}^{j,0}b_{1}^{j,1}\dots b_{1}^{j,\ell} \\ & = \sum_{i,j}a_{0}^{i,0}\dd a_{0}^{i,1}\wedge (b_{0}^{j,0}\dd b_{0}^{j,1})\wedge\dots \wedge \dd b_{0}^{j,1}\otimes a_{1}^{i,0}a_{1}^{i,1}b_{1}^{j,0}b_{1}^{j,1}\dots b_{1}^{j,\ell} \\ & = \Delta_{\Omega^{1}(A)}\tonde{\sum_{i}a^{i,0}\dd a^{i,1}}\Delta_{\Omega^{\ell}(A)}\tonde{\sum_{j}b^{j,0}\dd b^{j,1}\wedge\dots \wedge \dd b^{j,\ell}}} holds.  

The cases $k>1$ are proven by the same argumentation above by rewriting $\omega\wedge \eta $ in the form $\sum_{i}c^{i,0}\dd c^{i,1}\wedge\dots\wedge \dd c^{i,k+\ell}$ for some $c^{i,0},\dots, c^{i,k+\ell}$s in $A$. 
Thus, we have $(\Omega^{\bullet},\wedge,\dd)$  is a right $H-$covariant differential calculus.  
 	
 	\smallskip 
 	On the other hand assuming the differential calculus to be right $H-$covariant we want to show the assumption $1.$ holds. Consider $\omega\in \Omega^{k}$ and apply $\Delta_{\Omega^{k}(A)}$ to obtain  \eqa{\nonumber\Delta_{\Omega^{k}(A)}\tonde{\sum_{i}a^{i,0}\dd a^{i,1}\wedge\dots\wedge \dd a^{i,k}} & = \sum_{i}\Delta_{A}(a^{i,0})\Delta_{\Omega^{1}(A)}(\dd a^{i,1})\cdots \Delta_{\Omega^{1}(A)}(\dd a^{i,k}) \\ & = \sum_{i}\Delta_{A}(a^{i,0})(\dd \otimes \id)\Delta_{A}(a^{i,1})\cdots (\dd \otimes \id)\Delta_{A}(a^{i,k}) \\ & = \sum_{i} a_{0}^{i,0}\dd a_{0}^{i,1} \wedge \dots \wedge \dd a_{0}^{i,k} \otimes a_{1}^{i,0} a_{1}^{i,1}\dots a_{1}^{i,k}, }  since $\Delta_{\Omega^{\bullet}(A)}$ is an algebra map and $\dd$ is colinear.  In particular $\omega=0$ implies $\Delta^{k}(\omega)=0$ and thus we have the equivalence of point $1.$ with right $H-$covariance of the differential calculus.  
 	
 	For the point $2.$, if the differential calculus is right $H-$covariant with respect to the graded right $H-$coaction $\Delta_{\Omega^{\bullet}(A)}:\Omega^{\bullet}\rightarrow \Omega^{\bullet}\otimes H$ we have $\Delta_{\Omega^{k}(A)}$ is differentiable with $\dd \Delta_{\Omega^{k}(A)}=\Delta_{\Omega^{k+1}(A)}$. By assumption $\Delta_{\Omega^{\bullet}(A)}$ is an algebra map and thus $\Delta_{\Omega^{k}(A)}(a\cdot\omega)=\Delta_{A}(a)\Delta_{\Omega^{k}(A)}(\omega)$ for every $a\in A$ and $\omega\in \Omega^{k}$, proving the first implication. 

On the other hand by left linearity it is straightforward to prove $\Delta_{\Omega^{k}(A)}$ is of the required form for all positive values of $k$. Therefore, as the assumption of $1.$ are satisfied we have right $H-$covariance of the calculus. \qed
 \end{prop}
 
 \section{The maximal prolongation}
 \label{maximal prolongation}

Let us consider an associative unital algebra $A$ over a field $\B{k}$. Starting from $(\Gamma,\dd)$  a first order differential calculus over $A$, we want to construct a differential calculus over $A$ as a prolongation of the first order.

Let us consider the tensor bundle algebra \eqa{\Gamma^{\otimes_{A}}=\bigoplus_{i=0}^{\infty}\Gamma^{\otimes_{A}^{n}} = A\oplus\Gamma\oplus(\Gamma\otimes_{A}\Gamma)\oplus\dots,} that is a graded, associative, unital algebra with product $\otimes_{A}$ and unit $1\in A$.  Consider the graded ideal $$S^{\wedge}=\bigoplus_{n=0}^{\infty}S^{\wedge}\cap \Gamma^{\otimes_{A}^{n}}\subseteq \Gamma^{\otimes_{A}},$$ generated by elements $\sum_{i}\dd a^{i}\otimes_{A}\dd b^{i}$ where $a^{i},b^{i}\in A$, such that $\sum_{i}a^{i}\dd b^{i}=0$. 

 We define accordingly the graded associative unital algebra $\Gamma^{\wedge}:=\Gamma^{\otimes_{A}}/S^{\wedge}$, with induced product $\wedge$. The quotient algebra $\Gamma^{\wedge}$ is graded, indeed  as \eqa{\nonumber\Gamma^{\wedge} & = \frac{\bigoplus_{n=0}^{\infty} \Gamma^{\otimes_{A}^{n}}}{\bigoplus_{n=0}^{\infty} S^{\wedge}\cap \Gamma^{\otimes_{A}^{n}}} \\ & = \bigoplus_{n=0}^{\infty}x\frac{\Gamma^{\otimes_{A}^{n}}}{S^{\wedge}\cap\Gamma^{\otimes_{A}^{n}}}\\ & \cong \bigoplus_{n=0}^{\infty} \frac{\Gamma^{\otimes_{A}^{n}}+S^{\wedge}}{S^{\wedge}}, }where in the last row we exploited the second isomorphism theorem. Accordingly the product $\wedge:\Gamma^{\wedge}\otimes \Gamma^{\wedge}\rightarrow\Gamma^{\wedge}$ is induced as \eqa{\nonumber(\Gamma^{\wedge^{j}})\wedge (\Gamma^{\wedge^{k}}) & = \frac{(\Gamma^{\otimes_{A}^{j}}+S^{\wedge})}{S^{\wedge}} \otimes_{A} \frac{(\Gamma^{\otimes_{A}^{k}}+S^{\wedge})}{S^{\wedge}}\\ & \subseteq \frac{(\Gamma^{\otimes^{k+j}_{A}}+S^{\wedge})}{S^{\wedge}}.}
 
We report the following important result. 
\begin{thm}[\cite{Durdevic-1},Appendix B]\label{Maximal prolongation: extension of the differential}The differential $\dd:A\rightarrow \Gamma$ uniquely extends to a differential on $\Gamma^{\wedge}$ such that $(\Gamma^{\wedge},\wedge,\dd)$ is a differential calculus on $A$.
	\proof To show existence, we notice that the assignment $$\dd (\sum_{i}a_{i}\dd b_{i}) = \sum_{i}\dd a_{i}\wedge \dd b_{i}$$  defines a linear map $\dd:\Gamma\rightarrow\Gamma^{\wedge}$. This map is well defined since whenever $\sum_{i}a_{i}\dd b_{i}=0$ we have $\sum_{i}\dd a_{i}\wedge \dd b_{i}=0$ as we quotiented out the ideal $S^{\wedge}$. 
	
	The $\dd^{2}(a)=0$ property for each $a\in A$ is easily checked as \eqa{ \label{dd=0 for the maximal prolongation}\dd^{2}a= \dd(\dd a) = \dd(1)\wedge \dd a=0 .}  
	
	Moreover, let $a\in A$ and let $\theta=\sum_{i}a_{i}\dd b_{i}\in \Gamma$. We find \eqa{\nonumber
& \dd (a\theta) &&= \dd (\sum_{i}aa_{i}\dd b_{i})  \\
&&&= \sum_{i} \dd (aa_{i})b_{i} \\
&&&= \sum_{i} \dd (a)a_{i} \wedge \dd b_{i} + a \dd (a_{i})\wedge\dd b_{i} \\
&&&= \dd(a)\wedge\sum_{i}a_{i}\dd b_{i} + a \sum_{i}\dd a \wedge\dd b_{i} \\
&&&= \dd(a)\wedge\theta + a\dd(\theta); } and \eqa{ \nonumber\dd (\theta a) &&& = \dd (\sum_{i}a^{i}\dd b^{i} a) \\
&&&= \dd\sum_{i}( a_{i}\dd(b_{i}a)- a_{i}b_{i}\dd a) \\
&&&= \sum_{i}(\dd a_{i} \wedge\dd (b_{i}a) -  \dd(a_{i}b_{i}) \wedge\dd a) \\
&&&= \sum_{i} \dd a_{i} \wedge\dd b_{i} \ a + \dd a_{i} b_{i} \wedge\dd a - \dd a_{i} b_{i} \wedge\dd a - a_{i} \dd b_{i} \wedge \dd a \\
&&&= \sum_{i}(\dd a_{i}\wedge \dd b_{i})a - \sum_{i}(a_{i}\dd b_{i})\wedge\dd a  \\
&&&= \dd(\theta)a - \theta \wedge\dd a.
} 
 From the last two properties we have the map $\dd$ has the unique extension $\dd^{\wedge}:\Gamma^{\otimes_{A}^{2}}\rightarrow\Gamma^{\wedge}$, satisfying  \eqa{\label{extension of d to the maximal prolongation} \dd^{\wedge}(w\otimes_{A}u)=\dd(w)\wedge\pi(u)+(-1)^{|w|}\pi(w)\wedge\dd(u),}  where $\pi:\Gamma^{\otimes_{A}}\rightarrow\Gamma^{\wedge}$ is the projection map.  
 
 Equations \eqref{dd=0 for the maximal prolongation} and \eqref{extension of d to the maximal prolongation} imply elements of $S^{\wedge}$ are in the kernel of the differential $\dd^{\wedge}$, since \eqa{ \dd^{\wedge}\tonde{\sum_{i}\dd a_{i}\otimes_{A}\dd b_{i}} =\sum_{i}\dd^{2}a \wedge \pi(\dd b_{i})+(-1)^{|\dd a|}\ \pi(\dd a)\wedge \dd^{2}b=0.} 
 
 This extension can be generalised inductively to higher order forms. Consequently there exists a unique map $\dd:\Gamma^{\wedge}\rightarrow\Gamma^{\wedge}$ defined as a factorisation of the map $\dd^{\wedge}:\Gamma^{\otimes_{A}}\rightarrow \Gamma^{\wedge}$ through $\pi$ as $$\begin{tikzcd}
 \Gamma^{\otimes_{A}} \arrow[dd,"\pi",swap] \arrow[rr,"\dd^{\wedge}"] & & \Gamma^{\wedge} \\ \\ \Gamma^{\wedge}	\arrow[uurr,"\dd",swap] 
 \end{tikzcd}.$$
This map features all the required properties.   \qed \end{thm}

In example \ref{UDC} we introduced the universal first order differential calculus. Then in Theorem \ref{Woronowicz classification} we showed how any first order differential calculus can be induced as the quotient of the universal calculus. One may ask if this is also true when extended to higher order forms. 

We recall the universal differential calculus $(\Omega^{\bullet}_{u},\wedge_{u},\dd_{u})$ of Example \ref{universal differential calculus}.  Here the space of forms is $\Omega_{u}^{\bullet}:=\cap_{k=1}^{n}\ker \mu_{k}$, where $\mu_{k}:=\id_{A^{\otimes(k-1)}}\otimes \mu \otimes \id_{A^{\otimes(n-k)}}:A^{\otimes(n+1)}\rightarrow A^{\otimes n}$ denotes the multiplication of two subsequent tensor factors in the position $k$, and where the differential is given by $\dd_{u}(a^{0}\otimes\dots \otimes a^{n}):=\sum_{k=0}^{n+1}(-1)^{k}a^{0}\otimes \dots \otimes a^{k-1}\otimes 1 \otimes a^{k}\otimes \dots \otimes a^{n}$. 
The wedge product of forms is $(a^{0}\otimes \dots \otimes a^{k})\wedge_{u}(b^{0}\otimes \dots \otimes b^{\ell})=a^{0}\otimes \dots \otimes a^{k}b^{0}\otimes \dots b^{\ell}$, for all $a^{0},\dots,a^{k},b^{0},\dots,b^{\ell}\in A$  and the bimodule structure on $\Omega_{u}^{\bullet}$ is given as left multiplication on the first tensor factor and right multiplication on the last tensor factor. 
	
	\begin{lemma}
		Let $(\Omega^{\bullet},\wedge,\dd)$ be a differential calculus on $A$. Then $(\Omega^{\bullet},\wedge,\dd)$  is induced by a surjective morphism of differential graded algebras $\Omega_{u}^{\bullet}\rightarrow \Omega^{\bullet}$. In other words, any differential calculus over an algebra is a quotient of the universal differential calculus. 
		\proof 
		Following the lines in the proof of Theorem \ref{Woronowicz classification} we define a map \eqa{\Phi^{\bullet}:\Omega^{\bullet}_{u}\rightarrow \Omega^{\bullet}, \quad a^{0}\otimes \dots \otimes a^{k} \mapsto a^{0}\dd a^{1} \wedge \dots \wedge \dd a^{k}.} For the action on the space of $n-$forms we write \eqa{\Phi^{n}:\Omega^{n}_{u}\rightarrow\Omega^{n}, \quad \omega\mapsto \cdot\circ ( \id_{A}\wedge \dd \wedge \cdots \wedge \dd)\circ \omega,} where $\\cdot{A}:A\otimes \Omega^{\bullet}_{u}\rightarrow\Omega^{\bullet}_{u}$ is the left $A-$action on $\Omega^{\bullet}_{u}$. 
		
		We notice that \eqa{\Phi^{n}(0\otimes \dots \otimes 0) & = \cdot\circ(0\otimes\dots\otimes 0)\\ &  = 0\wedge\dots \wedge 0,} therefore $\Phi^{n}$ and thus $\Phi:=\bigoplus_{n\in \B{N}} \Phi^{n}$ are well defined maps.   
		
		Define the projection map \eqa{\pi:A\otimes A \rightarrow \ker \mu, \quad a\otimes b \mapsto a\otimes b - ab\otimes 1.} Accordingly, we define the $k-$th projection between two successive adjacent tensor factors  as $ \pi^{k}:=\id_{A^{\otimes(k-1)}}\otimes \pi\otimes \id_{A^{\otimes(n-k-1)}}$.  Let also $\xi_{n}=\pi_{n}\cdots \pi_{1}(a^{0}\otimes \dots \otimes a^{n})\in A^{\otimes(n+1)}$. As the various projections map to kernel of multiplications, we can infer that $\xi_{n}$ is an element of $\Omega^{n}_{u}$. We find \eqa{\nonumber\Phi^{n}(\xi_{n}) & = \Phi^{n}(\pi_{1}\dots \pi_{n}(a^{0}\otimes \dots \otimes a^{n}))\\ & = \cdot \circ (\id_{A}\otimes \dd \otimes \dots \otimes \dd)(\pi_{n}\dots \pi_{1}(a^{0}\otimes \dots \otimes a^{n})) \\  & = \cdot \circ(\id_{A}\otimes \dd \otimes \dots \otimes \dd)\circ (a^{0}\otimes a^{1}-a^{0}a^{1}\otimes 1)\otimes\dots \otimes (a^{n-1}\otimes a^{n}-a^{n-1}a^{n}\otimes 1) \\ & = \cdot \circ(\id_{A}\otimes \dd \otimes \dots \otimes \dd)(a^{0}\otimes \dots \otimes a^{n}) =\\ & = a^{0} \dd a^{1}\wedge\dots\wedge \dd a^{n},  }  showing the surjectivity of $\Phi^{n}$ and thus of $\Phi^{\bullet}$.  
		
		 The map $\Phi^{\bullet}$ is left $A-$linear, since \eqa{\nonumber\Phi^{n}(a (a^{0}\otimes\dots\otimes a^{n})) & = a a^{0}\dd a^{1}\wedge\dots \wedge \dd a^{n}\\ & = a \Phi^{n}(a^{0}\otimes \dots \otimes a^{n}).}
		
		 For right $A-$linearity consider \eqa{\nonumber\Phi^{n}((a^{0}\otimes\dots \otimes a^{n})) & = a^{0}\dd a^{1} \wedge \dots \wedge \dd(a^{n}a) \\ & = (a^{0}\dd a^{1}\wedge\dots\wedge \dd a^{n})a + a^{0}\dd a^{1}\wedge\dots \wedge \dd a^{n-1}a^{n}\wedge \dd a \\ &= \Phi^{n}(a^{0}\otimes\dots \otimes a^{n})a + \omega \wedge\dd(a), } and so right $A-$linearity follows if and only if the form $\omega$ vanishes, which is the case since \eqa{ \nonumber\omega & = a^{0}\dd a^{1}\wedge\dots \wedge \dd (a^{n-1})a^{n} \\ & = a^{0}\dd a^{1}\wedge\dots\wedge \dd (a^{n-1}a^{n}) - a^{0}\dd a^{1}\wedge\dots \wedge \dd (a^{n-2})a^{n-1}\wedge \dd a^{n}  \\ & = -a^{0}\dd a^{1}\wedge\dots \wedge \dd a^{n-2} a^{n-1} \wedge \dd a^{n} \\ & = (-1)^{2}a^{0}\dd a^{1}\wedge\dots \wedge\dd(a^{n-3})a^{n-2}\wedge\dd a^{n-1}\wedge\dd a^{n} \\ & = (-1)^{n-1}a^{0}a^{1}\dd a^{2}\wedge\dots \wedge \dd a^{n} \\ & = 0, } where we used that elements of $\Omega_{u}^{\bullet}$ are in the intersection of the kernel of all adjacent multiplication maps. 
		
		The map $\Phi$ is a morphism of graded algebra since, given $a^{0}\otimes \dots \otimes a^{n}\in \Omega^{n}_{u}$ and $b^{0}\otimes\dots\otimes b^{m}\in\Omega^{m}$ we find \eqa{ \nonumber\Phi(a^{0}\otimes \dots \otimes a^{n})\wedge\Phi(b^{0}\otimes \dots \otimes b^{m}) & = a^{0}\dd a^{1}\wedge\dd a^{n}b^{0}\wedge\dd b^{1}\wedge\dots \wedge b^{m} \\ & = a^{0}\dd a^{1}\wedge\dots\wedge \dd(a^{n}b^{0}) \wedge\dots \wedge\dd b^{m} \\ & \ \ - a^{0}\dd a^{1} \wedge\dd (a^{n-1})a^{n} \wedge\dd b^{0}\wedge\dots \dd b^{m} \\ & = \Phi\tonde{(a^{0}\otimes \dots \otimes a^{n})\wedge_{u}(b^{0}\otimes \dots \otimes b^{m})}\\ & \ \ - \omega\wedge \dd b^{0} \wedge\dots \wedge \dd b^{n} \\ & = \Phi\tonde{(a^{0}\otimes \dots \otimes a^{n})\wedge_{u}(b^{0}\otimes \dots \otimes b^{m})}. }
		
		Finally we have \eqa{\nonumber \Phi^{n+1}\circ\dd_{u}(a^{0}\otimes \dots \otimes a^{n}) & = \Phi^{n+1}\circ \tonde{\sum_{k=0}^{n+1}(-1)^{k}a^{0}\otimes \dots a^{k-1}\otimes 1\otimes a^{k}\otimes \dots \otimes a^{n}} \\ & =\Phi^{n+1}\tonde{1\otimes a^{0}\otimes\dots\otimes a^{n}} \\ & \ \ + \Phi^{n+1}\tonde{\sum_{k=1}^{n+1} (-1)^{k} a^{0}\otimes  \dots \otimes a^{k-1} \otimes 1 \otimes a^{k}\otimes \dots \otimes a^{n}} \\ & = \dd a^{0}\wedge \dd a^{1}\wedge \dots \wedge \dd a^{n}\\ & = \dd \circ \Phi^{n}(a^{0}\otimes \dots \otimes a^{n}),} where we used that $1\in \ker\dd$. This proves that $\Phi$ respects the differentials.  
	Therefore we proven that there exists a surjective morphism of differential graded algebra $\Phi:\Omega^{\bullet}_{u}	\rightarrow\Omega^{\bullet}$. It follows that $(\Omega_{u}^{\bullet}/\ker\Phi,\wedge_{u},\dd_{u})\cong (\Omega^{\bullet},\wedge,\dd)$.  \qed 
	\end{lemma}
	\begin{prop}
		The maximal prolongation of the first order universal differential calculus $(\Gamma_{u},\dd_{u})$ is its tensor algebra, in other words the higher order universal differential calculus $(\Omega_{u}^{\bullet},\wedge_{u},\dd_{u})$. 
		\proof 
		Let $\M{T}A$ be the tensor algebra built from the first order universal differential calculus $\Gamma_{u}$ on $A$, i.e. $$\M{T}A:=\bigoplus_{n\in\B{N}}\M{T}^{n}A=A\oplus \Gamma_{u}\oplus \Gamma_{u}\otimes_{A}\Gamma_{u}\oplus\cdots.$$ Define $S^{\wedge}$ as the ideal of $2-$forms generated by elements $\sum_{i}\dd_{u} a^{i}\otimes_{A} \dd_{u} b^{i}$ such that $\sum_{i}a^{i}\dd_{u} b^{i}=0$, where $\dd_{u}:A\rightarrow \Gamma$ is the usual universal differential.  
		
		To construct the maximal prolongation of $(\Gamma_{u},\dd_{u})$ we quotient out the ideal $S^{\wedge}$. We now show $\M{T}A/S^{\wedge}$ is isomorphic to $\Omega_{u}^{\bullet}$. 
		
		 From the Leibinz rule we see \eqa{\dd (\sum_{i}a^{i}b^{i})=\sum_{i}\dd_{u} a^{i} b^{i} + \sum_{i}a^{i}\dd_{u}b^{i},} which equivalently tells us that $\sum_{i}\dd_{u} a^{i}b^{i}$ is equal to zero since $\sum_{i}a^{i}b^{i}=0$ in the very definition of $\Gamma_{u}$.  
		Accordingly we find \eqa{ \sum_{i}\dd_{u} a^{i}b^{i}& = \sum_{i}a^{i}\otimes b^{i} -  1\otimes a^{i}b^{i}= \sum_{i} a^{i}\otimes b^{i}=0, \\ \sum_{i}a^{i}\dd_{u} b^{i} & = \sum_{i}a^{i}b^{i}\otimes 1-a^{i}\otimes b^{i}= -\sum_{i}a^{i}\otimes b^{i}=0. } Therefore, the ideal $S^{\wedge}$ is generated by $0$, and the quotient $\M{T}A/S^{\wedge}$ is naturally isomorphic to $\Omega_{u}^{\bullet}$ itself. Moreover, by theorem \ref{Maximal prolongation: extension of the differential} we know that the differential $\dd_{u}:\Omega^{\bullet}\rightarrow\Omega^{\bullet}$ is the unique extension of the differential $\dd:A\rightarrow \Gamma$.  \qed 
	\end{prop}
	
	\begin{thm}[\cite{Durdevic-1},Appendix B] Let $(\Omega^{\bullet},\widetilde{\wedge},\widetilde{\dd})$ be any differential calculus on $A$ such that $\Omega^{1}=\Gamma$ and $\widetilde{d}|_{A}=\dd$. There exists a surjective morphism $\Gamma^{\wedge}\rightarrow \Omega^{\bullet}$ of differential graded algebras. In particular, $(\Omega^{\bullet},\widetilde{\wedge},\widetilde{\dd})$ is a quotient of $(\Gamma^{\wedge},\wedge,\dd)$.    \qed

\end{thm}

\chapter{Quantum principal bundles}\label{Quantum principal bundles chap}
This chapter represents the core of this thesis work. The \DJ ur\dj evi\'c theory of quantum principal bundles introduced in \cite{Durdevic-2,Durdevic-1} is presented adopting a modern language and notation. Essentially, a quantum principal bundle, or principal comodule algebra, is understood as a faithfully flat Hopf-Galois extension $B\subseteq A$, where $A$ is a right $H-$comodule algebra, $H$ is a Hopf algebra and $B$ is the subalgebra of coinvariant elements in $A$ under the coaction of $H$. In section \ref{preliminary notions} we introduce a left covariant first order differential calculus on the Hopf algebra $H$, the corresponding subalgebra of left invariant forms and the corresponding right $H-$module structure. We extend to the maximal prolongation of this first order calculus providing a differential calculus over $H$. We discuss the space of left invariant forms of this calculus.  We then fix a bicovariant first order differential calculus over $H$. 

In section \ref{Quantum principal bundles} we provide the definition of quantum principal bundle as a faithfully flat Hopf-Galois extension. We discuss the space of vertical forms, which turns out to be a differential calculus under a suitable choice of wedge product and differential.

In section \ref{Total space calculus} we introduce the notion of complete calculus over a quantum principal bundle, being a differential calculus over the totale space such that the right $H-$coaction lifts  to a morphism of differential graded algebras. We develop a theory for complete differential calculi over right $H$-comodule algebras. A definition of horizontal forms is provided; it is shown that for first order differential calculi there is a short exact sequence of $A-$modules involving horizontal, vertical and total space forms. We show how in general such exact sequence fails for higher order calculi. In subsection \ref{Base space calculus} we discuss differential calculus over the base space, i.e. over the right $H-$coinvariant elements. We show how the pull-back calculus by \cite{beggs-majid} and the one by \DJ ur\dj evi\'c are different in general. However, we want to stress that in all explicit examples in this thesis the pull-back calculus coincides with the base space calculus. In subsection \ref{equivalence with bm} we investigate the relation between quantum principal bundle of \cite{BrzMjd,beggs-majid} and \DJ ur\dj evi\'c, providing a condition under which the two definition coincide. In particular, every quantum principal bundle in the sense of \cite{BrzMjd} on a faithfully flat Hopf-Galois extension is automatically first order complete in the sense of \DJ ur\dj evi\'c. 

In section \ref{Examples} we provide some non-trivial examples of quantum principal bundles in the \DJ ur\dj evi\'c theory. The non-commutative torus, the quantum Hopf fibration and crossed product calculi turn out to be complete differential calculi in our previous definition. Moreover, in these three examples, the base space calculus is generated by the base itself.
 
\section{ Preliminary notions}\label{preliminary notions}
Let $H$ be a Hopf algebra $(H,\mu,\eta,\Delta,\epsilon,S)$. Recall the maximal prolongation $\Gamma^{\wedge}$ introduced in section \ref{maximal prolongation} obtained as a quotient of the tensor bundle algebra on a first order calculus $\Gamma$ over $H$ by the ideal $S^{\wedge}$. For now let $\Gamma$ be a left-covariant first order differential calculus over $H$; let ${}_{\Gamma}\Delta:\Gamma\rightarrow H\otimes \Gamma$ be the corresponding left coaction of $H$.  Denote by $\Lambda^{1}$ the space of left-invariant elements of $\Gamma$, namely $$\Lambda^{1}=\{\theta \in \Gamma : {}_{\Gamma}\Delta(\theta)=1\otimes \theta\}.$$ Consider the quantum Maurer-Cartan form $\varpi:H^{+}\rightarrow \Lambda^{1}$ given by $\varpi(a)=S(a_{1})\dd a_{2}$ where $H^{+}=\ker\epsilon$, and let furthermore $I\subseteq \ker(\epsilon)=H^{+}$ be the right $H-$ideal which corresponds canonically to $\Gamma$ under the identification provided by Theorem \ref{Woronowicz classification section}. We have $\ker\varpi = I$, and because of this there exists a natural isomorphism $\Lambda^{1} \cong H^{+}/I.$ This isomorphism induces a right $H-$module structure on $\Lambda^{1}$ which will be denoted as $\leftharpoonup $. 
Explicitely, the action is given as  \eqa{\nonumber\varpi (a)\leftharpoonup b & =\varpi (ab-\epsilon(a)b) \\ & = \varpi(ab)-\epsilon(a)\varpi(b) \\ & = S(a_{1}b_{1})\dd(a_{2}b_{2})- \epsilon(a)S(b_{1})\dd b_{2} \\ & = S(b_{1})S(a_{1})(\dd(a_{2})b_{2} + a_{2}\dd b_{2})- \epsilon(a)S(b_{1})\dd b_{2} \\ & = S(b_{1}) \varpi(a)b_{2} + S(b_{1})\epsilon(a)\dd b_{2} - \epsilon(a) S(b_{1})\dd b_{2} \\ & =S(b_{1})\varpi(a)b_{2},}  
for each $a,b\in H$.  

\begin{rmk}
	
The Maurer-Cartan form satisfies the following relation among product of elements in $H^{+}$. \eqa{\label{maurer cartan along product} \varpi(ab) & = S(a_{1}b_{1}) \dd (a_{2}b_{2}) \\ & = S(b_{1})S(a_{1}) (\dd a_{2}b_{2} + a_{2} \dd b_{2}) \\ & = S(b_{1})\varpi(a)b_{2} + S(b_{1})\epsilon(a)\dd b_{2} \\ & = \varpi(a)\leftharpoonup b + \epsilon(a) \varpi(b),}  for every $a,b\in H^{+}$. \qed \end{rmk}

\begin{rmk}By left covariance of the calculus the maps $\Delta$ and ${}_{\Gamma}\Delta$ admit common extensions to homomorphisms $${}_{\Gamma^{\wedge}}\Delta:\Gamma^{\wedge}\rightarrow H\otimes \Gamma^{\wedge}\quad \text{and} \quad{}_{\Gamma^{\otimes}}\Delta:\Gamma^{\otimes} \rightarrow H\otimes \Gamma^{\otimes}$$ according to proposition \ref{equivalent definitions of covariance}. \qed \end{rmk}Let $\Lambda^{\bullet}$ be the differential graded subalgebra of left-invariant elements of $\Gamma^{\wedge}$, namely \eqa{ \Lambda^{\bullet}=\bigoplus_{k\geq 0} \graffe{ \omega \in \Omega^{k}(H) \ : {}_{\Gamma^{\wedge}}\Delta(\omega)=1\otimes \omega }.}   \begin{prop}
  	The map $\pi_{inv}:\Gamma^{\wedge}\rightarrow \Gamma^{\wedge}$ defined by $\omega\mapsto S(\omega_{-1})\omega_{0}$ maps onto the space of invariant forms $\Lambda^{\bullet}$. 
  	\proof Consider $\omega\in \Gamma^{\wedge}$, which can be written as $\omega=a^{0}\dd a^{1}\wedge\dots \wedge \dd a^{k}$. We find that \eqa{\nonumber\pi_{inv}(\omega) & = S(\omega_{-1})\omega_{0} \\ & = \mu\circ (S\otimes \id)\circ {}_{\Gamma^{\wedge}}\Delta(\omega) \\ & =  \mu\circ (S\otimes \id)\circ {}_{\Gamma^{\wedge}}\Delta(a^{0}\dd a^{1}\wedge\cdots \dd a^{k}). }   and so \eqa{\nonumber\pi_{inv}(\omega) & = \mu \circ ( S\otimes \id) \circ (a^{0}_{1}\otimes a^{0}_{2})(\wedge\otimes \wedge)(a^{1}_{1}\otimes \dd  a^{1}_{2}) (\wedge\otimes \wedge)\cdots (\wedge\otimes \wedge)(a^{k}_{1}\otimes \dd a^{k}_{2}) \\ & = \epsilon(a^{0})S(a^{1}_{1}\cdots a^{k}_{1})\dd a^{1}_{2}\wedge\cdots \wedge \dd a^{k}_{2}.} Applying ${}_{\Gamma^{\wedge}}\Delta$ we obtain  \eqa{\nonumber {}_{\Gamma^{\wedge}}\Delta(\pi_{inv}(\omega)) & = {}_{\Gamma^{\wedge}}\Delta(\epsilon(a^{0})S(a^{1}_{1}\cdots a^{k}_{1})\dd a^{1}_{2}\wedge\cdots \wedge \dd a^{k}_{2}) \\ & = \epsilon(a^{0})\Delta(S(a^{1}_{1}\cdots a^{k}_{1}))(\id\otimes \dd) \Delta(a^{2}_{2})(\wedge\otimes \wedge)\cdots (\wedge\otimes \wedge)(\id\otimes \dd)\Delta(a^{k}_{2}) \\ & = \epsilon(a^{0})(S(a_{2}^{1}\cdots a_{2}^{k})a^{1}_{3}\cdots a^{k}_{3}\otimes S(a_{1}^{1}\cdots a_{1}^{k})\dd a^{1}_{4}\wedge\cdots \wedge \dd a^{k}_{4}) \\ & = 1\otimes \epsilon(a^{0})S(a_{1}^{1}\cdots a_{1}^{k})\dd a^{1}_{2}\wedge\cdots \wedge \dd a^{k}_{2} \\ & = 1\otimes S(a^{0}_{1}\cdots a^{k}_{1})a^{0}_{2}\dd a^{1}_{2}\wedge\cdots \wedge \dd a^{k}_{2}. } 
  	Clearly, $\pi_{inv}$ is surjective, since for $\omega\in\Lambda^\bullet$ we obtain $\pi_{inv}(\omega)=\omega.$ \qed 
  \end{prop}
  \begin{prop}[\cite{beggs-majid},Proposition 2.31] 
Let $(\Gamma,\dd)$ be a left covariant first order differential calculus on a Hopf algebra $H$. Higher order coinvariant forms correspond to tensor products of coinvariant 1-forms quotiented by some relations. Precisely \eqa{\Lambda^{\bullet}=(\Lambda^{1})^{\otimes}/S^{\wedge}_{inv},} where $S^{\wedge}_{inv}:=\langle \varpi(\pi_{\epsilon}(a_{1}))\otimes \varpi(\pi_{\epsilon}(a_{2})) \ | \ a \in A\rangle$, where $\pi_{\epsilon}:H\rightarrow H^{+}$ is defined by $\pi_{\epsilon}(h)=h-\epsilon(h)1$. \qed  \end{prop} 

 \begin{prop} The right $H-$module structure $\leftharpoonup$ on $\Lambda^{1}$  can be  extended uniquely to $\Lambda^{\bullet}$, i.e. $\Lambda^{\bullet}$ is a graded right $H-$module algebra, with \eqa{ \label{hook action on the wedge} 1\leftharpoonup a& =\epsilon(a)1,\\  (\theta \wedge \eta)\leftharpoonup a& = (\theta \leftharpoonup a_{1})\wedge (\eta\leftharpoonup a_{2}),}for each $\theta,\eta$ either in $\Lambda^{\bullet}$, and every $a\in H$.
 	 Explicitly, $\theta \leftharpoonup a =S(a_{1})\theta a_{2}$.
 	\proof  First of all, for any $\theta\in \Lambda^{\bullet}$ and $a,b\in H$, we have \eqa{\nonumber(\theta \leftharpoonup a)\leftharpoonup b& = S(b_{1}) (\theta\leftharpoonup a)b_{2}\\ & = S(b_{1})S(a_{1}) \theta a_{2}b_{2} \\ & = S(a_{1}b_{1})\theta a_{2}b_{2} \\ & = \theta \leftharpoonup ab,  } i.e. $\leftharpoonup$ is a right $H-$action. 
 	
 	Moreover, let $\omega,\eta\in \Lambda^{\bullet}$. Then \eqa{ \nonumber(\theta\leftharpoonup h_{1})\wedge (\eta\leftharpoonup h_{2}) & = (S(h_{11})\theta h_{12})\wedge (S(h_{21})\eta h_{22}) \\ & = S(h_{1})\theta h_{2}\wedge S(h_{3})\eta h_{4} \\ & = S(h_{1})\theta \wedge h_{2}S(h_{3})\eta h_{4} \\ & = S(h_{1})\theta\wedge \eta h_{2} \\ & = (\theta\wedge\eta)\leftharpoonup h.}
  \qed 
 \end{prop}
 \begin{prop}
 	The algebra $\Lambda^{\bullet}\subseteq \Gamma^{\wedge}$ is $\dd-$invariant. \label{d-invariance of coinvariant forms}
 	\proof  Let $\dd:\Gamma^{\wedge}\rightarrow\Gamma^{\wedge}$ be the differential extending $\dd:H\rightarrow \Gamma$ constructed in \ref{Maximal prolongation: extension of the differential}. Let $\omega\in \Lambda^{\bullet}$ be a coinvariant form. By definition ${}_{\Gamma^{\wedge}}\Delta(\omega)=1\otimes \omega$. This means $a_{-1}^{0}\cdots a_{-1}^{k}=1$.  Since $\dd$ is left $H-$colinear, it follows that \eqa{ {}_{\Gamma^{\wedge}}\circ \Delta(\dd\omega) & =(\id\otimes \dd)\circ {}_{\Gamma^{\wedge}}\Delta(\omega)\\ & =1\otimes \dd\omega.} 
 \end{prop}
 \begin{prop}
 	The "Maurer-Cartan equation" $\dd \varpi(\pi_{\epsilon}(a)) =- \varpi(\pi_{\epsilon}(a_{1}))\varpi(\pi_{\epsilon}(a_{2}))$ holds.
 \proof By a direct calculation we find \eqa{ \nonumber \varpi(\pi_{\epsilon}(a_{1}))\wedge \varpi(\pi_{\epsilon}(a_{2})) & = S(\pi_{\epsilon}(a_{11}))\dd \pi_{\epsilon}(a_{12}) \wedge S(\pi_{\epsilon}(a_{21}))\dd \pi_{\epsilon}(a_{22}) \\ & = S(\pi_{\epsilon}(a_{1}))\dd \pi_{\epsilon}(a_{2})\wedge S(\pi_{\epsilon}(a_{3}))\dd \pi_{\epsilon}(a_{4}) \\ & = S(\pi_{\epsilon}(a_{1}))\dd \pi_{\epsilon}(a_{2})S(\pi_{\epsilon}(a_{3}))\wedge \dd \pi_{\epsilon}(a_{4}) \\ & = S(\pi_{\epsilon}(a_{1})) (\dd(\pi_{\epsilon}(a_{2})S(\pi_{\epsilon}(a_{3})))-\pi_{\epsilon}(a_{2})\dd S(\pi_{\epsilon}(a_{3})))\wedge \dd \pi_{\epsilon}(a_{4}) \\ & = S(\pi_{\epsilon}(a_{1}))\dd(1)\wedge \dd \pi_{\epsilon}(a_{4}) - S(\pi_{\epsilon}(a_{1}))\pi_{\epsilon}(a_{2}) \dd S(\pi_{\epsilon}(a_{3})) \wedge \dd \pi_{\epsilon}(a_{4}) \\ & = - S(\pi_{\epsilon}(a_{1}))\pi_{\epsilon}(a_{2})\dd(S(\pi_{\epsilon}(a_{3})))\wedge \dd \pi_{\epsilon}(a_{4}) \\ &=  - \epsilon(a_{1})\dd (S(\pi_{\epsilon}(a_{11}))) \wedge \dd \pi_{\epsilon}(a_{12}) \\ & = -\dd(S(\epsilon(a_{1})a_{2})))\wedge \dd \pi_{\epsilon}(a_{3}) \\ & = - \dd(S(\pi_{\epsilon}(a_{1})))\wedge \dd \pi_{\epsilon}(a_{2}) \\ & = -\dd (\varpi(\pi_{\epsilon}(a))),} since $\dd(\varpi(\pi_{\epsilon}(a)))=\dd(S(\pi_{\epsilon}(a_{1}))\dd \pi_{\epsilon}(a_{2})) =\dd(S(\pi_{\epsilon}(a_{1})))\wedge \dd \pi_{\epsilon}(a_{2})$.  	\qed 
 \end{prop}
 \begin{lemma} The differential of the right $H-$module structure on coinvariant differential forms reads
 	\eqa{ \dd(\theta\leftharpoonup a) = (\dd(\theta)\leftharpoonup a)  - \varpi(\pi_{\epsilon}(a_{1}))(\theta\leftharpoonup a_{2})+ (-1)^{|\theta|}(\theta \leftharpoonup a_{1})\varpi (\pi_{\epsilon}(a_{2})),} for $a\in H$ and 
 	\proof The left hand side reads \eqa{ \nonumber\dd(\theta\leftharpoonup a) & = \dd( S(a_{1})\theta a_{2}) \\ & = \dd(S(a_{1}))\theta a_{2} + S(a_{1}) \dd\theta a_{2} + (-1)^{|\theta|} S(a_{1})\theta \dd a_{2}.} 
 	Terms on the right hand side read  
 	\eqa{\nonumber \dd\theta\leftharpoonup a & =S(a_{1})\dd\theta a_{2}; \\ \varpi(\pi_{\epsilon}(a_{1}))(\theta\leftharpoonup a_{2}) & = S(\pi_{\epsilon}(a_{11}))\dd(\pi_{\epsilon}(a_{12}))S(a_{21})\theta a_{22} \\ & = S(\pi_{\epsilon}(a_{1}))\dd \pi_{\epsilon}(a_{2}) S(a_{3}) \theta a_{4} \\ & = S(\pi_{\epsilon}(a_{1}))(\dd(\pi_{\epsilon}(a_{2})S(a_{3}))-\pi_{\epsilon}(a_{2})\dd S(a_{3}))\theta a_{4} \\ & =  -S(\pi_{\epsilon}(a_{1}))\pi_{\epsilon}(a_{2})\dd S(a_{3})\theta a_{4} \\ & = - \dd S(\pi_{\epsilon}(a_{1})) \theta a_{2}; \\ (-1)^{|\theta|}(\theta\leftharpoonup a_{1})\varpi(\pi_{\epsilon}(a_{2})) & = (-1)^{|\theta|}S(a_{11})\theta a_{12}S(\pi_{\epsilon}(a_{21}))\dd \pi_{\epsilon}(a_{22}) \\ & = (-1)^{|\theta|}S(a_{1})\theta a_{2} S(\pi_{\epsilon}(a_{3}))\dd \pi_{\epsilon}(a_{4}) \\ & = (-1)^{|\theta|} S(a_{1}) \theta \dd \pi_{\epsilon}(a_{2}),} therefore \eqa{ \nonumber \dd(\theta)\leftharpoonup a  - \varpi(\pi_{\epsilon}(a_{1}))(\theta\leftharpoonup a_{2})+ (-1)^{|\theta|}(\theta \leftharpoonup a_{1})\varpi (\pi_{\epsilon}(a_{2})) & = S(\pi_{\epsilon}(a_{1}))\dd\theta a_{2} \\ & + \dd S(a_{1}) \theta a_{2} \\ & +  (-1)^{|\theta|} S(a_{1}) \theta \dd \pi_{\epsilon}(a_{2}).}\qed  
 \end{lemma}
  In the following we fix $\Gamma$ a bicovariant first order differential calculus over $H$ and consider its maximal prolongation to a higher order calculus $\Gamma^{\wedge}$.    \begin{prop}
	Let $H$ be a Hopf algebra and let $(\Omega^{\bullet}(H),\wedge,\dd)$ be a bicovariant differential calculus on $H$. Then $\Omega^{\bullet}(H)$ is a graded Hopf algebra with coproduct $\Delta^{\bullet}:\Omega^{k}(H)\rightarrow \Omega^{m}(H)\otimes \Omega^{n}(H), $ where $m+n=k$, sending $\omega\mapsto \omega_{(1)}\otimes \omega_{(2)}=:\omega_{-1}\otimes \omega_{0}+\omega_{0}\otimes \omega_{1}=({}_{\Gamma^{\wedge}}\Delta+\Delta_{\Gamma^{\wedge}})(\omega)$, counit $\epsilon^{\bullet}:\Omega^{\bullet}(H)\rightarrow \B{k}$, with $\epsilon(\omega)=0$ for any $\omega$ such that $|\omega|>1$, and antipode $S^{\bullet}:\Omega^{\bullet}\rightarrow\Omega^{\bullet}$ defined as $S^{\bullet}(\omega)=-S(\omega_{-1})\omega_{0}S(\omega_{1})$.  
	\proof 
	We start by checking coassociativity of $\Delta^{\bullet}$. Let $\omega\in \Omega^{k}(H)$. We have \eqa{\nonumber(\Delta^{\bullet}\otimes \id)\circ \Delta^{\bullet}(\omega) & = (\Delta^{\bullet}\otimes \id)(\omega_{-1}\otimes \omega_{0}+ \omega_{0}\otimes \omega_{1}) \\ & = (\omega_{-1})_{1}\otimes (\omega_{-1})_{2} \otimes \omega_{0} \\ & + (\omega_{0})_{-1}\otimes(\omega_{0})_{0}\otimes \omega_{1}\\ &  + (\omega_{0})_{0}\otimes (\omega_{0})_{1} \otimes \omega_{1}\\ & = \omega_{-2}\otimes \omega_{-1}\otimes \omega_{0} + \omega_{-1}\otimes \omega_{0}\otimes \omega_{1} + \omega_{0}\otimes \omega_{1}\otimes \omega_{2} \\ & = (\id\otimes \Delta^{\bullet})\circ \Delta^{\bullet}(\omega).} Next we have \eqa{ \nonumber(\epsilon\otimes \id)\circ \Delta^{\bullet}(\omega) & = (\epsilon\otimes \id) ( \omega_{-1}\otimes \omega_{0}+ \omega_{0}\otimes \omega_{-1}) \\ & = \epsilon(\omega_{-1})\otimes \omega_{0} + \epsilon(\omega_{0})\otimes \omega_{1} \\ & = \epsilon(\omega_{-1})\omega_{0}\otimes 1 \\ & = \omega\otimes 1 \\ & = \omega.} Finally \eqa{\nonumber\omega_{(1)}S(\omega_{(2)}) & = \omega_{-1}S(\omega_{0}) +\omega_{0}S(\omega_{1}) \\ & = -\omega_{-1}S((\omega_{0})_{-1})(\omega_{0})_{0}S((\omega_{0})_{1}) + \omega_{0}S(\omega_{1}) \\ & = -\omega_{-2}S(\omega_{-1})\omega_{0}S(\omega_{1}) + \omega_{0}S(\omega_{1}) \\ & = -\omega_{0}S(\omega_{1}) + \omega_{0}S(\omega_{1}) \\ & = 0\\ & = \epsilon(\omega)1 . }  Similarly for $S(\omega_{(1)})\omega_{(2)}$.  \qed 
\end{prop}

\begin{rmk}\label{delta bullet not a morphism of DGAs in general} 
	In \cite{Durdevic-2} it is stated that the lifted graded coproduct  $\Delta^{\bullet}={}_{\Gamma^{\wedge}}\Delta+\Delta_{\Gamma^{\wedge}}:\Gamma^{\wedge}\rightarrow \Gamma^{\wedge}\otimes \Gamma^{\wedge}$ is a morphism of differential graded algebras. However, it turns out that  $\Delta^{\bullet}$ is a morphism of graded algebras, but is not a morphism of differential graded algebras, since it is not compatible with the differential in general. Consider for instance a $2-$form $\omega\in \Omega^{2}(H)$. Then, if we assume $\Delta^{\bullet}$ to be a morphism of differential graded algebras we get \eqa{\nonumber \Delta^{\bullet}( \omega)&  =\Delta^{\bullet}(a\dd b\wedge \dd c) \\ & = \Delta(a)\Delta^{1}(\dd a) \Delta^{1}(\dd b) \\ & = (a_{1}\otimes a_{2})(\dd\otimes \id + \id \otimes \dd)(b_{1}\otimes b_{2})(\dd\otimes \id + \id \otimes \dd)(c_{1}\otimes c_{2}) \\ & = (a_{1}\otimes a_{2})(-b_{1}\dd c_{1}\otimes \dd b_{2}c_{2} + \dd b_{1} \wedge \dd c_{1}\otimes b_{2}c_{2} + \dd b_{1}c_{1}\otimes b_{2}\dd c_{2} + b_{1} c_{1}\otimes \dd b_{2}\wedge \dd c_{2}). \nonumber }On the other hand, according to \DJ ur\dj evi\'c's definition we would get \eqa{\nonumber \Delta^{\bullet}(\dd \omega) & = ({}_{\Gamma^{\wedge}}\Delta+\Delta_{\Gamma^{\wedge}})(a\dd b \wedge \dd b) \\ & = {}_{\Gamma^{\wedge}}\Delta(a\dd b \wedge \dd c) + \Delta_{\Gamma^{\wedge}}(a \dd b \wedge \dd c) \\ & = a_{1}b_{1}c_{1}\otimes a_{2}\dd b_{2}\wedge \dd c_{2} + a_{1}\dd b_{1}\wedge \dd c_{1}\otimes a_{2}b_{2}c_{2}.} Since two terms are missing, we conclude $\Delta^{\bullet}={}_{\Gamma^{\wedge}}\Delta+\Delta_{\Gamma^{\wedge}}$ is not a morphism of differential graded algebras. \qed 
\end{rmk}

 Thus, we choose to restrict our attention to bicovariant differential calculi over Hopf algebras such that the coproduct lifts uniquely to a morphism $\Delta^{\bullet}:\Gamma^{\wedge}\rightarrow\Gamma^{\wedge}\otimes \Gamma^{\wedge}$ of differential graded algebras. \begin{notation}
 	
 The action of $\Delta^{\bullet}:\Gamma^{\wedge}\rightarrow\Gamma^{\wedge}\otimes \Gamma^{\wedge}$ on a given form $\omega\in \Gamma^{\wedge}$ is denoted as $\Delta^{\bullet}(\omega)=\omega_{[1]}\otimes \omega_{[2]}.$  \end{notation}
 
\section{Quantum principal bundles}\label{Quantum principal bundles}
 We now introduce the definition of quantum principal bundle in \DJ ur\dj evi\'c' theory as a faithfully flat Hopf-Galois extension.  
\begin{defn}
	Let $B$ be a ring and let $X,Y,Z$ be $B-$modules. A $B-$module $A$ such that  $$0\longrightarrow X\longrightarrow Y \longrightarrow Z \longrightarrow 0$$ is short exact if and only if $$0\longrightarrow X\otimes_{B} A \longrightarrow Y\otimes_{B} A\longrightarrow Z\otimes_{B} A \longrightarrow 0 $$ is short exact, is called a faithfully flat left $B-$module.  
\end{defn} 
 
\begin{defn}
	A quantum principal bundle, or principal comodule algebra, is a faithfully flat Hopf-Galois extension $B\subseteq A$, where $A$ is a right $H-$comodule algebra and $B:=A^{coH}$ is the space of coinvariant forms in $A$ under the coaction of $H$.\end{defn} We consider differential calculi over such quantum principal bundles. 
	
Recall that in classical differential geometry the vertical bundle and the horizontal bundle are vector bundles associated with a smooth fiber bundle. Specifically, given a smooth fiber bundle $\pi: E \to B$, the vertical bundle $VE$ and horizontal bundle $HE$ are subbundles of the tangent bundle satisfying $VE \oplus HE \cong TE$, the horizontal bundle choice being fixed by the connection. This implies that, over each point $e \in E$, the fibers $V_e E$ and $H_e E$ form complementary subspaces of the tangent space $T_e E$. The vertical bundle comprises all vectors tangent to the fibers.
	\subsection{Vertical forms}

	The following definition provides the non-commutative analogue of vertical forms of classical differential geometry. 
\begin{defn} 
	Let $B\subseteq A$ be a quantum principal bundle. Let $\Gamma$ be a bicovariant first order differential calculus over $H$,  $\Gamma^\wedge$ the maximal prolongation and $\Lambda^\bullet$ the space of left coinvariant forms on $H$. The graded vector space  $\ver^{\bullet}(A)=A\otimes \Lambda^{\bullet}$ defines the space of vertical forms. 
\end{defn} The space of vertical forms is a differential calculus thanks to the following lemma.  
\begin{lemma}\label{vertical forms as a DC} 
	The space of vertical forms $(\ver^{\bullet}(A),\wedge_{\ver},\dd_{\ver})$ is a differential calculus over $A$ by the assignments \eqa{ (a\otimes \theta)\wedge_{\ver}(b\otimes \eta)=a b_{0}\otimes (\theta \leftharpoonup b_{1})\wedge\eta, } 
	\eqa{\label{hermitean differential} \dd_{\ver}(a\otimes \theta)=a\otimes \dd\theta +a_{0}\otimes \varpi(\pi_\epsilon(a_1))\wedge \theta,} for every $\theta,\eta\in \Lambda^{\bullet}$ and $a,b\in A$. 
	\proof During the proof we will call write $(\wedge,\dd)$ for $(\wedge_{\ver},\dd_{\ver})$ to shorten the notation. We shall prove that $\wedge$ is an associative product that is unital, and moreover the differential satisfies the graded Leibniz rule and squares to zero. 
	The product of vertical forms is associative, indeed \eqa{\nonumber((a\otimes \theta)\wedge(b\otimes \eta))\wedge(c\otimes \omega) & = (a_{0}b_{0}\otimes (\theta\leftharpoonup b_{1})\wedge\eta)\wedge\omega \\ & = ab_{0}c_{0}\otimes((\theta\leftharpoonup b_{1})\wedge \eta)\leftharpoonup c_{1})\wedge \omega \\ & = a b_{0}c_{0}\otimes ((\theta \leftharpoonup b_{1})\leftharpoonup c_{11})\wedge (\eta \leftharpoonup c_{12})\wedge\omega \\  &  = a b_{0}c_{0}\wedge(\theta \leftharpoonup b_{1}c_{1})\wedge (\eta \leftharpoonup c_{2})\wedge \omega  \\ & = a b_{0} c_{00} \otimes (\omega\leftharpoonup b_{1}c_{01})\wedge (\eta \leftharpoonup c_{1}) \wedge\omega  \\ & = (a\otimes \theta)\wedge (b c_{0} \otimes (\eta\leftharpoonup c_{1})\wedge \omega \\ & = (a\otimes \theta)\wedge ((b\otimes \eta)\wedge (c\otimes \omega)).} Next we have $\ver(P)$ is a unital algebra with unit $1\otimes 1$, indeed   
	 \eqa{\nonumber (1\otimes 1)\wedge (a\otimes \theta) & = a_{0} \otimes (1\leftharpoonup a_{1})\wedge \theta \\ & = a_{0} \otimes \epsilon(a_{1}) \theta \\ & = a_{0}\epsilon(a_{1})\otimes \theta \\ & = a\otimes \theta} and \eqa{\nonumber(a\otimes \theta)\wedge (1\otimes 1) & = a\otimes (\theta\leftharpoonup 1) \\ & = a\otimes S(1)\theta \\ & = a\otimes \theta. } Now we use the Maurer-Cartan equation $\dd\varpi(\pi_{\epsilon}(h)) = -\varpi(\pi_{\epsilon}(h_{1}))\wedge\varpi(\pi_{\epsilon}(h_{2}))$ and the property $$ \dd(\theta \leftharpoonup h)=(\dd(\theta)\leftharpoonup h) -\varpi(\pi_{\epsilon}(h_{1}))\wedge(\theta\leftharpoonup h_{2})+(-1)^{|\theta|}(\theta \leftharpoonup h_{1})\wedge\varpi(\pi_{\epsilon}(h_{2}))$$ to show we have a differential on the algebra  $\ver(A)$ defined by \eqref{hermitean differential}. 
	  We find \eqa{ \label{differential over vertical} \dd((a\otimes\theta)\wedge(b\otimes \eta)) & =\dd(ab_{0}\otimes (\theta\leftharpoonup b_{1})\wedge \eta) \\ & =ab_{0}\otimes \dd((\theta\leftharpoonup b_{1})\wedge\eta)\\ & + a_{0}b_{0}\otimes \varpi(\pi_{\epsilon}(a_{1}b_{1}))\wedge(\theta\leftharpoonup b_{2})\wedge \eta \\ & =  ab_{0}\otimes (\dd(\theta\leftharpoonup b_{1})\wedge\eta) \\ & + (-1)^{|\theta|}ab_{0}\otimes [(\theta \leftharpoonup b_{1}) \wedge \dd \eta] \\ & + a_{0}b_{0}\otimes \varpi (\pi_{\epsilon}(a_{1}b_{1})) \wedge(\theta\leftharpoonup b_{2})\wedge\eta \\  & = ab_{0} \otimes (\dd \theta \leftharpoonup b_{1})\wedge \eta \\  & - ab_{0}\otimes \varpi(\pi_{\epsilon}(b_{1}))\wedge (\theta \leftharpoonup b_{2})\wedge \eta  \\  & + (-1)^{|\theta|}ab_{0}\otimes (\theta\leftharpoonup b_{1}) \wedge \varpi(\pi_{\epsilon}(b_{2}))\wedge \eta \\ & + (-1)^{|\theta|}ab_{0}\otimes (\theta\leftharpoonup b_{1})\wedge \dd \eta \\ & + a_{0}b_{0}\otimes \varpi(\pi_{\epsilon}(a_{1}b_{1}))\wedge(\theta\leftharpoonup b_{2})\wedge \eta. }  
	  	  
On the other hand
	  \eqa{\nonumber\dd(a\otimes\theta) \wedge (b\otimes \eta)  & + (-1)^{|\theta|}(a\otimes \theta)\wedge \dd(b\otimes \eta) \\  & = (a\otimes \dd\theta+a_{0}\otimes \varpi(\pi_{\epsilon}(a_{1}))\wedge\theta)\wedge(b\otimes \eta) \\ & + (-1)^{|\theta|}(a\otimes \theta)\wedge (b\otimes \dd \eta+b_{0}\otimes \varpi(\pi_{\epsilon}(b_{1})\wedge \eta)) \\ & = ab_{0}\otimes (\dd\theta\leftharpoonup b_{2})\wedge\eta \\  & + (-1)^{|\theta|}a b_{0}\otimes (\theta\leftharpoonup b_{1})\wedge(\varpi (\pi_{\epsilon}(b_{2}))\wedge\eta)\\  & +(-1)^{|\theta|}ab_{0}\otimes (\theta\leftharpoonup b_{1})\wedge\dd\eta \\  & + a_{0}b_{0}\otimes ((\varpi(\pi_{\epsilon}(a_{1}))\wedge\theta) \leftharpoonup b_{1})\wedge\eta . }  
	  We work out the last term of this expression as following. \eqa{a_{0}b_{0}  \otimes ((\varpi(\pi_{\epsilon}(a_{1}))& \wedge\theta) \leftharpoonup b_{1})\wedge\eta  \\ \text{Equation}\  \eqref{hook action on the wedge} & = a_{0}b_{0}[(\varpi(\pi_{\epsilon}(a_{1}))\leftharpoonup b_{1})\wedge (\theta\leftharpoonup b_{2})]\wedge \eta   \\ & = a_{0}b_{0} \otimes [\varpi(\pi_{\epsilon}(a_{1}b_{1}))\wedge(\theta\leftharpoonup b_{2})- \epsilon(a_{1})\varpi(\pi_{\epsilon}(b_{1}))\wedge(\theta\leftharpoonup b_{2})]\wedge \eta \\ \text{Equation} \ \eqref{maurer cartan along product} & = a_{0}b_{0}\otimes \varpi(\pi_{\epsilon}(a_{1}b_{1})) \wedge (\theta\leftharpoonup b_{2})\wedge \eta \\ & - ab_{0}\otimes \varpi(\pi_{\epsilon}(b_{1}))\wedge (\theta\leftharpoonup b_{2})\wedge \eta. \nonumber  } Therefore we have 
	    \eqa{\nonumber\dd(a\otimes\theta) \wedge (b\otimes \eta)  & + (-1)^{|\theta|}(a\otimes \theta)\wedge \dd(b\otimes \eta) \\ & = ab_{0} \otimes (\dd \theta \leftharpoonup b_{1})\wedge \eta \\  & - ab_{0}\otimes \varpi(\pi_{\epsilon}(b_{1}))\wedge (\theta \leftharpoonup b_{2})\wedge \eta  \\  & + (-1)^{|\theta|}ab_{0}\otimes (\theta\leftharpoonup b_{1}) \wedge \varpi(\pi_{\epsilon}(b_{2}))\wedge \eta \\ & + (-1)^{|\theta|}ab_{0}\otimes (\theta\leftharpoonup b_{1})\wedge \dd \eta \\ & + a_{0}b_{0}\otimes \varpi(\pi_{\epsilon}(a_{1}b_{1}))\wedge(\theta\leftharpoonup b_{2})\wedge \eta,} coinciding with Equation \eqref{differential over vertical}.

	    Moreover 
	  \eqa{ \dd^{2}(a\otimes \theta) & =   \dd(a\otimes \dd \theta + a_{0}\otimes \varpi(a_{1})\wedge \theta) \\ & = \dd(a\otimes \dd \theta) + \dd(a_{0}\otimes \varpi(a_{1})\wedge \theta)\\ & = a\otimes \dd^{2}\theta + a_{0}\otimes \varpi(a_{1})\wedge \dd \theta + a_{0} \otimes \dd(\varpi(a_{1})\wedge\theta) +a_{00}\otimes \varpi(a_{01})\wedge (\varpi(a_{1})\wedge\theta) \\ & = 0,} by associativity of the $\wedge$ product and the Maurer-Cartan equation.
	
	 To show that the surjectivity axiom holds start by  considering $\ver^{1}(A)=A\otimes \Lambda^{1}$. It follows that $\ver^{1}(A)=A\dd_{\ver}A$, indeed \eqa{\dd(a\otimes 1) & = a\otimes \dd(1) + a_{0}\otimes \varpi(\pi_{\epsilon}(a_{1})) =a_{0}\otimes \varpi(\pi_{\epsilon}(a_{1})),} but $\varpi(\pi_{\epsilon}(a_{1}))$ maps surjectively onto $\Lambda^{1}$. Thus $a_{0}\otimes \varpi(a_{1})$ is  an element of $A\otimes \Lambda^{1}$, and by a suitable multiplication from left by elements of $A$ we span $\ver^{1}(A)$. 
	 
	  Next we show that $\ver^{i+1}(A)$ is linearly generated by $\ver^{i}(A)$. Let $(a\otimes \theta)\in \ver^{i}(A)$. Differentiating we get $$\dd(a\otimes \theta)=a\otimes \dd\theta + a_{0}\otimes \varpi(a_{1})\wedge \theta $$ by definition of $\dd_{\ver}$. Since the algebra $\Lambda^{\bullet}$ is $\dd$-invariant according to proposition \ref{d-invariance of coinvariant forms} we conclude $\dd\theta\in \Lambda^{k+1}$. \qed 
\end{lemma}
The next lemma characterises the lift of the right $H-$coaction on $A$ to the space of vertical forms. We remind $\Gamma^{\wedge}$ is a graded Hopf algebra with graded coproduct $\Delta^{\bullet}:\Gamma^{\wedge}\rightarrow\Gamma^{\wedge}\otimes \Gamma^{\wedge}$, and that we denote $\Delta^{\bullet}(\omega)=\omega_{[1]}\otimes \omega_{[2]}$. 

\subsection{Total space calculus} \label{Total space calculus} Now we introduce the notion of total calculus on a quantum principal bundle. Fix the usual bicoviariant first order differential calculus $\Gamma$ on $H$. Denote the corresponding maximal prolongation  by $\Omega^{\bullet}(H)$, for convenience.  
\begin{defn}
	A complete differential calculus over the quantum principal bundle $B=A^{coH}\subseteq A$ is a differential calculus $\Omega^{\bullet}(A)$ such that the coaction $\Delta_{A}:A\rightarrow A\otimes H$ extends to a morphism 
	
	$$ \Delta_{A}^{\wedge} : \Omega^{\bullet}(A)\rightarrow \Omega^{\bullet}_{A} \otimes \Omega^{\bullet}(H)$$ of differential graded algebras.  
\end{defn}
From now on we will focus on complete differential calculi. We still continue to restrict our attention to those calculi where $\Delta^{\bullet}:\Omega^{\bullet}(H)\rightarrow\Omega^{\bullet}(H)\otimes\Omega^{\bullet}(H)$ is a morphism of differential graded algebras. The next result shows that $\Delta_{A}^{\wedge}$ is understood as a graded coaction itself. 
\begin{lemma}
	We have \eqa{\label{right H-colinearity of the extension}(\Delta_{A}^{\wedge}\otimes \id)\circ \Delta_{A}^{\wedge} = (\id \otimes \Delta^{\bullet})\circ \Delta_{A}^{\wedge}.} Moreover $\Delta_{A}^{\wedge}$ is uniquely determined as extension of $\Delta_{A}$.   
	\proof For $a\in \Omega^{0}(A)=A$ everything follows as $\Delta_{A}^{\wedge}$ extends $\Delta_{A}$. Let $a^{0}\dd a^{1}\wedge\dots \wedge \dd a^{k+1}$ be in $\Omega^{k+1}(A)$ and assume by induction that equation \eqref{right H-colinearity of the extension} holds on $\Omega^{k}(A)$. We have \eqa{ (\Delta_{A}^{\wedge}\otimes \id) \Delta_{A}^{\wedge} (a^{0}\dd a^{1}\wedge\dots \wedge\dd a^{k+1}) & = (\Delta_{A}^{\wedge}\otimes \id)\Delta_{A}^{\wedge}(a^{0}\dd a^{1}\wedge\dots \wedge \dd a^{k})(\wedge_{A}\otimes \wedge_{H})\Delta_{A}^{\wedge}(\dd a^{k+1}) \\ & = (\id\otimes\Delta^{\bullet})\Delta_{A}^{\wedge}(a^{0}\dd a^{1}\wedge\dots \wedge \dd a^{k})(\wedge_{A}\otimes \wedge_{H})\Delta_{A}^{\wedge}(\dd a^{k+1})\\ & = (\id\otimes \Delta^{\bullet})\Delta_{A}^{\wedge}(a^{0}\dd a^{1}\wedge\dots \wedge\dd a^{k}\wedge\dd a^{k+1}).\nonumber}  \qed
\end{lemma}
\begin{lemma}
	 $\Omega^{\bullet}(A)$ is right $H-$covariant differential calculus with respect to the coaction $$\Delta_{\Omega^{\bullet}(A)}=(\id\otimes \pi_{0})\Delta_{A}^{\wedge}:  \Omega^{\bullet}(A)\rightarrow \Omega^{\bullet}(A)\otimes H.$$ In other words the following identities hold: \eqa{ & ( \Delta_{\Omega^{\bullet}(A)} \otimes \id) \circ \Delta_{\Omega^{\bullet}(A)}  = (\id\otimes \Delta)\circ \Delta_{\Omega^{\bullet}(A)}; \\ & (\id \otimes \epsilon)\circ \Delta_{\Omega^{\bullet}(A)}  = \id ;\\ &  (\dd\otimes \id)\circ \Delta_{\Omega^{\bullet}(A)}  = \Delta_{\Omega^{\bullet}(A)} \circ \dd. }
\proof Follows from theorem \ref{equivalent definitions of covariance}. \qed  
\end{lemma}

The next proposition is used to construct a non-commutative analogue of the verticalising homeomorphism of classical differential geometry. 
\begin{prop}\label{surjectivity of piver}
	Let $\Omega^\bullet(A)$ be a complete calculus on $A$. Then there is a surjective morphism of differential graded algebras
$$
\pi_{\ver}\colon\Omega^\bullet(A)\to\ver^\bullet(A)
$$
extending the identity $\id\colon A\to A$. It is determined on homogeneous elements by
\begin{equation}
\pi_{\ver}|_{\Omega^k(A)}:=(\mathrm{id}\otimes(\pi_{inv}\circ\pi_{k}))\circ\Delta_A^\wedge|_{\Omega^k(A)}\colon\Omega^k(A)\to\ver^k(A)
\end{equation}
for all $k>0$, where $\pi_{k}\colon\Omega^\bullet(H)\to\Omega^k(H)$ is the obvious projection.  Explicitly, we have
\begin{equation}\label{piverExpl}
\pi_\ver(a^0\mathrm{d}a^1\wedge\cdots\wedge\mathrm{d}a^k)
=a^0_0a^1_0\ldots a^k_0\otimes S(a^0_1a^1_1\ldots a^k_1)a^0_2\mathrm{d}a^1_2\wedge\cdots\wedge\mathrm{d}a^k_2
\end{equation}
for all $a^0,\ldots,a^k\in A$.
\proof Consider $\pi_{inv}:\Omega^{\bullet}(H)\rightarrow \Lambda^{\bullet}$ defined by $\pi_{inv}(\omega)=S(\omega_{-1})\omega_{1}$. We may equivalently write $$ \pi_{inv}(\omega)= \triangleright\circ (S\otimes \id) \circ {}_{\Gamma^{\wedge}}\Delta(\omega).$$  
Let $\omega=a\in H$. We find \eqa{\nonumber\pi_{inv}(\omega)& = \pi_{inv}(a) \\ & = \triangleright \circ (S\otimes \id) \circ {}_{\Gamma^{\wedge}}\Delta(a) \\ & = \triangleright \circ (S\otimes \id) \circ \Delta(a) \\ & = \triangleright \circ (S\otimes \id) \circ (a_{1}\otimes a_{2}) \\ & = S(a_{1})a_{2}\\ & = \epsilon(a) } for any $\omega\in \Omega^{\bullet}(A)$.  
 Let $\omega=\dd b \in \Omega^{1}(H)=\Gamma$ for some $b\in H$.  We have \eqa{\nonumber\pi_{inv}(\omega) & = \pi_{inv}(\dd b) \\ & = \triangleright\circ (S\otimes \id) \circ {}_{\Gamma^{\wedge}}\Delta (\dd b) \\ & = \triangleright \circ (S\otimes \id) \circ  {}_{\Gamma^{\wedge}}\Delta(\dd b) \\ & =  \triangleright \circ (S\otimes \id) (\id\otimes \dd)(b_{1}\otimes b_{2}) \\ & =  \triangleright \circ (b_{1}\otimes \dd b_{2}) \\ & = S(b_{1})  \dd b_{2} \\ & = S(b_{1})\dd b_{2} \\ &  = \varpi(b).   }
Let now $\omega\in \Omega^{\bullet}(H)$.  Then $\omega=a^{0}\dd a^{1}\wedge\cdots \wedge \dd a^{k}$ for some $a^{0},\dots , a^{k}\in A$.  We find \eqa{\nonumber  \pi_{inv}(\omega) & = \pi_{inv}(a^{0}\dd a^{1}\wedge\cdots \wedge\dd a^{k}) \\ & = \triangleright\circ (S\otimes \id) \circ {}_{\Gamma^{\wedge}}\Delta (a^{0}\dd a^{1}\wedge\cdots \wedge \dd a^{k}) \\ & = \triangleright\circ (S\otimes \id) \circ \Delta(a^{0})(\wedge\otimes \wedge) (\id\otimes \dd)\Delta(a^{1})(\wedge\otimes \wedge)\cdots (\wedge\otimes \wedge)(\id\otimes \dd)\Delta(a^{k})\\ & = \triangleright\circ (S\otimes \id)\circ (a^{0}_{1}\otimes a^{0}_{2})(\wedge\otimes \wedge)(a^{1}_{1}\otimes \dd a^{1}_{2})(\wedge\otimes \wedge) \cdots (\wedge\otimes \wedge)(a^{k}_{1}\otimes \dd a^{k}_{2}) \\ & = S(a^{0}_{1}a^{1}_{1}\cdots a^{k}_{1}) a^{0}_{2}\dd a^{1}_{2}\wedge\cdots \wedge \dd a^{k}_{2} \\ & = \epsilon(a^{0})S(a_{1}^{1}a_{1}^{2}\cdots a_{1}^{k}) \dd a_{2}^{1}\wedge\cdots \wedge \dd a_{2}^{k}.} Let us now consider a  morphism $\pi_{\ver}:\Omega^{\bullet}(A)\rightarrow\ver^{\bullet}(A)$ by \eqa{\pi_{\ver}(\omega)= (\mathrm{id}\otimes(\pi_{inv}\circ\pi^k))\circ\Delta_{A}^{\wedge}(\omega). } $\pi_{\ver}$ is obviously well defined. Moreover, given $a\in A$ we find 
$$ \pi_{\ver}(a)= (\id\otimes \pi_{inv}\circ \pi_{0})\Delta_{A}(a) = ( a_{1}\otimes \pi_{inv}\circ \pi_{0}(a_{2})) = a\otimes 1 = a, $$ 
and so $\pi_{\ver}|_{A}=\id_{A}.$  

The morphism $\pi_{\ver}$ is degree preserving. Indeed, given a $k$-form $\omega=a^{0}\dd a^{1}\wedge\cdots \wedge\dd a^{k}\in \Omega^{k}(A)$, we find \eqa{\nonumber (\id & \otimes \pi_{k})\Delta_{A}^{\wedge}(a^{0} \dd a^{1}\wedge\dots \wedge\dd a^{k})  \\ & = (\id \otimes \pi_{k})\Delta_{A}(a^{0})(\id\otimes \dd + \dd \otimes \id)\Delta_{A}(a^{1}) (\wedge\otimes \wedge)\cdots (\wedge\otimes \wedge)(\id \otimes \dd + \dd \otimes \id)\Delta_{A}(a^{k}) \\ & =  a^{0}_{1}\cdots a^{k}_{1} \otimes a_{2}^{0} \dd a_{2}^{1}\wedge\cdots \wedge\dd a_{2}^{k}.  } Moreover \eqa{\nonumber \pi_{\ver}(\omega) & = a^{0}_{1}\cdots a_{1}^{k} \otimes \pi_{inv}(a^{0}_{2}\dd a^{1}_{2}\wedge\cdots \wedge \dd a_{2}^{k}) \\ & =   a^{0}_{1}\cdots a_{1}^{k} \otimes \epsilon (a^{0}_{2})S(a_{2_{1}}^{1}a_{2_{1}}^{2}\cdots a_{2_{1}}^{k}) \dd a_{2_{2}}^{1}\wedge\cdots \wedge \dd a_{2_{2}}^{k} \\ & = a_{1}^{0}\epsilon(a_{2}^{0}) \cdots a_{1}^{k} \otimes S(a_{1}^{1}\cdots a_{1}^{k}) \dd a^{1}_{2}\wedge\cdots \wedge \dd a^{k}_{2} \\ & = a^{0}a_{0}^{1}\cdots a^{k}_{0}\otimes S(a_{1}^{1}\cdots a_{1}^{k}) \dd a^{1}_{2}\wedge\cdots \wedge \dd a^{k}_{2}, } providing an explicit formula for $\pi_{\ver}$. 

The second tensor factor in the above formula defines indeed a coinvariant form in $\Lambda^{\bullet}$. By induction:\eqa{\nonumber{}_{\Gamma^{\wedge}}\Delta[S(a^{1}_{1}\cdots a^{k+1}_{1})(\dd a^{1}_{2}\wedge\cdots &\wedge \dd a^{k+1}_{2})] \\ &  = {}_{\Gamma^{\wedge}}\Delta[S(a^{k+1}_{1})S(a_{1}^{1}\cdots a_{1}^{k})(\dd a_{2}^{1}\wedge\cdots \wedge\dd a^{k+1}_{2})]\\ & = {}_{\Gamma^{\wedge}}\Delta(S(a^{k+1}_{1})) {}_{\Gamma^{\wedge}}\Delta[S(a_{1}^{1}\cdots a_{1}^{k})(\dd a_{2}^{1}\wedge\cdots \wedge\dd a^{k}_{2})]{}_{\Gamma}\Delta(\dd a^{k+1}_{2}) \\ & = \Delta(S(a^{k+1}_{1}) )[1\otimes S(a_{1}^{1}\cdots a_{1}^{k})\dd a^{1}_{2}\wedge\cdots \wedge \dd a^{k}_{2}] (\id\otimes \dd)(a^{k+1}_{0}\otimes a^{k+1}_{1}) \\ & = (S(a_{2}^{k+1})\otimes S(a^{k+1}_{1}))[1\otimes S(a_{1}^{1}\cdots a_{1}^{k})\dd a^{1}_{2}\wedge\cdots \wedge \dd a^{k}_{2}] (\id\otimes \dd)(a^{k+1}_{3}\otimes a^{k+1}_{4})\\ & = 1\otimes S(a_{1}^{1}\cdots a_{1}^{k+1})\dd a_{2}^{1}\wedge\cdots \wedge\dd a_{2}^{k+1}.} 
We now show that $\pi_{\ver}$ is a morphism of differential graded algebras. Thanks to lemma \ref{morphism of DGAs on fixed degree forms} this amounts to show \eqa{\label{piver is a morphism of DGA}\pi_{\ver}(a^{0}\dd a^{1}\wedge\cdots \wedge\dd a^{k}) & = \pi_{\ver}(a^{0})\wedge_{\ver}\dd_{\ver}\circ\pi_{\ver}(a^{1})\wedge_{\ver}\cdots \wedge_{\ver}\dd_{\ver}\circ \pi_{\ver}(a^{k}) \\ & =  (a^{0}\otimes 1)\wedge_{\ver}(a^{1}_{0}\otimes \varpi(a^{1}_{1}))\wedge_{\ver}\cdots \wedge_{\ver}(a^{k}_{0}\otimes \varpi(a^{k}_{1})). }  We proceed by induction. For a $1-$form $\omega=a\dd b\in \Omega^{1}(A)$ we have $\pi_{\ver}(a \dd b) = a b_{0}\otimes \varpi(b_{1})$, whereas \eqa{\nonumber \pi_{\ver}(a)\wedge_{\ver} \dd_{\ver}\circ \pi_{\ver}(b_{1}) & = (a\otimes 1) \wedge_{\ver}(b_{0}\otimes \varpi(b_{1})) \\ & = a b_{0}\otimes (1\leftharpoonup b_{1})\wedge\varpi(b_{2}) \\ & = a b_{0}\otimes S(b_{1})b_{2} \wedge \varpi(b_{3}) \\ & = a b_{0} \otimes \varpi(b_{1}).}  

Assuming \eqref{piver is a morphism of DGA} hold for $k-$forms we have \eqa{ \pi_{\ver}(a^{0})& \wedge_{\ver} \dd_{\ver}\circ \pi_{\ver}(a^{k}) \wedge_{\ver}\cdots \wedge_{\ver} \dd_{\ver}\circ \pi_{\ver}(a^{k+1}) \\ & = \pi_{\ver}(a^{0}\dd a^{1}\wedge\cdots \wedge \dd a^{k}) \wedge_{\ver} \dd_{\ver}\circ \pi_{\ver}( a^{k+1}) \\ & = (a^{0}a^{1}_{0}\cdots a^{k}_{0})\otimes S(a^{1}_{1}\cdots a_{1}^{k}) \dd a_{2}^{1}\wedge\cdots \wedge\dd a^{k}_{2} \wedge_{\ver}(a^{k+1}_{0} \otimes \varpi(a_{1}^{k+1})) \\ & = (a^{0}a^{1}_{0}\cdots a^{k+1}_{0})\otimes (S(a^{1}_{1}\cdots a_{1}^{k}) \dd a_{2}^{1}\wedge\cdots \wedge\dd a^{k}_{2}\leftharpoonup a_{1}^{k+1}) \wedge\varpi(a_{2}^{k+1}) \\ & = a^{0}a^{1}_{0}\cdots a^{k+1}_{0} \otimes S(a_{1}^{k+1})S(a^{1}_{1}\cdots a_{1}^{k}) \dd a_{2}^{1}\wedge\cdots \wedge\dd a^{k}_{2} a_{2}^{k+1} \wedge S(a_{3}^{k+1})\dd a_{4}^{k+1} \\ & = a^{0}a^{1}_{0}\cdots a^{k+1}_{0} \otimes S(a_{1}^{k+1})S(a^{1}_{1}\cdots a_{1}^{k}) \dd a_{2}^{1}\wedge\cdots \wedge\dd a^{k}_{2} \wedge a_{2}^{k+1} S(a_{3}^{k+1})\dd a_{4}^{k+1} \\ & = a^{0}a^{1}_{0}\cdots a^{k+1}_{0} \otimes S(a^{1}_{1}\cdots a^{k+1}_{1}) \dd a^{1}_{2}\wedge\cdots \wedge \dd a^{k+1}_{2}.}  
By proposition \ref{surjective morphism of DGA} we now know that $\pi_{\ver}:\Omega^{\bullet}(A)\rightarrow \ver^{\bullet}(A)$ is  unique as extension of $\id_{A}:A\rightarrow A$. Moreover $\pi_{\ver}$ is also surjective as $\id_{A}$ is. 

\qed 
\end{prop} 
\begin{prop}\label{deltaver}
	The vertical forms $\ver^\bullet(A)$ are complete, i.e. the right $H$-coaction $\Delta_A\colon A\to A\otimes H$ extends to a morphism $\Delta_\ver^\wedge\colon\ver^\bullet(A)\to\ver^\bullet(A)\otimes\Omega^\bullet(H)$ of differential graded algebras. Moreover, the diagram $$\begin{tikzcd}
 \Omega^{\bullet}(A)\arrow[dd," \pi_{\ver}",swap] \arrow[rr," \Delta_{A}^{\wedge} "] & & \Omega^{\bullet}(A) \otimes \Omega^{\bullet}(H)	\arrow[dd,"(\pi_{\ver}\otimes \id) "] \\ \\  \ver^{\bullet}(A) \arrow[rr, " \Delta_{\ver}^{\wedge}"] & & \ver^{\bullet}(A)\otimes \Omega^{\bullet}(H)
 \end{tikzcd} $$ 
commutes, i.e. $\Delta_{\ver}^{\wedge}\circ \pi_{\ver} = (\pi_{\ver}\otimes \id)\Delta_{A}^{\wedge}.$ 
\proof 
We want to define $\Delta_{\ver}^{\wedge}$ by commutativity of the above diagram. To that purpose, let us consider a morphism \eqa{ \widetilde{\Delta}_{\ver}^{\wedge}: \Omega^{\bullet}(A)\rightarrow \ver^{\bullet}\otimes \Omega^{\bullet}(H), \quad \omega \mapsto (\pi_{\ver}\otimes \id)\circ \Delta_{A}^{\wedge}(\omega).}  Let $\omega=a^{0}\dd a^{1}\wedge\cdots \wedge \dd a^{k}\in \Omega^{k}(A)$. If in particular $\omega\in \ker\pi_{\ver}$ we have \eqa{\label{piver zero}\pi_{\ver}(\omega)= a^{0}_{0}a^{1}_{0}\dots a^{k}_{0}\otimes S(a_{1}^{0}\dots a^{k}_{1})a^{0}_{2} \dd a^{1}_{2}\wedge\cdots \wedge \dd a^{k}_{2}=0 .} We now show that even $\widetilde{\Delta}_{\ver}^{\wedge}(\omega)=0$.  Explicitly it reads  

\eqa{\label{summands under deltaver} \widetilde{\Delta}_{\ver}^{\wedge}(\omega) & = (\pi_{\ver}\otimes \id)\circ \Delta_{A}^{\wedge}(a^{0}\dd a^{1}\wedge\cdots \wedge \dd a^{k}) \\ & = \pi_{\ver}(a^{0}_{0}\dd a^{1}_{0}\wedge\cdots\wedge \dd a^{k}_{0})\otimes a^{0}_{1}\cdots a^{k}_{1} \\ & + \pi_{\ver}(a^{0}_{0}a^{1}_{0}\dd a^{2}_{0}\wedge\cdots \wedge \dd a^{k}_{0})\otimes a_{1}^{0} \dd a^{1}_{1}\cdots a^{k}_{0} \\ & + \cdots \\  & + \pi_{\ver}(a^{0}_{0}\cdots a^{k-1}_{0} \dd a^{k}_{0}) \otimes a_{1}^{0}\dd a^{1}_{1}\wedge\cdots \wedge \dd a^{1}_{k-1} a^{k}_{0} \\ & + \pi_{\ver}(a^{0}_{0}\cdots a^{k}_{0}) \otimes a^{0}_{1}\dd a^{1}_{1}\wedge\cdots \wedge \dd a^{k}_{1}, } where the dots stand for all the other possible combinations in $\Omega^{r}(A)\otimes \Omega^{s}(H)$ such that $s+r=k$. 

Considering the first tensor factor of the summand in the above expression we find \eqa{ \nonumber (\Delta_{A}&\otimes \Delta_{\Omega^{k}(A)})\circ  \pi_{\ver}(a^{0}_{0}\dd a^{1}_{0}\wedge\cdots\wedge \dd a^{k}_{0})\\ & = (\Delta_{A}\otimes \Delta_{\Omega^{k}(A)}) \circ (a^{0}_{0}\cdots a^{k}_{0}\otimes S(a^{0}_{1}\cdots a^{k}_{1}) a^{0}_{2}\dd a^{1}_{2}\wedge\cdots \wedge \dd a^{k}_{2} )     \\ & = a^{0}_{0}\cdots a^{k}_{0}\otimes a^{0}_{1}\cdots a^{k}_{1} \otimes S(a^{0}_{1}\cdots a^{k}_{1})_{0} a^{0}_{20}\dd a^{1}_{20}\wedge\cdots \dd a^{k}_{20} \otimes S(a^{0}_{1}\cdots a^{k}_{1})_{1}a^{0}_{21}\cdots a^{k}_{21}  \\ & = a^{0}_{0}\cdots a^{k}_{0}\otimes a^{0}_{1}\cdots a^{k}_{1}\otimes S(a^{0}_{11}\cdots a^{k}_{11})a^{0}_{20}\dd a^{1}_{20}\wedge\cdots \wedge a^{k}_{20}\otimes S(a^{0}_{10}\cdots a^{k}_{10})a^{0}_{21}\cdots a^{k}_{21}\\      & = a^{0}_{0}\cdots a^{k}_{0}\otimes a^{0}_{1}\cdots a^{k}_{1} \otimes S(a^{0}_{3}\cdots a^{k}_{3})a^{0}_{4}\dd a^{1}_{4}\wedge\cdots \wedge \dd a^{k}_{4} \otimes S(a^{0}_{2}\cdots a^{k}_{2})a^{0}_{5}\cdots a^{k}_{5} .    } Multiplying the second and last tensor factors in the result of the last expression we have \eqa{a^{0}_{0}\cdots a^{k}_{0}\otimes S(a^{0}_{1}\cdots a^{k}_{1}) a^{0}_{2}\dd a^{1}_{2}\wedge\cdots \wedge \dd a^{k}_{2} \otimes a^{0}_{3}\cdots a^{k}_{3}=0,} according to Equation \eqref{piver zero}. Therefore \eqa{a^{0}_{0}\cdots a^{k}_{0}\otimes S(a^{0}_{1}\cdots a^{k}_{1}) a^{0}_{2}\dd a^{1}_{2}\wedge\cdots \wedge \dd a^{k}_{2} \otimes a^{0}_{3}\cdots a^{k}_{3}\otimes a^{0}_{4}\cdots a^{k}_{4}=0. }  
The above reasoning can be repeated for every summand in Equation \eqref{summands under deltaver}. Therefore we conclude  $\widetilde{\Delta}_{\ver}^{\wedge}(\omega)=0.$ Thus, $\widetilde{\Delta}_{\ver}^{\wedge}:\Omega^{\bullet}(A)\rightarrow \ver^{\bullet}(A)\otimes\Omega^{\bullet}(H)$ descends to the morphism \eqa{\Delta_{\ver}^{\wedge}:\ver^{\bullet}(A)\cong \Omega^{\bullet}(A)/\ker\pi_{\ver}\rightarrow \Omega^{\bullet}(A)\otimes \Omega^{\bullet}(H), }  which makes the diagram commute. Since every other map in the diagram is a morphism of differential graded algebras also $\Delta_{\ver}^{\wedge}$ is.  \qed  
\end{prop}
\begin{exmp}
	We consider a 2-form $\omega=a\dd b\wedge \dd c\in \Omega^{2}(A)$ in the kernel of $\pi_{\ver}$, i.e. \eqa{\label{piverex} a_{0}b_{0}c_{0}\otimes S(a_{1}b_{1}c_{1})a_{2}\dd b_{2}\wedge\dd c_{2}=0,} and show that $\widetilde{\Delta}_{\ver}^{\wedge}(\omega)=0$. We have \eqa{ \nonumber \Delta_{A}^{\wedge}(a \dd b \wedge \dd c) & = a_{0}b_{0}c_{0}\otimes a_{1}\dd b_{1}\wedge \dd c_{1} \\ & + a_{0} \dd b_{0} \wedge \dd c_{0} \otimes a_{1}b_{1}c_{1} \\ & + a_{0} \dd b_{0}c_{0}\otimes a_{1}b_{1} \dd c_{1} \\ & - a_{0}b_{0}\dd c_{0}\otimes a_{1}\dd b_{1} c_{1}.}  Accordingly \eqa{\label{piverexterms}\widetilde{\Delta}_{\ver}^{\wedge}(\omega) & = (\pi_{\ver}\otimes \id)\circ \Delta_{A}^{\wedge}(a \dd b \wedge \dd c) \\ (i) \ & = a_{0}b_{0}c_{0}\otimes a_{1}b_{1}c_{1}\otimes a_{2}\dd b_{2}\wedge \dd c_{2} \\ (ii) \ & + a_{0}b_{0}c_{0}\otimes S(a_{1}b_{1}c_{1})a_{2}\dd b_{2} \wedge \dd c_{2} \otimes a_{3}b_{3}c_{3} \\(iii) \  & + a_{0}b_{0}c_{0}\otimes S(a_{1}b_{1}c_{1})a_{2}\dd b_{2} c_{2}\otimes a_{3}b_{3} \dd c_{3} \\(iv) \ & - a_{0}b_{0}c_{0}\otimes S(a_{1}b_{1}c_{1})a_{2}b_{2}\dd c_{2} \otimes a_{3}\dd b_{3}c_{3}. }	Considering Equation \eqref{piverex} we coact on the first tensor factor to obtain \eqa{\nonumber (\Delta_{A}\otimes \id)(a_{0}b_{0}c_{0}&\otimes S(a_{1}b_{1}c_{1})a_{2}\dd b_{2}\wedge\dd c_{2})=0 \\ & = a_{0}b_{0}c_{0}\otimes a_{1}b_{1}c_{1}\otimes S(a_{2}b_{2}c_{2})a_{3}\dd b_{3}\wedge\dd c_{3}.}  Multiplying the first and last tensor factor we have \eqa{a_{1}b_{1}c_{1}\otimes a_{2}\dd b_{2}\wedge \dd c_{2}=0.} Therefore term $(i)$ of equation \eqref{piverexterms} vanishes. 	Term $(ii)$ is zero by the derivation in Proposition \eqref{deltaver}. Finally, considering again Equation \eqref{piverex} we have \eqa{ \nonumber (\Delta_{A}\otimes \pi_{1}\otimes \pi_{1})&\circ (\id\otimes \Delta^{\bullet})\Delta_{A}(a_{0}b_{0}c_{0})\otimes S(a_{1}b_{1}c_{1})a_{2}\dd b_{2}\wedge\dd c_{3} )= 0\\ & = \Delta_{A}(a_{0}b_{0}c_{0}) \otimes (\pi_{1}\otimes \pi_{1})\Delta^{\bullet}(S(a_{1}b_{1}c_{1})a_{2}\dd b_{2}\wedge \dd c_{2}) \\ & = \Delta_{A}(a_{0}b_{0}c_{0})\otimes S(a_{1}b_{1}c_{1})_{1}a_{2}\dd b_{2}c_{2}\otimes S(a_{1}b_{1}c_{1})_{2}a_{3}b_{3}\dd c_{3} \\ & - \Delta_{A}(a_{0}b_{0}c_{0})\otimes S(a_{1}b_{1}c_{1})_{1}a_{2}b_{2}\dd c_{2}\otimes S(a_{1}b_{1}c_{1})_{2}a_{3}\dd b_{3}c_{3} \\ & = a_{0}b_{0}c_{0}\otimes a_{1}b_{1}c_{1}\otimes S(a_{22}b_{22}c_{22})a_{3}\dd b_{3}c_{3}\otimes S(a_{21}b_{21}c_{21})a_{4}b_{4}\dd c_{4} \\ & - a_{0}b_{0}c_{0}\otimes a_{1}b_{1}c_{1}\otimes S(a_{22}b_{22}c_{22})a_{3}b_{3}\dd c_{3}\otimes S(a_{21}b_{21}c_{21})a_{4}\dd b_{4}c_{4} \\ & = a_{0}b_{0}c_{0}\otimes a_{1}b_{1}c_{1}\otimes S(a_{3}b_{3}c_{3})a_{4}\dd b_{4}c_{4}\otimes S(a_{2}b_{2}c_{2})a_{5}b_{5}\dd c_{5} \\ & - a_{0}b_{0}c_{0}\otimes a_{1}b_{1}c_{1}\otimes S(a_{3}b_{3}c_{3})a_{4}b_{4}\dd c_{4}\otimes S(a_{2}b_{2}c_{2})a_{5}\dd b_{5}c_{5}. } Multiplying the second and last tensor factors in the last expression's result we have \eqa{\nonumber a_{0}b_{0}c_{0}\otimes S(a_{1}b_{1}c_{1})a_{2}\dd b_{2}c_{2}\otimes a_{3}b_{3}\dd c_{3} - a_{0}b_{0}c_{0}\otimes S(a_{1}b_{1}c_{1})a_{2}b_{2}\dd c_{2}\otimes a_{3}\dd b_{3}c_{3}=0,  } and this holds if and only both summands are zero. Therefore terms $(iii)$ and $(iv)$ of Equations \eqref{piverexterms} vanish. Accordingly $\widetilde{\Delta}_{\ver}^{\wedge}(\omega)=0$. \qed  
\end{exmp}
\begin{rmk}
	It is also possible to work out the explicit form for $\pi_{\ver}$ by assuming it to be the unique morphism of differential graded algebras extending the identity $\id:A\rightarrow A$.  Let $a\dd b\in \Omega^{1}(A)$.  We find \eqa{ \nonumber\pi_{\ver}(a\dd b) & = a\pi_{\ver}(\dd b) \\ & = a \dd_{\ver}(b\otimes 1) \\ & = a(b_{0}\otimes \varpi(b_{1})) \\ & = ab_{0}\otimes S(b_{11}) \dd (b_{12}) \\ & = a b_{0}\otimes S(b_{1}) \dd b_{2}.}  Let now $a\dd b\wedge\dd c\in \Omega^{2}(A)$. We have \eqa{\nonumber \pi_{\ver}(a\dd b \wedge\dd c) & =(ab_{0}\otimes S(b_{1})\dd b_{2})\wedge_{\ver}(c_{0}\otimes \varpi(c_{1})) \\ & = ab_{0}c_{0}\otimes (\varpi(b_{1})\leftharpoonup c_{1})\wedge \varpi(c_{2}) \\ & = a b_{0}c_{0} \otimes S(c_{10})\varpi(b_{1})c_{11}\wedge\varpi(c_{2}) \\ & = ab_{0}c_{0}\otimes S(c_{1})S(b_{1})\dd b_{2}c_{2}\wedge S(c_{3})\dd c_{4} \\ & = ab_{0}c_{0}\otimes S(b_{1}c_{1}) \dd b_{2}\wedge c_{2}S(c_{3})\dd c_{4} \\ & = ab_{0}c_{0}\otimes S(b_{1}c_{1}) (\dd b_{2}\wedge \dd c_{2}).  } Generalising for $\omega=a^{0}\dd a^{1}\wedge\cdots \wedge \dd a^{k}\in \Omega^{k}(A)$ it is easy to see\eqa{ \label{alternative form for piver} \pi_{\ver}(a^{0}\dd a^{1}\wedge\cdots \wedge\dd a^{k}) = a^{0}a^{1}_{0}\cdots a^{k}_{0}\otimes S(a^{1}_{1}\cdots a^{k}_{1}) (\dd a^{1}_{2}\wedge\cdots \wedge\dd a^{k}_{2}).}  \qed 
\end{rmk}
\medskip
We now introduce the non-commutative analogue of horizontal forms of classical geometry.
\begin{defn}
	Let $\Omega^{\bullet}(A)$ be a complete differential calculus on the quantum principal bundle $B:=A^{coH}\subseteq A$. The horizontal forms of the bundle are defined as elements of the $\Omega^{\bullet}(A)-$subalgebra \eqa{ \hor^{\bullet}(A):=(\Delta_{A}^{\wedge})^{-1}(\Omega^{\bullet}(A)\otimes H).} 
\end{defn}
Notice that since since $\Omega^{0}(A)=A$ we  have $(\Delta_{A}^{\wedge})^{-1}(A\otimes H)=(\Delta_{A})^{-1}(A\otimes H) =A, $ i.e. $\hor^{0}(A)=A$. This also means $\Omega^{0}(A)=\ver^{0}(A)=\hor^{0}(A).$
\begin{prop}
	$\hor^{\bullet}(A)$ is a graded algebra. In other words, given $\omega\in\hor^{k}(A)$ and $\gamma\in \hor^{\ell}(A)$ we have $\omega\wedge\gamma \in \hor^{k+\ell}(A)$.  \proof Given $\omega$ and $\gamma$ as above we find \eqa{\Delta_{A}^{\wedge}(\omega\wedge\gamma) & = \Delta_{A}^{\wedge}(\omega)(\wedge\otimes \id) \Delta_{A}^{\wedge}(\gamma) \in (\Omega^{k}(A)\wedge\Omega^{\ell}(A)) \otimes H = \Omega^{k+\ell}(A)\otimes H,} therefore the wedge product of horizontal forms gives an horizontal forms.  \qed   
\end{prop}
\begin{lemma}
	The graded algebra $\hor^{\bullet}(A)$ is $\Delta_{A}^{\wedge}-$invariant, namely \eqa{\Delta_{A}^{\wedge}(\hor^{\bullet}(A))\subseteq \hor^{\bullet}(A)\otimes H.  } In other words, $\hor^{\bullet}(A)$ is a right $H-$comodule algebra. 
	\proof 
	
	Let $\varphi\in \hor^{\bullet}(A)$. By the definition of horizontal forms we have \eqa{(\Delta_{A}^{\wedge}\otimes \id)\circ \Delta_{A}^{\wedge}(\varphi) & = (\id\otimes \Delta^{\bullet}) \circ \Delta_{A}^{\wedge}(\varphi) \\ & = (\id \otimes \Delta)\circ \Delta_{\Omega^{\bullet}_{A}}(\varphi),} which is an element in $\Omega^{\bullet}(A)\otimes H\otimes H$; accordingly also $(\Delta_{A}^{\wedge}\otimes \id)\circ \Delta_{A}^{\wedge}(\varphi)$ is an element of the same space. It follows that $\Delta_{A}^{\wedge}(\varphi) \in \hor^{\bullet}(A)\otimes H$, since then $$\Delta_{A}^{\wedge}(\hor^{\bullet}(A))\otimes H =\Delta_{A}^{\wedge}((\Delta_{A}^{\wedge})^{-1}(\Omega^{\bullet}(A)\otimes H))\otimes H = \Omega^{\bullet}(A)\otimes H \otimes H. $$  Accordingly $\hor^{\bullet}(A)$ is a right $H-$comodule algebra by the coaction $$\Delta_{\hor}:=\Delta_{A}^{\wedge}|_{\hor^{\bullet}(A)}:\hor^{\bullet}(A)\rightarrow\hor^{\bullet}(A)\otimes H.$$  \qed 
\end{lemma}

\begin{rmk}
	 Horizontal forms are not closed under the differential. This is true already for horizontal $1-$forms. Indeed  \eqa{ \nonumber \Delta_{A}^{\wedge}(\dd\omega) & = (\id\otimes \dd + \dd \otimes \id)\Delta_{A}^{\wedge}(\omega) \\ & = (\id\otimes \dd + \dd \otimes \id)\Delta_{A}^{\wedge}(a\dd b) \\ & = (\id\otimes \dd + \dd \otimes \id)(a_{0}\dd b_{0}\otimes a_{1}b_{1}) \\ & =- a_{0}\dd b_{0}\otimes \dd(a_{1}b_{1}) +   \dd a_{0}\wedge \dd b_{0}\otimes a_{1}b_{1},  }  so $\Delta_{A}^{\wedge}(\dd\omega)\in (\Omega^{1}(A)\otimes \Omega^{1}(H))\oplus (\Omega^{2}(A)\otimes H$).  \qed 
\end{rmk}
The next lemma provides an exact sequence of $A-$modules involving vertical, horizontal and total space forms. The result we provide holds only for $1-$forms and is false for higher orders. This, as we will see, happens precisely because forms in the kernel of $\pi_{\ver}$ are not necessarily horizontal for $k>1.$ 
\begin{lemma} There is a short-exact sequence of  $A-$modules given by $$\begin{tikzcd}
 0 \arrow[r] & \hor^{1}(A) \arrow[r,"\iota",hook] & \Omega^{1}(A)\arrow[r,"\pi_{\ver}",twoheadrightarrow] & \ver^{1}(A) \arrow[r] & 0.	
\end{tikzcd}$$  
\proof Injectivity of $\iota$ is obvious and surjectivity of $\pi_{\ver}$ follows from proposition \ref{surjectivity of piver}. Therefore we only need to prove that the kernel of $\pi_{\ver}$ equals the image of $\iota$, i.e. that $\ker(\pi_{\ver})=\hor^{\bullet}(A)$.  

Let $\varphi\in \hor^{k}(A)$. We find \eqa{\nonumber\pi_{\ver}(\varphi) & = (\id\otimes \pi_{inv})(\id\otimes \pi_{k})\Delta_{A}^{\wedge}(\varphi) \\ & =(\id\otimes \pi_{inv})(\id\otimes \pi_{k})(\theta\otimes h) \\ & =  \theta\otimes 0 = 0, }  so $\hor^{\bullet}(A)\subseteq \ker(\pi_{\ver})$, and in particular $\hor^{1}(A)\subseteq \ker(\pi_{\ver}).$    

Let $\omega\in \ker(\pi_{\ver})\cap \Omega^{1}(A)$. We have \eqa{\nonumber\pi_{\ver}(\omega) & = (\id\otimes \pi_{inv})(\id\otimes \pi_{1})\Delta_{A}^{\wedge}(\omega)= 0, } which explicitely reads \eqa{\nonumber ab_{0}\otimes S(b_{1})\dd b_{2}=0.} Applying $(\id \otimes \triangleright)\circ (\id\otimes\Delta_{A})$ to $\pi_{\ver}(\omega)$ we find \eqa{\nonumber(\id \otimes \triangleright)\circ (\id\otimes\Delta_{A}) \circ (ab_{0}\otimes S(b_{1})\dd b_{2}) & = (\id\otimes \triangleright)\circ( (a_{0}\otimes a_{1})(b_{0}\otimes b_{1})\otimes S(b_{2})\dd b_{3})\\ &  =a_{0}b_{0}\otimes a_{1}\dd b_{1}= 0. } But then \eqa{\nonumber\Delta_{A}^{\wedge}(a\dd b)= a_{0}b_{0}\otimes a_{1}\dd b_{1} + a_{0}\dd b_{0}\otimes a_{1}b_{1}=a_{0}\dd b_{0}\otimes a_{1}b_{1},} which means $\Delta_{A}^{\wedge}(\omega)\in \Omega^{1}(A)\otimes H,$ i.e. $\omega\in \hor^{1}(A)$.       \qed

\end{lemma} 
\begin{rmk} As we discussed this exact sequence fails on $\Omega^{k}(A)$ for $k>1$. Indeed, given $\omega\in\ker(\pi_{\ver})\cap \Omega^{k}(A)$ we have $$ \Delta_{A}^{\wedge}(\omega)\in (A\otimes \Omega^{k}(H))\oplus (\Omega^{1}(A)\otimes \Omega^{k-1}(H))\oplus\cdots \oplus (\Omega^{k}(A)\otimes H), $$ and the fact that $\omega$ is also in the kernel of $\pi_{\ver}$ is sufficient to obtain \eqa{\nonumber\pi_{\ver}(a^{0}\dd a^{1}\wedge\cdots \wedge \dd a^{k}) & = a^{0} a^{1}_{0}\cdots a^{k }_{0} \otimes S(a^{1}_{1}\cdots a^{k}_{1})\dd a^{1}_{2}\wedge\cdots \wedge \dd a^{k}_{2}=0, } from which \eqa{ \nonumber (\Delta_{A}\otimes \id) \circ \pi_{\ver}(\omega) = a^{0}_{0}\cdots a^{k}_{0}\otimes a^{0}_{1}\cdots a^{k}_{1}\otimes  S(a_{2}^{0}\cdots a_{1}^{k}a^{0}_{2}) a^{0}_{3}\dd a^{1}_{3}\wedge\cdots \wedge \dd a^{k}_{3}=0.  } Multiplying the second and last tensor factor we obtain \eqa{\nonumber a^{0}_{0}\cdots a^{k}_{0}\otimes  a^{0}_{1}\dd a^{1}_{1}\wedge\cdots \wedge \dd a^{k}_{1}=0,}which is not enough to conclude that $\omega$ is also horizontal. \qed      \end{rmk}  
 \subsection{Base space calculus}\label{Base space calculus}
 \begin{defn}
 	Given a complete calculus $\Omega^{\bullet}(A)$ on a quantum principal bundle $B:=A^{coH}\subseteq A$ we define the base space forms as \eqa{\Omega^{\bullet}(B):=\graffe{\omega\in \Omega^{\bullet}(A) \ : \ \Delta^{\wedge}_{A}(\omega)=\omega\otimes 1}.} The differential on $\Omega^{\bullet}(B)$ is the restriction of the differential on $\Omega^{\bullet}(A)$.  
 \end{defn}
 Recall that a first order differential calculus is called \textit{connected} if the kernel of the differential is solely made by scalars, i.e. $\ker \dd = \B{k}$. The next proposition gives a characterisation of the base space forms under connectedness assumption on the total space calculus.  
 \begin{prop}
 	Given a connected, bicovariant, first order differential calculus $(\Gamma,\dd)$ on $H$ we have the identification $\Omega^{\bullet}(B)=\hor^{\bullet}(A)\cap \dd^{-1}(\hor^{\bullet}(A))$, in other words the base space calculus consists of horizontal forms which are mapped to horizontal forms by the differential.  
 	\proof 
 	Let $\omega\in \Omega^{k}(B)$. Then $\Delta_{A}^{\wedge}(\omega)=\omega \otimes 1$, i.e. $\Delta^{\wedge}_{A}(\omega)\in \Omega^{k}(A)\otimes H$, that implies $\omega\in \hor^{k}(A)$.  Moreover \eqa{\nonumber\Delta_{A}^{\wedge}(\omega) & = (\dd \otimes \id + (-1)^{k}\id \otimes \dd) \Delta_{A}^{\wedge}(\omega) \\ & = (\dd \otimes \id + (-1)^{k}\id \otimes \dd)(\omega\otimes 1) \\ & = \dd\omega \otimes 1,} and so $\dd\omega \in \hor^{k+1}(A)$ and accordingly $\omega\in \hor^{k}(A)\cap \dd^{-1}(\hor^{k+1}(A))$.  
 	
 	Let now $\omega\in \hor^{k}(A)\cap \dd^{-1}(\hor^{k+1}(A))$.  Then \eqa{ \Delta_{A}^{\wedge}(a^{0}\dd a^{1}\wedge\cdots \wedge \dd a^{k}) =  a^{0}_{0}\dd a^{1}_{0}\wedge\cdots \wedge \dd a^{k}_{0}\otimes a^{0}_{1}\cdots a^{k}_{1},} and moreover \eqa{ \nonumber \Delta_{A}^{\wedge}(\dd\omega) & = (\dd\otimes \id + (-1)^{k}\id \otimes \dd) \Delta_{A}^{\wedge}(\omega) \\ & = (\dd\otimes \id + (-1)^{k}\id \otimes \dd) (a^{0}_{0}\dd a^{1}_{0}\wedge\cdots \wedge \dd a^{k}_{0}\otimes a^{0}_{1}\cdots a^{k}_{1}) \\ & = \dd a^{0}_{0}\wedge\cdots \wedge \dd a^{k}_{0}\otimes a_{1}^{0}\cdots a^{k}_{1} \\ & +(-1)^{k} a^{0}_{0}\dd a^{1}_{0}\wedge\cdots \wedge \dd a^{k}_{0}\otimes \dd(a^{0}_{1}\cdots a^{k}_{1}) \\ & = \dd a^{0}_{0}\wedge\cdots \wedge \dd a^{k}_{0}\otimes a_{1}^{0}\cdots a^{k}_{1}, } and therefore $a^{0}_{0}\dd a^{1}_{0}\wedge\cdots \wedge \dd a^{k}_{0}\otimes \dd(a^{0}_{1}\cdots a^{k}_{1})=0$.  By connectedness assumption on the first order differential calculus $(\Gamma,\dd)$ on $H$ we have $a_{1}^{0}\cdots a^{k}_{1}\in \B{k}$.  Accordingly $$ \Delta_{A}^{\wedge}(\omega)\in \Omega^{k}(A)\otimes \B{k} \quad \Rightarrow \quad \Delta^{\wedge}_{A}(\omega)= \omega \otimes 1 \quad \Rightarrow \quad \omega\in \Omega^{k}(B).$$  \qed 
 \end{prop}
\begin{prop} \label{inclusion of BdB into coinvariant forms} Elements of the form $B\dd B$ are contained in $\Omega^{1}(B)$. 

\proof Let $b,b'\in B$ and consider $b \dd b'$. As $B:=\graffe{a\in A \ : \ \Delta_{A}(a)=a\otimes 1},$ we find \eqa{\nonumber\Delta_{A}^{\wedge}(b\dd b') & = \Delta_{A}(b)(\wedge\otimes \wedge) \Delta_{A}(\dd b') \\ & = (b\otimes 1) (\wedge\otimes \wedge) (\dd \otimes \id + \id \otimes \dd) \Delta_{A}(b') \\ & = (b\otimes 1) (\wedge\otimes \wedge) (\dd b' \otimes 1 ) \\ & = (b\dd b' \otimes 1),} i.e. $b\dd b'\in \Omega^{1}(B)$.  This reasoning can be carried for higher order forms as follows. Consider $b^{0}\dd b^{1}\wedge\cdots \wedge \dd b^{k},$ where $b^{0},\dots, b^{k}\in B$.  Then \eqa{ \nonumber\Delta_{A}^{\wedge}(b^{0}\dd b^{1}\wedge\cdots \wedge\dd b^{k}) & = \Delta_{A}(b^{0})(\wedge\otimes \wedge) (\dd\otimes \id + \id\otimes \dd) \Delta_{A}(b^{1})(\wedge\otimes \wedge)\cdots \\ & \cdots (\wedge\otimes \wedge)(\dd\otimes \id + \id\otimes \dd) \Delta_{A}(b^{k})\\ & = b^{0}\dd b^{1}\wedge\cdots \dd b^{k}\otimes 1, } and so $ B\dd B\wedge\cdots \wedge \dd B\subseteq\Omega^k(B).$ \qed\end{prop} 

\begin{rmk}The other inclusion, $\Omega^{1}(B)\subseteq B\dd B$ does not hold in general. Indeed, let $\omega\in \Omega^{1}(A)$ be such that $\Delta_{A}^{\wedge}(\omega)=\Delta_{A}^{\wedge}(a\dd a')=a\dd a'\otimes 1.$ We have \eqa{\Delta_{A}^{\wedge}(a\dd a') & = \Delta_{A}(a)(\wedge\otimes \wedge) (\dd \otimes \id + \id \otimes \dd) \Delta_{A}(a') \\ & = a_{0} \dd a'_{0}\otimes a_{1}a'_{1} + a_{0}a'_{0}\otimes a_{1}\dd a'_{1}\\ & = a\dd a'\otimes 1  } if and only if $a_{0} \dd a'_{0}\otimes a_{1}a'_{1} = a\dd a'\otimes 1$ and $ a_{0}a'_{0}\otimes a_{1}\dd a'_{1}=0.$  \end{rmk}

\begin{prop}
	Differential forms over the base space $B$ are exactly the horizontal, right $H-$coinvariant forms; in other words $\Omega^{\bullet}(B)=\hor^{\bullet}(A)\cap \Omega^{\bullet}(A)^{coH}$.
	\proof Let $\omega\in \Omega^{k}(B)$ be a differential $k-$form on the base space. By definition $\omega=a^{0}\dd a^{1}\wedge\cdots \wedge \dd a^{k}$ such that $\Delta_{A}^{\wedge}(\omega)=\omega\otimes 1$. Explicitly we have \eqa{ \label{horizontal coinvariant are the base} \Delta_{A}^{\wedge}(a^{0}\dd a^{1}\wedge\cdots \wedge \dd a^{k})& = a^{0}_{0}\dd a^{1}_{0}\wedge\cdots \wedge \dd a^{k}_{0} \otimes a^{0}_{1}\cdots a^{k}_{1} \\ & + (-1)^{k-1}a^{0}_{0} \dd a^{1}_{0}\wedge\cdots \wedge \dd a^{k-1}_{0} a^{k}_{0}\otimes a^{0}_{1}\cdots a^{k-1}\dd a^{k}_{1} \\&  + \dots \\ & - a^{0}_{0}\dd a^{1}_{0}a^{2}_{0}\cdots a^{k}_{0}\otimes a^{0}_{1}a^{1}_{1}\dd a^{2}_{1}\wedge \cdots \wedge \dd a^{k}_{1} \\ & +a^{0}_{0}\cdots a^{k}_{0}\otimes a^{0}_{1}\dd a^{1}_{1}\wedge\cdots \wedge \dd a^{k}_{1},} and since $$\Delta_{A}^{\wedge}(a^{0}\dd a^{1}\wedge\cdots \wedge \dd a^{k})= a^{0}\dd a^{1}\wedge\cdots \dd a^{k}\otimes 1,$$ every summand in the above equation, except the first, vanish. This implies $\Delta_{A}^{\wedge}(a^{0}\dd a^{1}\wedge\cdots \wedge \dd a^{k}) \in \Omega^{k}\otimes H$, and so $\omega$ is horizontal and right $H-$coinvariant. By the same calculation we find that every form $\omega$ that is both horizontal and right $H-$coinvariant, is in $\Omega^{\bullet}(B)$.  \qed  
\end{prop}
\subsection{Comparison with the Brzezi\'nski-Majid approach}\label{equivalence with bm}
The following section supplies a comparison between the Brzezi\'nski-Majid and \DJ ur\dj evi\'c approaches to quantum principal bundles. 
\begin{defn}[\cite{BrzMjd,beggs-majid}]
	Let $(\Gamma_{A},\dd_{A})$ be a right $H-$covariant first order differential calculus over a right $H-$comodule algebra $A$, and let $(\Gamma_{H},\dd_{H})$ be a bicoviariant first order differential calculus over a Hopf algebra $H$. We have a Brzezi\'nski-Majid quantum principal bundle if the map \eqa{\ver_{BM}:\Gamma_{A}\rightarrow A\otimes \Lambda^{1},\quad \ver_{BM}(a\dd_{A}a')=aa'_{0}\otimes \varpi(a'_{1})=aa'_{0}\otimes S(a'_{1})\dd_{H}(a'_{2})} is well defined, and the sequence \eqa{0\rightarrow A\dd_{A}(B)A \hookrightarrow \Gamma_{A}\stackrel{\ver_{BM}}{\longrightarrow} A\otimes \Lambda^{1}\rightarrow 0} is exact. We call $A\dd_{A}(B)A$ the Brzezi\'nski-Majid horizontal forms.
\end{defn}We talk of a quantum principal bundles referring to \DJ ur\dj evi\'c, and of a Brzezi\'nski-Majid quantum principal bundle referring to the last definition. 
	\begin{defn}
		A differential calculus on a quantum principal bundle is called first order complete if the right $H-$coaction $\Delta_{A}:A\rightarrow A\otimes H$ is $1-$differentiable, i.e. if there is a morphism of first order differential calculi \eqa{\Delta_{A}^{1}:\Omega^{1}(A\otimes H)\rightarrow (\Omega^{1}(A)\otimes H)\oplus (A\otimes \Omega^{1}(H))} extending $\Delta_{A}$: $$\begin{tikzcd}
 \Omega^{1}(A)\arrow[rr,"\Delta_{A}^{1}"] & & \Omega^{1}(A\otimes H) \\ \\ A \arrow[uu,"\dd "]	\arrow[rr,"\Delta_{A}"] & & A\otimes H \arrow[uu,"\dd_{\otimes}"]
 \end{tikzcd}$$ 

	\end{defn}
We want to clarify the relation between the definition of quantum principal bundle in the Brzezi\'nski-Majid and \DJ ur\dj evi\'c approaches.
As the Brzezi\'nski-Majid approach only assumes a right $H-$comodule algebra, whereas the \DJ ur\dj evi\'c's requires a faithfully flat Hopf-Galois extension, for a fair comparison we need to consider faithfully flat Hopf-Galois extensions in both scenarios. 
\begin{prop}
	Let $B=A^{coH}\subseteq A$ be a faithfully flat Hopf-Galois extension. \begin{enumerate}
		\item If there is a Brzezi\'nski-Majid quantum principal bundle on $A$, then the maximal prolongation of $\Gamma_{A}$ is first order complete. 
		\item  If there is a first order complete differential calculus $\Omega^{\bullet}(A)$ on $A$ with corresponding calculus  $\Omega^{\bullet}(H)$ on the structure Hopf algebra, then the first order truncation $(\Omega^{1}(A),\Omega^{1}(H))$ is a Brzezi\'nski-Majid quantum pricipal bundle if and only if the horizontal $1-$forms of Brzezi\'nski-Majid and \DJ ur\dj evi\'c' coincide. 
	\end{enumerate}
	\proof 
	Given a Brzezi\'nski-Majid  quantum principal bundle $(\Gamma_{A},\Gamma_{H})$ we show that the maps \eqa{\Delta_{\Gamma_{A}}:\Gamma_{A}\rightarrow \Gamma_{A}\otimes H, \quad a\dd_{A}a'\mapsto a_{0}\dd_{A}a'_{0}\otimes a_{1}a'_{1},} \eqa{\ver:\Gamma_{A}\rightarrow A\otimes \Gamma_{H}, \quad a\dd_{A}a'\mapsto a_{0}a'_{0}\otimes a_{1}\dd_{H}a'_{1}} are well defined. This gives exactly first order completeness of the maximal prolongation of $\Gamma_{A}$. The map $\Delta_{\Gamma_{A}}$ is well defined, since $\Gamma_{A}$ is a right $H-$covariant first order differential calculus.  Moreover, by assumption the map $$\ver_{BM}:\Gamma_{A}\rightarrow A\otimes \Gamma_{H},\quad a\dd_{A}(a')\mapsto aa'_{0}\otimes S(a'_{1})\dd_{H}a'_{2}=a_{0}a'_{0}\otimes S(a_{1}a'_{1})a_{2}\dd_{H}a'_{2}$$ is well defined, and thus $$(\id\otimes \triangleright)\circ (\Delta_{A}\otimes \id)\circ \ver_{BM}:\Gamma_{A}\rightarrow H\otimes \Gamma_{H}$$ is well defined. We have that the map $\ver$ is well defined, since \eqa{\nonumber(\id\otimes \triangleright)\circ (\Delta_{A}\otimes \id)\circ \ver_{BM}(a\dd_{A}a') & =   (\id\otimes \triangleright)\circ (\Delta_{A}\otimes \id)(a_{0}a'_{0}\otimes S(a_{1}a'_{1})a_{2}\dd_{H}a'_{2}) \\ & = (\id \otimes \triangleright)(a_{0}a'_{0}\otimes a_{1}a'_{1}\otimes S(a_{2}a'_{2})a_{3}\dd_{h}a'_{3})  \\ & = a_{0}a'_{0}\otimes \triangleright(a_{1}a'_{1}\otimes S(a_{2}a'_{2})a_{3}\dd_{H}a'_{3})\\ & = a_{0}a'_{0}\otimes a_{1}\dd_{H}a'_{1} \\ & = \ver(a\dd_{A}a').   } 
	
	For the other point, let us consider a first order complete differential calculus $\Omega^{\bullet}(A)$, and let $\Omega^{\bullet}(H)$ be the corresponding differential calculus on $H$. In particular $\Omega^{\bullet}(A)$ is right $H-$covariant and $\Omega^{\bullet}(H)$ is bicovariant. The vertical map $\ver:\Omega^{1}(A)\rightarrow A\otimes \Omega^{1}(A)$ is well defined, and so \eqa{\ver_{BM}=(\id\otimes \triangleright)\circ (\id\otimes S\otimes \id)\circ (\Delta_{A}\otimes \id)\circ \ver} is well defined.  Moreover, we have that $\pi_{\ver}|_{\Omega^{1}(A)}:\Omega^{1}(A)\rightarrow A\otimes \Lambda^{1}$ explicitly reads $$\pi_{\ver}(a\dd_{A}a') = a_{0}a_{1}\otimes S(a_{1}a'_{1})a_{2}\dd a'_{2},$$ i.e. it coincides with $\ver_{BM}(a\dd_{A}a')$. This means $\ver_{BM}$ is surjective, since $\pi_{\ver}$ is. The morphism sending horizontal forms into total space form is automatically injective, therefore it is only left to show that $\ker\ver_{BM}=A\dd_{A}(B)A$. Since $\hor^{1}(A)=\ker\pi_{\ver}=\ker\ver_{BM}$ this follows if and only if horizontal 1-forms of the \DJ ur\dj evi\'c' and 	Brzezi\'nski-Majid approaches are equal. \qed  
\end{prop}


\section{Examples}\label{Examples}

In this section we provide some non trivial examples of quantum principal bundles and complete differential calculi. 	Let $\Omega^{\bullet}(A)$ be a complete differential calculus over $A$. The morphism $\Delta_{A}^{\wedge}:\Omega^{\bullet}(A)\rightarrow \Omega^{\bullet}(A)\otimes \Omega^{\bullet}(H)$ maps elements of fixed degree $k$ in the direct sum decomposition of $\Omega^{\bullet}(A)$ as $\Delta_{A}^{\wedge}:\Omega^{k}(A)\rightarrow \Omega^{s}(A)\otimes \Omega^{r}(H)$, with $r+s=k$. 
\begin{prop}\label{decomposition of vertical maps}
Let $\Omega^\bullet(A)$ be a complete calculus.
Let us write 
$$
\ver^{k,\ell}:=(\pi_A^k\otimes\pi_H^\ell)\circ\Delta_A^\wedge|_{\Omega^{k+\ell}(A)}\colon\Omega^{k+\ell}(A)\to\Omega^k(A)\otimes\Omega^\ell(H)
$$
for the graded components of $\Delta_A^\wedge$.
Then
$$
\ver^{k,\ell}(\omega\wedge\eta)
=\sum_{m=0,\dots,|\omega|}\ver^{m , |\omega| - m}(\omega)\ver^{k - m , |\eta| - k + m}(\eta)
$$
for all $\omega,\eta\in\Omega^\bullet(A)$ such that $|\omega|+|\eta|=k+\ell$.
\proof Let us consider $\omega=a^{0}\dd a^{1}\wedge\cdots \wedge \dd a^{r}$ and $\eta=b^{0}\dd b^{1}\wedge \cdots \wedge\dd b^{s}$  in $\Omega^{\bullet}(A).$ We write $\dd_{\otimes}=(\id\otimes \dd + \dd \otimes \id)$ as a shorthand notation for this proof. On the left hand side we have\eqa{ \nonumber\ver^{k,\ell}(\omega\wedge \eta) & = ( \pi_{k}\otimes \pi_{\ell})\Delta_{A}^{\wedge}(\omega\wedge \eta) \\ & = (\pi_{k}\otimes \pi_{\ell}) (\Delta_{A}(a^{0})\dd_{\otimes}\Delta_{A}(a^{1})\wedge\cdots \wedge \dd_{\otimes}\Delta_{A}(a^{r})\\ & \wedge \Delta_{A}(b^{0})\dd_{\otimes}\Delta_{A}(b^{1}) \wedge\cdots \wedge \dd_{\otimes}\Delta_{A}(b^{r})). } On the right hand side we have \eqa{\nonumber\sum_{m=0,\dots,|\omega|}\ver^{m,|\omega|-m}(\omega)&\ver^{k-m,|\eta|-k+m}(\eta) \\ & =\sum_{m=0,\dots,|\omega|}(\pi^{m} \otimes \pi^{|\omega|-m})(\Delta_{A}(a^{0})\dd_{\otimes}\Delta_{A}(a^{1})\wedge\cdots \dd_{\otimes}\Delta_{A}(a^{r})) \\ & \quad \quad \quad  (\pi^{k-m}\otimes \pi^{|\eta|-k+m})(\Delta_{A}(b^{0})\dd_{\otimes}\Delta_{A}(b^{1})\wedge\cdots \dd_{\otimes}\Delta_{A}(b^{s})).} Since the projection maps fix the degree of left and right tensor factors of the above formulas we have the claim.

\qed 
\end{prop}
 \begin{exmp}
	Let $H$ be a Hopf algebra and let $(\Gamma,\dd)$ be a first order differential calculus on $H$. Consider the maximal prolongation of $(\Gamma,\dd)$. The subalgebra of coinvariant elements $B=H^{coH}$ is made of elements that are invariant under the coaction of $H$ on itself, i.e. the coproduct $\Delta:H\rightarrow H\otimes H$. This subalgebra is naturally isomorphic to $\B{k}$. Therefore we consider the Hopf-Galois map $\chi:H\otimes H \rightarrow H\otimes H$ sending $h\otimes g \mapsto h g_{1}\otimes g_{2}$, where no balanced tensor product is required after the established isomorphism $B\cong \B{k}$. This map is a bijection and in fact we can provide explicitely the inverse $\chi^{-1}:H\otimes H \rightarrow H \otimes H$ sending $h\otimes g \mapsto hS(g_{1})\otimes h_{2}$, indeed \eqa{ \chi\circ\chi^{-1}(h\otimes g) = \chi (hS(g_{1})\otimes g_{2}) = hS(g_{1}) g_{2}\otimes g_{3} = h \otimes g ,} \eqa{\chi^{-1}\circ \chi(h\otimes g)=\chi^{-1}(hg_{1}\otimes g_{2}) = hg_{1}S(g_{2})\otimes g_{3} = h\otimes g.} Therefore we have a quantum principal bundle.\qed 
\end{exmp} 
\begin{exmp}
	Consider the group algebra $H=\B{C}[\B{Z}]$. $H$ is generated by an element $g$, i.e. every element of $H$ can be written as $\alpha_{k}g^{k}$, where $\alpha_{k}\in \B{C}$ and $k\in \B{Z}$. The coproduct, counit and antipode read \eqa{ \Delta(g)=g\otimes g, \quad \epsilon(g)=1, \quad S(g)=g^{-1}. }  Consider the first order differential calculus $\Gamma$ over $H$ defined as $\Gamma= \span_{H}\{\dd g\}=H\dd g,$ with $(g^{n}\dd g)g :=q g^{n+1}\dd g$, where $q$ is not a root of unity. The differential $\dd : H\rightarrow \Gamma$ maps $\dd(g):=\dd g$ and in general \eqa{ & \dd(g^{n}) :=(1+q+\dots + q^{n-1})g^{n-1} \dd g, \quad &&\text{if} \ n\geq 0, \\ & \dd (g^{n}) :=-(q^{n}+\dots + q^{-1}) g^{n-1} \dd g, \quad && \text{if} \ n<0, \\ &\dd f(t):= \frac{(f(qt)-f(t))}{(t(q-1))},\quad &&\text{for a rational function}. } In particular $\dd (g^{-1}) =-g^{-1} (\dd g) g^{-1}$  much like a chain rule. 
	We define right and left $H-$coactions on $\Gamma$ as \eqa{  &\Delta_{\Gamma}: \Gamma\rightarrow \Gamma\otimes H, \quad \quad  && {}_\Gamma\Delta:\Gamma\rightarrow H\otimes \Gamma. \\ &  g^{n}\dd g \mapsto g^{n}\dd g\otimes g^{n+1}, \quad \quad &&g^{n}\dd g \mapsto g^{n+1}\otimes g^{n}\dd g.  } Since $\Delta$ extends to the first order calculus we have that $(\Gamma,\dd)$ is a bicovariant first order differential calculus over $H$.
	
	Moreover, since $ \dd(g^{-1})+q^{-1} g^{-2} \dd g=0,$ we find \eqa{ \nonumber(\dd g^{-1}+ q^{-1}g^{-2} \dd g)\otimes \dd g & =  q^{-1}(\dd g^{-1} g^{-1}+ g^{-1} \dd g^{-1})\otimes \dd g \\ & = -q^{-1}(g^{-1} \dd g g^{-2}+ g^{-2}\dd g g^{-1})\otimes \dd g \\ & =- q^{-2}(1+ q^{-1}) g^{-3} \dd g \otimes \dd g. } The last expression vanishes on the quotient defined by the maximal prolongation. Therefore on this quotient we have $[\dd g \otimes \dd g]=0$, as we assumed $q$ not a root of unity. Accordingly we have no non-zero $k-$forms for $k>1$, and thus a complete differential calculus. \qed \end{exmp}

\subsection{The noncommutative algebraic 2-torus}
	In the following example we will discuss the quantum principal bundle given by the Hopf algebra $H=\M{O}_{q}(\text{U}(1))=\B{C}[t,t^{-1}]$ and right $H-$comodule algebra $A=\M{O}_{\theta}(\B{T}^{2})=\B{C}[u,u^{-1},v,v^{-1}]/\langle uv-e^{i\theta}vu\rangle$,  where $\theta\in \B{R}$. This algebra $A$ is known as the non-commutative algebraic 2-torus. The corresponding Hopf algebra structure of $H$ is given by coproduct $\Delta(t)=t\otimes t$, counit $\epsilon(t)=1$ and antipode $S(t)=t^{-1}$. Let $\Delta_{A}:A\rightarrow A\otimes H$ be the map assigning \eqa{ \begin{pmatrix}
		 u \\ v 
	\end{pmatrix} \mapsto \begin{pmatrix}
		u \\ v
	\end{pmatrix} \otimes \begin{pmatrix}
		 t \\ t^{-1}
	\end{pmatrix}.}  This is a right $H-$coaction on $A$, indeed \eqa{\label{coaction on the torus 1}(\Delta_{A}\otimes \id)\circ \Delta_{A}\begin{pmatrix}
		u \\ v
	\end{pmatrix} =  (\Delta_{A}\otimes \id)\begin{pmatrix}
		u \\ v
	\end{pmatrix} \otimes \begin{pmatrix}
		 t \\ t^{-1} \end{pmatrix} = \begin{pmatrix}
		u \\ v
	\end{pmatrix} \otimes \begin{pmatrix}
		 t \\ t^{-1} \end{pmatrix}\otimes \begin{pmatrix}
		 t \\ t^{-1} \end{pmatrix},} whereas \eqa{ \label{coaction on the torus 2}(\id \otimes \Delta)\circ \Delta_{A} \begin{pmatrix}
		 	u \\ v
		 \end{pmatrix} = (\id \otimes \Delta)\begin{pmatrix}
		 	u \\ v 
		 \end{pmatrix}\otimes \begin{pmatrix}
		 	 t \\ t^{-1} 
		 \end{pmatrix}= \begin{pmatrix}
		 	u \\ v
		 \end{pmatrix}\otimes \begin{pmatrix}
		 	 t \\ t^{-1}
		 \end{pmatrix}\otimes \begin{pmatrix}
		 	 t \\ t^{-1}
		 \end{pmatrix}.} Moreover \eqa{\label{coaction on the torus 3}(\id \otimes \epsilon)\circ \Delta_{A}\begin{pmatrix}
		 	u \\ v
		 \end{pmatrix} = (\id\otimes \epsilon)\begin{pmatrix}
		 	 u \\ v
		 \end{pmatrix} \otimes \begin{pmatrix}
		 	 t\\ t^{-1}
		 \end{pmatrix} =  \begin{pmatrix}
		 	u \\ v
		 \end{pmatrix}\otimes 1 \cong \begin{pmatrix}
		 	 u \\ v 
		 \end{pmatrix}.} The space of right $H-$coinvariant forms is given by $$B:=A^{coH}=\span_{\B{C}}\{(uv)^{k},\ k\in \B{Z}\}.$$ 
		 \begin{prop}
		The noncommutative algebraic 2-torus under the coaction of $\M{O}(\text{U}(1))$ is a cleft extension $B\subseteq A$. In particular it is a quantum principal bundle.
		\proof  Let us consider $j:H\rightarrow A$ sending $$\begin{pmatrix}
			 t^{k}\\ t^{-k} 
		\end{pmatrix} \mapsto \begin{pmatrix}
			 u^{k} \\ v^{k} 
		\end{pmatrix},$$ for all $k\geq 0.$  We show that $j:H\rightarrow A$ is a convolution invertible morphism of right $H-$comodules. \eqa{ & \Delta_{A}(j(t)) = \Delta_{A}(u) = u \otimes t =(j\otimes \id)\Delta(t),\\   & \Delta_{A}(j(t^{-1})) = \Delta_{A}(v) = v \otimes t^{-1} =(j\otimes \id)\Delta(t^{-1}).} 
		 We provide an explicit inverse for the map $j:H\rightarrow A$, that is $j^{-1}:H\rightarrow A$ sending $$ \begin{pmatrix}
		 	t^{k} \\ t^{-k}
		 \end{pmatrix} \mapsto \begin{pmatrix}
		 	 u^{-k} \\ v^{-k}
		 \end{pmatrix};$$ indeed \eqa{ \nonumber j* j^{-1}(t^{k}) & = \mu \circ (j\otimes j^{-1})\circ \Delta(t^{k}) \\ & = j(t^{k})j^{-1}(t^{k})\\ & = u^{k}u^{-k}\\ & = 1 , \\  j* j^{-1}(t^{-k}) & = \mu \circ (j\otimes j^{-1})\circ \Delta(t^{-k}) \\ & = j(t^{-k})j^{-1}(t^{-k}) \\ & = v^{k}v^{-k}\\ & = 1, } and similarly $j^{-1}* j = 1$ is easily verified. 
		 
		 Since every cleft extension is in particular a faithfully flat Hopf-Galois extension \cite{BRZ2} we have the claim. \qed   
		 \end{prop} We now investigate in depth the differential calculi over $H=\M{O}(\text{U}(1))$ and $A=\M{O}_{\theta}(\B{T}^{2})$. 
		 	 For the differential structure over $H$, we consider a bicovariant first order differential calculus $(\Omega^{1}(H)=\span_{H}\{\dd t\},\dd)$ with bimodule relations \eqa{t\dd t & = q^{\alpha}\dd tt , \\  t^{-1}\dd t & = q^{-\alpha}\dd t t^{-1}.} Here $q$ is a understood as a deformation parameter and $\alpha\in \B{R}$.  Differentiating these relations we find \eqa{ \nonumber\dd t\wedge \dd t = - q^{\alpha} \dd t \wedge \dd t.} Accordingly there are no non-zero higher order forms. 
		 
	 Proceeding with the calculus on $A$, we define a bicovariant first order differential calculus $(\Omega^{1}(A)=\span_{A}\{\dd u,\dd v\},\dd)$ via relations \eqa{\label{relations for U(1)} \dd u u = u \dd u,\quad \dd v v = v \dd v,\quad \dd uv=e^{-i\theta} v\dd u, \quad \dd vu = e^{i\theta}u\dd v. } Differentiating, we find \eqa{ & \dd u \wedge\dd v + e^{-i\theta}\dd v\wedge\dd u=0, \\ &\dd u \wedge \dd u=0, \\ &\dd v\wedge \dd v = 0,  } and clearly higher order forms vanish.
	 
	 Let us consider $\alpha=0$ in equation \eqref{relations for U(1)}, so that $q^{\alpha}=1$.
	
		 \begin{prop}\label{torus calculus complete 1}
		 	$\Omega^{1}(A):=\span_{A}\{\dd u, \dd v\}$ is a right $H-$covariant differential calculus on the noncommutative algebraic 2-torus. Moreover the map $\ver:=\ver^{0,1}$ is well defined.   
		 	\proof We show the coaction $\Delta_{A}:A\rightarrow A\otimes H$ lifts to a morphism  $\Delta_{A}^{1}:\Omega^{1}(A)\rightarrow \Omega^{1}(A)\otimes \Omega^{1}(H)$.  We have $$\Delta_{A}^{1}=\Delta_{\Omega^{1}(A)}+\ver:\Omega^{1}(A)\rightarrow (A\otimes \Omega^{1}(H)) \oplus (\Omega^{1}(A)\otimes H). $$ To properly define this map we need \eqa{ &  \Delta_{\Omega^{1}(A)}(\dd u) = (\dd \otimes \id)\Delta_{A}(u)=\dd u \otimes t,\\  & \Delta_{\Omega^{1}(A)}(\dd v)=(\dd \otimes \id)\Delta_{A}(v)= \dd v \otimes t^{-1},}  and then we extend $\Delta_{\Omega^1(A)}$ to $\Omega^1(A)$ by
\eqa{\Delta_{\Omega^1(A)}(a\dd u+a'\dd v)=\Delta_A(a)\Delta_{\Omega^1(A)}(\dd u)+\Delta_A(a')\Delta_{\Omega^1(A)}(\dd v).} 
	 Moreover $\Delta_{\Omega^{1}(A)}$ is well defined, since  \eqa{ \nonumber\Delta_{\Omega^{1}(A)}(\dd u. v) & =  \Delta_{\Omega^{1}(A)}(\dd u) \Delta_{A}(v) \\ & = (\dd u \otimes t)(v\otimes t^{-1}) \\ & = \dd uv\otimes 1,  }  \eqa{ \nonumber\Delta_{\Omega^{1}(A)}(e^{-i\theta}v\dd u) & = e^{-i\theta} \Delta_{A}(v)\Delta_{\Omega^{1}(A)}(\dd u) \\ & = e^{-i\theta}(v\otimes t^{-1})(\dd u \otimes t) \\ & = e^{-i\theta} v\dd u \otimes 1 \\ & = \Delta_{\Omega^{1}(A)}(\dd uv), } 
	 and similarly $\Delta_{\Omega^{1}(A)}( \dd v.u)=\Delta_{\Omega^{1}(A)}(e^{i\theta}u.\dd v)$. To read more explicitly $\ver:\Omega^{1}(A)\rightarrow A\otimes \Omega^{1}(H)$ we consider $a,a'\in A$ and develop \eqa{ \Delta_{A}^{\wedge}(a \dd a') & = \Delta_{A}(a)(\dd \otimes \id + \id \otimes \dd) \Delta_{A}(a') \\ & = a_{0}\dd a'_{0} \otimes a_{1}a'_{1} + a_{0}a'_{0}\otimes a_{1}\dd a'_{1}.} Accordingly we define on the basis elements\eqa{ \ver(\dd u)=u \otimes \dd t, \quad \ver(\dd v)=v\otimes \dd t^{-1}.  } Therefore we have \eqa{ \nonumber\ver( \dd u v) & = \ver(\dd u)\Delta_{A}(v) \\ & = (u\otimes \dd t)(v\otimes t^{-1}) \\ & = uv \otimes \dd t t^{-1}  } \eqa{\nonumber\ver(e^{-i\theta}v \dd u) & = e^{-i\theta} \Delta_{A}(v)\ver(\dd u) \\ & = e^{-i\theta} (v\otimes t^{-1})(u\otimes \dd t) \\ & = e^{-i\theta}vu\otimes t^{-1} \dd t \\ & = e^{-i\theta} vu\otimes q^{-\alpha}\dd t t^{-1}. } Similarly \eqa{\nonumber\ver(\dd v u) & = \ver(\dd v)\Delta_{A}(u) \\ & = (v\otimes \dd t^{-1})(u\otimes t) \\ & = vu \otimes \dd t^{-1}t \\ & = - vu \otimes t^{-2} \dd tt \\ & = -vu\otimes  t^{-2}t\dd t \\ & = -vu\otimes  t^{-1}\dd t ,} \eqa{\nonumber\ver(e^{i\theta}u\dd v) & = e^{i\theta} \Delta_{A}(u)\ver(\dd v) \\ & = e^{i\theta}(u\otimes t)(v\otimes \dd t^{-1}) \\ & = e^{i\theta}uv\otimes t\dd t^{-1} \\ & =- e^{i\theta}uv\otimes t t^{-2} \dd t \\ & = - e^{i\theta}uv\otimes t^{-1} \dd t. } 
	  The map $\ver:\Omega^{1}(A)\rightarrow A\otimes \Omega^{1}(H)$ is thus well defined along the relations between the generators of the first order differential calculus over $A$. We extend to $\Omega^{1}(A)$ as \eqa{\ver(a \dd u + a' \dd v) =\Delta_A(a)\ver(\dd u)+\Delta_A(a')\ver(\dd v).} \qed \end{prop}
	  
		 \begin{prop}\label{torus calculus complete 2}		 	$\Omega^{2}(A):=\span_{A}\{\dd u\wedge\dd v\}$ is a right $H-$covariant differential calculus on the noncommutative algebraic 2-torus. Moreover maps $\ver^{1,1}$ and $\ver^{0,2}$ are well defined. 
		 	\proof We show that $\Delta_{A}:A\rightarrow A\otimes H$ extends to a morphism  $$\Delta_{A}^{2}:\Omega^{2}(A)\rightarrow \underbrace{(\Omega^{2}(A)\otimes H )}_{\Delta_{\Omega^{2}(A)}}\oplus \underbrace{(\Omega^{1}(A)\otimes \Omega^{1}(H))}_{\ver^{1,1}}\oplus \underbrace{(A\otimes \Omega^{2}(H))}_{\ver^{0,2}}.$$  Let us consider $\omega=a \dd a'\wedge\dd a''\in \Omega^{2}(A)$ for $a,a',a''\in A$. We have \eqa{\nonumber\Delta_{A}^{\wedge}(a\dd a' \wedge \dd a'') & = \Delta_{A}(a) (\dd \otimes \id + \id \otimes \dd) \Delta_{A}(a') (\id\otimes \dd + \dd \otimes \id)\Delta_{A}(a'') \\ & = (a_{0}\otimes a_{1})(\dd a'_{0}\otimes a_{1} + a'_{0}\otimes \dd a'_{1})(\dd a''_{0}\otimes a''_{1} + a''_{0}\otimes \dd a''_{1}) \\ & = (a_{0}\otimes a_{1}) ( \dd a'_{0} \wedge \dd a''_{0} \otimes a'_{1} a''_{1} + \dd a'_{0} a''_{0} \otimes a'_{1}\dd a''_{1} - a'_{0} \dd a''_{0}\otimes \dd a'_{1} a''_{1} ), } where we exploited that there are no non-zero two forms on $H$. Accordingly we consider \begin{equation}\begin{aligned}
 	 \Delta_{\Omega^{2}(A)}(\mathrm{d}u\wedge \mathrm{d}v ) & := \mathrm{d}u\wedge \mathrm{d}v, \\ 
 	\mathrm{ver}^{1,1}(\mathrm{d}u\wedge \mathrm{d}v) & := \mathrm{d}uv\otimes t\mathrm{d}(t^{-1}) - u \mathrm{d} v\otimes \mathrm{d} tt^{-1},
 	\end{aligned}
 \end{equation} and $\ver^{0,2}$ to be the zero map.  These maps are well defined, since  \eqa{ \nonumber\Delta_{\Omega^{2}(A)}(\dd u \wedge \dd v) & = \Delta_{\Omega^{1}(A)}(\dd u) \Delta_{\Omega^{1}(A)}(\dd v) \\ & = (\dd u \otimes t) (\dd v \otimes t^{-1}) \\ & = (\dd u \wedge\dd v) \otimes 1. } 
		 	\eqa{\nonumber \ver^{1,1}(\dd u \wedge \dd v) & = \dd u v \otimes t \dd t^{-1}- u \dd v \otimes \dd t t^{-1} \\ & = -\dd u v \otimes t^{-1}\dd t - u\dd v \otimes \dd tt^{-1} \\ & = -(\dd u v + u  \dd v )\otimes \dd t t^{-1} \\ & = - \dd (uv)\otimes \dd t t^{-1};} 
		 	\eqa{ \nonumber\ver^{1,1}(\dd v\wedge \dd u) & = \dd v u \otimes t^{-1} \dd t - v\dd u \otimes \dd t^{-1} t \\ & = \dd v  u \otimes t^{-1} \dd t + v \dd u \otimes t^{-2} \dd t t  \\  & = \dd v u \otimes \dd t t^{-1} + v \dd u \otimes t^{-1} \dd t \\ & = \dd v u\otimes \dd t t^{-1} + v\dd u\otimes \dd tt^{-1} \\ & = \dd(vu)\otimes \dd t t^{-1},} 
		 	and so \eqa{\nonumber \ver^{1,1}(\dd u \wedge \dd v)+e^{-i\theta} \ver^{1,1}(\dd v\wedge \dd u) & = - \dd(uv)\otimes \dd t t^{-1} +e^{-i\theta}\dd(vu)\otimes \dd t t^{-1} \\ &  -\dd (uv)\otimes \dd t t^{-1} + \dd (e^{-i\theta}vu)\otimes \dd t t^{-1} \\ & = -\dd (uv)\otimes \dd t t^{-1} + \dd(uv)\otimes \dd t t^{-1}. } 
		 	We extend to $\Omega^{2}(A)$ as \begin{equation}
 	\begin{aligned}
 		\Delta_{\Omega^{2}(A)}(a\mathrm{d}u\wedge \mathrm{d}v)&= \Delta_{A}(a)\Delta_{\Omega^{2}(A)}(\mathrm{d} u\wedge \mathrm{d} v), \\  \mathrm{ver}^{1,1}(a\mathrm{d}u\wedge\mathrm{d}v) &= \Delta_{A}(a) \mathrm{ver}^{1,1}(\mathrm{d}u\wedge \mathrm{d} v), \\  \mathrm{ver}^{0,2}(a\mathrm{d}u\wedge \mathrm{d}v)&=0.
 	\end{aligned}
 \end{equation}\qed 
		 \end{prop}
		 \begin{thm}
		 	The differential calculus $\Omega^{\bullet}(A)$ on the noncommutative algebraic 2-torus is complete. \proof This follows directly from Proposition \ref{torus calculus complete 1}-\ref{torus calculus   complete 2}\qed 
		 \end{thm}
		  \begin{prop}
		  The base space 1-forms are generated by $B$. In other words $B\dd B=\Omega^{1}(B)$. 
		  	\proof According to lemma \ref{inclusion of BdB into coinvariant forms} the non trivial statement is that $B\dd B \supseteq \Omega^{1}(B)$. A generic element in $\omega\in \Omega^{1}(B)$ can be written as $\omega= \alpha_{k\ell}u^{k}v^{\ell}\dd u  + \beta_{mn}u^{m}v^{n}\dd v,$ for $\alpha_{k\ell},\beta_{mn}\in \mathbb{C}$. We have \eqa{\nonumber\Delta_{A}^{\wedge}&(\alpha_{k\ell}u^{k}v^{\ell} + \beta_{mn}u^{m}v^{n}) =  \alpha_{k\ell}\Delta_{A}(u^{k}v^{\ell})\Delta_{A}^{1}(\dd u) +  \beta_{mn}\Delta_{A}(u^{m}v^{n})\Delta_{A}^{1}(\dd v) \\ & = \alpha_{k\ell}(u^{k}v^{\ell})\otimes t^{k-\ell}(\dd \otimes \id + \id \otimes \dd)\Delta_{A}(u) +\beta_{mn}(u^{m}v^{n})\otimes t^{m-n}(\dd \otimes \id + \id \otimes \dd)\Delta_{A}(v) \\ & = \alpha_{k\ell}u^{k}v^{\ell}( \dd u \otimes t^{k-\ell +1} + u \otimes t^{k-\ell}\dd t) + \beta_{mn}u^{m}v^{n}( \dd v\otimes t^{m-n -1} + v \otimes t^{m-n}\dd t^{-1})  \\ & = \alpha_{k\ell}u^{k}v^{\ell}\dd u \otimes t^{k-\ell +1} + \beta_{mn}u^{m}v^{n} \dd v\otimes t^{m-n -1} \\ & +\alpha_{k\ell}u^{k}v^{\ell}u \otimes t^{k-\ell}\dd t - \beta_{mn}u^{m}v^{n}v \otimes t^{m-n-2}\dd t \\ & = (\alpha_{k\ell}u^{k}v^{\ell}\dd u  + \beta_{mn}u^{m}v^{n}\dd v) \otimes 1 }if and only if \eqa{ \alpha_{k\ell}u^{k}v^{\ell}\dd u \otimes t^{k-\ell +1} + \beta_{mn}u^{m}v^{n} \dd v\otimes t^{m-n -1} & = (\alpha_{k\ell}u^{k}v^{\ell}\dd u  + \beta_{mn}u^{m}v^{n}\dd v) \otimes 1, \\  \alpha_{k\ell}u^{k}v^{\ell}u \otimes t^{k-\ell}\dd t - \beta_{mn}u^{m}v^{n}v \otimes t^{m-n-2}\dd t & =0.}  The first equation tells we must have $k-\ell +1 =0$ (or $\beta_{mn}=0$) and $m-n-1=0$ (or $\alpha_{k\ell}=0$), which leads to \eqa{\alpha_{k\ell}u^{k}v^{\ell}\dd u \otimes t^{k-\ell +1} + \beta_{mn}u^{m}v^{n} \dd v\otimes t^{m-n -1} = \alpha_{k,k+1}u^{k}v^{k+1}\dd u\otimes 1 + \beta_{n,n+1} u^{n+1}v^{n} \dd v \otimes 1 .} The second equation reads \eqa{ \nonumber\alpha_{k\ell}u^{k}v^{\ell}u \otimes t^{k-\ell}\dd t & - \beta_{mn}u^{m}v^{n}v \otimes t^{m-n-2}\dd t  \\ & = \alpha_{k,k+1}u^{k}v^{k+1}u \otimes t^{-1}\dd t - \beta_{n+1,n}u^{n+1}v^{n}v\otimes t^{-1}\dd t \\ & =(\alpha_{k,k+1}e^{i(k+1)\theta}u^{k+1}v^{k+1}- \beta_{n,n+1}u^{n+1}v^{n+1})\otimes t^{-1}\dd t, } which is zero if and only if \eqa{\alpha_{k,k+1}e^{i(k+1)\theta}u^{k+1}v^{k+1}- \beta_{n,n+1}u^{n+1}v^{n+1}=0,}  so we must have $n=k$, and accordingly $\alpha_{k,k+1}e^{i(k+1)\theta}=\beta_{k,k+1}$.  
		  	Therefore a general element $\omega\in \Omega^{1}(B)$ reads  \eqa{\nonumber\alpha_{k\ell}u^{k}v^{\ell}\dd u  + \beta_{mn}u^{m}v^{n}\dd v & = \alpha_{k,k+1} u^{k}v^{k+1}\dd u + \beta_{k,k+1}u^{k+1}v^{k}\dd v \\ & = \alpha_{k,k+1}(u^{k}v^{k}v \dd u   +e^{i(k+1)\theta} u^{k}u v^{k}\dd v ) \\ & =  \alpha_{k,k+1}u^{k}v^{k}(v\dd u + e^{i(k+1)\theta}u^{k}v^{k}\dd v u e^{-i(k+1)\theta}) \\ & =\alpha_{k,k+1} u^{k}v^{k}(v\dd u + \dd v u) \\ & = \alpha_{k,k+1} (uv)^{k} \dd (uv),} for any $k\in \B{Z}$. Since $B:=A^{coH}=\span_{\B{C}}\{(uv)^{k}| \ k\in \B{Z}\}$ we have the thesis.  \qed 
		  	\end{prop}

\subsection{Quantum Hopf fibration and the Podleś sphere}
Let us consider the Hopf algebra $H=\M{\text{O}}(\text{U}(1))$, and let us fix $\alpha=2$ in Equation \ref{relations for U(1)}. We introduce $\M{O}_{q}(\text{SU}(2))$ as the free algebra generated by elements $\alpha,\beta,\gamma,\delta$ modulo relations \eqa{& \beta\alpha=q\alpha\beta, \quad \gamma\alpha=q\alpha\gamma, \quad \delta\beta=q\beta\delta, \quad \delta\gamma=q\gamma\delta \\  & \gamma\beta=\beta\gamma, \quad \delta\alpha-\alpha\delta =(q-q^{-1})\beta \gamma, \quad   \alpha\delta-q^{-1}\beta \gamma =1.}

\begin{prop} $A=\M{O}_{q}(\text{SU}(2))$ is a right $H-$comodule algebra under the right $H-$coaction \eqa{\Delta_{A}:A\rightarrow A\otimes H, \quad \begin{pmatrix}
	\alpha & \beta \\ \gamma & \delta 
\end{pmatrix} \mapsto \begin{pmatrix}
	 \alpha & \beta \\ \gamma & \delta
\end{pmatrix}\otimes \begin{pmatrix}
	 t & 0 \\ 0 & t^{-1}
\end{pmatrix}= \begin{pmatrix}
	\alpha\otimes t & \beta \otimes t^{-1} \\ \gamma\otimes t & \delta \otimes t^{-1} 
\end{pmatrix}.}\proof We check that $\Delta_{A}$ satisfies the axioms of a right $H-$coaction on $A$, indeed \eqa{ \nonumber (\Delta_{A}\otimes \id)\circ \Delta_{A}\begin{pmatrix}
	\alpha & \beta \\ \gamma & \delta
\end{pmatrix} & = \Delta_{A}\begin{pmatrix}
	\alpha & \beta \\ \gamma & \delta 
\end{pmatrix}\otimes \begin{pmatrix}
	 t & 0 \\ 0 & t^{-1}
\end{pmatrix} \\ & = \begin{pmatrix}
	\alpha & \beta \\ \gamma & \delta 
\end{pmatrix}\otimes \begin{pmatrix}
	 t & 0 \\ 0 & t^{-1}\end{pmatrix}\otimes \begin{pmatrix}
	 t & 0 \\ 0 & t^{-1}\end{pmatrix}\\ & = \begin{pmatrix}
	\alpha & \beta \\ \gamma & \delta 
\end{pmatrix}\otimes \Delta\begin{pmatrix}
	 t & 0 \\ 0 & t^{-1}\end{pmatrix} \\ & = (\id \otimes \Delta) \circ \Delta_{A}\begin{pmatrix}
	\alpha & \beta \\ \gamma & \delta 
\end{pmatrix};  }  moreover
  \eqa{ (\id\otimes \epsilon)\circ \Delta_{A}\begin{pmatrix}
	\alpha & \beta \\ \gamma & \delta 
\end{pmatrix} = (\id\otimes \epsilon)\circ \begin{pmatrix}
	\alpha & \beta \\ \gamma & \delta 
\end{pmatrix} \otimes \begin{pmatrix}
	 t & 0 \\ 0 & t^{-1} 
\end{pmatrix} = \begin{pmatrix}
	\alpha & \beta \\ \gamma & \delta 
\end{pmatrix}.}  \qed \end{prop}
\begin{notation}
With $|\cdot|$ we indicate the degree of an element in $A$; $|f|$ is explicitly defind by the relation $\Delta_{A}(f)=f\otimes t^{|f|},$ for any $f\in \{\alpha,\beta,\gamma,\delta\}$.   
\end{notation}
\begin{rmk} The  subalgebra $B:=A^{coH}$ of coinvariant elements under the right $H-$coaction $\Delta_{A}$ is known as the Podleś sphere \cite{Podles1989DifferentialCO,Podles:1987wd}. By a quick calculation it is easy to deduce generators of $B$ must be of the form $fg$, where $f,g\in\{\alpha,\beta,\gamma,\delta\}$ and either $|f|=1$ and $|g|=-1$ or the opposite. Modulo commutation relations we have $3$ non equivalent generators for $B$ and we shall fix \eqa{ B_{+}=\alpha\beta, \quad B_{-}=\gamma\delta, \quad B_{0}=\gamma\beta.}   One finds \eqa{\nonumber B_{-}B_{0} & = \gamma\delta \gamma \beta = q\gamma \gamma \delta \beta \\ & = q^{2}\gamma\gamma \beta \delta = q^{2}\gamma\beta \gamma\delta \\ &= q^{2}B_{0}B_{-},} and similarly \eqa{ B_{-}B_{+} =q^{2}B_{0}(1-q^{2}B_{0}), \quad B_{+}B_{-} = B_{0}(1-B_{0}).}  \end{rmk}

\begin{prop}[\cite{BRZ2}]
	$B\subseteq A$ is a quantum principal bundle. \qed 
\end{prop}
We define a first order differential calculus $(\Omega^{1}(A),\dd)$ over $A$ as the free left $A-$module generated by \eqa{e^{+}=q^{-1}\alpha\dd\gamma-q^{-2}\gamma\dd \alpha. \quad e^{-}=\delta\dd \beta - q\beta \dd \delta, \quad e^{0}=\delta\dd \alpha-q\beta \dd \gamma,} with commutation relations \eqa{e^{\pm}f=q^{|f|}fe^{\pm},\quad e^{0}f=q^{2|f|}fe^{0},} where $f\in \{\alpha,\beta,\gamma,\delta\},$ and $|\alpha|=|\gamma|=1$, whereas $|\beta|=|\delta|=1$. 
 \begin{prop}
	$ \Omega^1(A)$ is  a right $H-$covariant differential calculus. Moreover the map $\ver=\ver^{0,1}$ is well-defined. \proof We need to show that $\Delta_{A}:A\rightarrow A\otimes H$ extends to a morphism $$\Delta_{A}^{1}=\Delta_{\Omega^{1}(A)}+\ver:\Omega^{1}(A)\rightarrow (\Omega^{1}(A)\otimes H) \oplus (A\otimes \Omega^{1}(H))$$restricting to the usual coaction on $A$. Such extension is completely determined by the generators. We want \eqa{\nonumber \Delta_{\Omega^{1}(A)}(e^{+}) & = \Delta_{\Omega^{1}(A)}(q^{-1}\alpha \dd \gamma - q^{-2} \gamma \dd \alpha) \\ & = q^{-1}\Delta_{A}(\alpha)\Delta_{\Omega^{1}(A)}(\dd\gamma)- q^{-2}\Delta_{A}(\gamma) \Delta_{\Omega^{1}(A)}(\dd\alpha) \\ & = q^{-1}(\alpha\otimes t)(\dd \otimes \id)(\gamma\otimes t) - q^{-2}(\gamma\otimes t)(\dd \alpha\otimes t) \\ & = e^{+}\otimes t^{2}} and similarly  $\Delta_{\Omega^{1}(A)}(e^{-})=e^{-}\otimes t^{-2}$. Moreover \eqa{ \nonumber\Delta_{\Omega^{1}(A)}(e^{0}) & = \Delta_{\Omega^{1}(A)}(\delta \dd \alpha -q\beta \dd\gamma) \\ & = \Delta_{A}(\delta)(\dd\otimes \id)\Delta_{A}(\alpha) - q \Delta_{A}(\beta)(\dd\otimes \id)\Delta_{A}(\gamma) \\ & = (\delta\otimes t^{-1})(\dd\alpha\otimes t) -q(\beta\otimes t^{-1})(\dd\gamma\otimes t) \\ & = (\delta\dd \alpha -q\beta\dd\gamma)\otimes 1 \\ & =e^{0}\otimes 1,} and so, on the generators of $\Omega^{1}(A)$, we define\eqa{  \Delta_{\Omega^{1}(A)}(e^{\pm}):=e^{\pm}\otimes t^{\pm 2},\quad \text{and} \quad \Delta_{\Omega^{1}(A)}(e^{0}):=e^{0}\otimes 1.} 

	Proceeding in a similar fashion for $\ver:\Omega^{1}(A)\rightarrow A\otimes \Omega^{1}(H)$ we need \eqa{\nonumber\ver(e^{0}) & =\ver(\delta \dd \alpha - q \beta \dd \gamma) \\ & = \delta\alpha \otimes t^{-1}\dd t - q \beta \gamma \otimes t^{-1}\dd t \\ & = (\delta \alpha - q\beta \gamma)\otimes t^{-1}\dd t;  } \eqa{\nonumber\ver(e^{+}) & = \ver(q^{-1}\alpha\dd \gamma - q^{-2} \gamma \dd \alpha) \\ & = q^{-1}\alpha\gamma \otimes t\dd t - q^{-2}\gamma \alpha \otimes t \dd t \\ & = (q^{-1}\alpha \gamma - q^{-2}\gamma\alpha) \otimes t \dd t\\ & = q^{-1}(\alpha\gamma -q^{-1}\gamma\alpha)\otimes t \dd t \\ & = 0 } 
	\eqa{\nonumber \ver(e^{-}) & = \ver(\delta \dd \beta -q \beta \dd \delta) \\ & = \delta \beta \otimes t^{-1}\dd t^{-1} - q \beta \delta\otimes t^{-1}\dd t^{-1} \\ & = (q\beta \delta - \delta \beta)\otimes t^{-1}\dd t^{-1} \\ & = 0, } where for $\ver(e^{\pm})$ we exploited the commutation relations. Accordingly, on the generators of $\Omega^{1}(A)$, we define \eqa{ \ver(e^{0}):= 1\otimes t^{-1}\dd t, \quad \ver(e^{\pm})=0.} With these definitions the map $\Delta_{A}^{1}:=\Delta_{\Omega^{1}(A)}+\ver$ extends correctly to $(\Omega^{1}(A),\dd)$ provided $\Delta_{\Omega^{1}(A)}$ and $\ver$ satisfy the commutation relations among $A$ and $\Omega^{1}(A)$. We have \eqa{\nonumber\Delta_{\Omega^{1}(A)}(e^{\pm}f) & = \Delta_{\Omega^{1}(A)} (e^{\pm})\Delta_{A}(f) \\ & =  (e^{\pm}\otimes t^{\pm 2})( f\otimes t^{|f|}) \\ & =  e^{\pm} f \otimes t^{\pm 2}t^{|f|} \\ & = q^{|f|} f e^{\pm}\otimes t^{\pm 2} t^{|f|}\\ & = q^{|f|} \Delta_{A}(f)\Delta_{\Omega^{1}(A)}(e^{\pm}) \\ & = \Delta_{\Omega^{1}(A)}(q^{f}f e^{\pm}),\\  \Delta_{\Omega^{1}(A)}(e^{0}f) & = \Delta_{\Omega^{1}(A)}(e^{0})\Delta_{A}(f) \\ & = (e^{0}\otimes 1) (f\otimes t^{|f|}) \\ & = e^{0}f \otimes t^{|f|} \\ & = q^{2|f|}f e^{0}\otimes t^{|f|} \\ & = q^{2|f|} \Delta_{A}(f)\Delta_{\Omega^{1}(A)}(e^{0}) \\ & = \Delta_{\Omega^{1}(A)}(q^{2|f|}fe^{0}), \\ \ver(e^{\pm}f) & = \ver(e^{\pm})\Delta_{A}(f) \\ & = 0 \\ & = q^{|f|}\Delta_{A}(f)\ver(e^{\pm}) \\ & = \ver(q^{|f|} f e^{\pm}),\\ \ver(e^{0}f) & = \ver(e^{0})\Delta_{A}(f) \\ & = (1\otimes t^{-1}\dd t)(f\otimes t^{|f|}) \\ & = f\otimes t^{-1}\dd t \ t^{|f|} \\ & =f\otimes t^{-1}\dd t \  t^{|f|} \\ & = q^{2|f|}f\otimes t^{|f|}t^{-1}\dd t ,}  and \eqa{\ver(q^{2|f|}fe^{0}) =q^{2|f|} f\otimes t^{|f|}t^{-1} \dd t. }  We extend to $\Omega^{1}(A)$ as \begin{equation}
	 	\begin{aligned}
	 		\Delta_{\Omega^{1}(A)}(ae^{+}+be^{-}+ce^{0}) & = \Delta_{A}(a)\Delta_{\Omega^{1}(A)}(e^{+})+ \Delta_{A}(b)\Delta_{\Omega^{1}(A)}(e^{-}) +\Delta_{A}(c)\Delta_{\Omega^{1}(A)}(e^{0}),\\ 
	 		\mathrm{ver}(ae^{+}+be^{-}+ce^{0}) & = \Delta_{A}(a)\mathrm{ver}(e^{+})+\Delta_{A}(b)\mathrm{ver}(e^{-})+\Delta_{A}(c)\mathrm{ver}(e^{0}).
	 	\end{aligned}
	 \end{equation} \qed    
\end{prop}
We define a second order differential calculus $(\Omega^{2}(A),\wedge,\dd)$ over $A$ as the free left $A-$module generated as $\Omega^{2}(A)=\span_{A}\{ e^{\pm}\wedge e^{0},e^{+}\wedge e^{-}\}$ with commutation relations \cite{beggs-majid} \eqa{e^{+}\wedge e^{-}& =-q^{-2}e^{-}\wedge e^{+}, \quad e^{\pm}\wedge e^{0}=-q^{\mp 4}e^{0}\wedge e^{\pm}, \quad \dd e^{0}=q^{3}e^{+}\wedge e^{-}, \\ &  \dd e^{\pm}=\mp q^{\pm 2}[2]_{q^{-2}}e^{\pm}\wedge e^{0}, \quad e^{\pm}\wedge e^{\pm}=e^{0}\wedge e^{0}=0,}  where $[2]_{q}=(1-q^{2})/(1-q)$. 
\begin{prop}
	 $\Omega^{2}(A)$ is a right $H-$covariant differential calculus. Moreover, maps \eqa{  \ver^{1,1}& : \Omega^{2}(A)\rightarrow \Omega^{1}(A)\otimes \Omega^{1}(H), \\ \ver^{0,2}& :\Omega^{2}(A)\rightarrow A\otimes \Omega^{2}(H)} are well defined.  \proof 
	The morphism $\ver^{0,2}$ is trivial since no non-zero 2-forms occur on $H$. On the generators we define \eqa{\Delta_{\Omega^{2}(A)}(e^{+}\wedge e^{-})& := e^{+}\wedge e^{-}\otimes 1, \quad \Delta_{\Omega^{2}(A)}(e^{+}\wedge e^{0}) := e^{+}\wedge e^{0}\otimes t^{2}, \\ & \Delta_{\Omega^{2}(A)}(e^{-}\wedge e^{0}):= e^{-}\wedge e^{0}\otimes t^{-2}.}  These maps are well-defined because compatible with all relations of 2-forms. In fact \eqa{\nonumber\Delta_{\Omega^{2}(A)}(e^{+}\wedge e^{-}) & = \Delta_{\Omega^{1}(A)}(e^{+})\Delta_{\Omega^{1}(A)}(e^{-}) \\ & = (e^{+}\otimes t^{2})(e^{-}\otimes t^{-2})\\ & = e^{+}\wedge e^{-} \otimes 1; \\ \Delta_{\Omega^{2}(A)}(e^{\pm}\wedge e^{0}) & = \Delta_{\Omega^{1}(A)}(e^{\pm}) \Delta_{\Omega^{1}(A)}(e^{0}) \\ & = (e^{\pm}\otimes t^{\pm 2})( e^{0}\otimes 1) \\ & = e^{\pm}\wedge e^{0}\otimes t^{\pm 2},\\ \Delta_{\Omega^{2}(A)}(-q ^{2} e^{-}\wedge e^{+}) & = -q^{2}\Delta_{\Omega^{2}(A)}(e^{-}\wedge e^{+}) \\ & = -q^{2}\Delta_{\Omega^{1}(A)}(e^{-})\Delta_{\Omega^{1}(A)}(e^{+}) \\ & = -q^{2} (e^{-}\otimes t^{-2})(e^{+}\otimes t^{2}) \\ & = -q^{2}e^{-}\wedge e^{+} \otimes 1 \\ & =  e^{+}\otimes e^{-}\otimes 1 \\ & = \Delta_{\Omega^{2}(A)}(e^{+}\wedge e^{-});  }  
	\eqa{\nonumber \Delta_{\Omega^{2}(A)}(-q ^{4} e^{0}\wedge e^{\pm}) & = -q^{4} \Delta_{\Omega^{1}(A)}(e^{0})\Delta_{\Omega^{1}(A)}(e^{\pm}) \\ & =  -q^{4}(e^{0}\otimes 1 )(e^{\pm}\otimes t^{\pm 2}) \\ & = -q^{4}e^{0}\wedge e^{\pm}\otimes t^{\pm 2} \\ & = e^{\pm}\wedge e^{0}\otimes t^{\pm} \\ & = \Delta_{\Omega^{2}(A)}(e^{\pm}\wedge e^{0}).\\ \nonumber \Delta_{\Omega^{2}(A)}(\dd e^{0}) & = (\dd \otimes \id)\Delta_{\Omega^{1}(A)}(e^{0}) \\ & = \dd e^{0} \otimes 1 \\ & = e^{+}\wedge e^{-}\otimes 1; \\ 
\Delta_{\Omega^{2}(A)}(\dd e^{+}) & = (\dd \otimes \id) \Delta_{\Omega^{1}(A)}(e^{+}) \\ & =  \dd e^{+} \otimes t^{2} \\ & = -q^{2}[2]_{q^{-2}} e^{+}\wedge e^{0}\otimes t^{2} \\ & = -q^{2}[2]_{q^{-2}}(e^{+}\otimes t^{2})(e^{0}\otimes 1) \\ & = -q^{2}[2]_{q^{-2}} \Delta_{\Omega^{1}(A)}(e^{+})\Delta_{\Omega^{1}(A)}(e^{0}) \\ & = \Delta_{\Omega^{2}(A)}(-q^{2}[2]_{q^{-2}}e^{+}\wedge e^{0});  \\ \Delta_{\Omega^{2}(A)}(\dd e^{-}) & = (\dd \otimes \id) \Delta_{\Omega^{1}(A)}(e^{-}) \\ & =  \dd e^{-} \otimes t^{-2} \\ & = q^{-2}[2]_{q^{-2}} e^{-}\wedge e^{0}\otimes t^{-2} \\ & = q^{-2}[2]_{q^{-2}}(e^{-}\otimes t^{-2})(e^{0}\otimes 1) \\ & = q^{-2}[2]_{q^{-2}} \Delta_{\Omega^{1}(A)}(e^{-})\Delta_{\Omega^{1}(A)}(e^{0}) \\ & = \Delta_{\Omega^{2}(A)}(q^{-2}[2]_{q^{-2}}e^{-}\wedge e^{0}). }
	
Moving to the vertical map $\ver^{1,1}$ we have \eqa{ \nonumber \ver^{1,1}(e^{+}\wedge e^{-}) & = \Delta_{\Omega^{1}(A)}(e^{+})\ver(e^{-}) + \ver(e^{+})\Delta_{\Omega^{1}(A)}(e^{-}) \\ & = 0 \\ & = \ver^{1,1}(\dd e^{0}),\\  \ver^{1,1}(e^{\pm}\wedge e^{0}) & = \Delta_{\Omega^{1}(A)}(e^{\pm})\ver(e^{0}) + \ver(e^{\pm})\Delta_{\Omega^{1}(A)}(e^{0}) \\ & = \Delta_{\Omega^{1}(A)}(e^{\pm})\ver(e^{0}) \\ & = (e^{\pm}\otimes t^{\pm 2}) (1\otimes t^{-1}\dd t) \\ & =  e^{\pm} \otimes t^{\pm 2 -1}\dd t,} and so we define \eqa{ \ver^{1,1}(e^{+}\wedge e^{-}):=0, \quad \ver^{1,1}(e^{\pm}\wedge e^{0}):= e^{\pm}\otimes t^{\pm 2 -1}\dd t.}  
	With these definitions on the generators the map $\ver^{1,1}$ extends correctly to $\Omega^{2}(A)$ provided the commutation relations among generators of $\Omega^{2}(A)$ are preserved. 
	
The relation $\ver^{1,1}(-q^{2}e^{-}\wedge e^{+})=\ver^{1,1}(e^{+}\wedge e^{-})=0$ is obvious. Moreover, in \cite{beggs-majid} (page 109) the following relations are provided: \eqa{ \dd \alpha= \alpha e^{0}+q\beta e^{+},\quad \dd \beta=\alpha e^{-}-q^{-2}\beta e^{0}, \quad \dd \gamma=\gamma e^{0}+q\delta e^{+},\quad \dd \delta =\gamma e^{-}-q^{-2}\delta e^{0}.}  
 Therefore \eqa{\nonumber\ver^{1,1}\circ \dd (e^{+}) & = \ver^{1,1}\circ \dd (q^{-1}\alpha\wedge \dd \gamma -q^{-2}\gamma \dd \alpha) \\ & = \ver^{1,1}(q^{-1}\dd \alpha \wedge \dd \gamma - q^{-2} \dd \gamma \wedge \dd \alpha) \\ & = q^{-1}\ver^{1,1}(\dd\alpha\wedge \dd \gamma) - q^{-2}\ver^{1,1}(\dd \gamma \wedge \dd \alpha)\\ & = q^{-1}\ver^{1,1}[(\alpha e^{0}+q\beta e^{+})\wedge (\gamma e^{0}+q \delta e^{+})] - q^{-2}\ver^{1,1}[ (\gamma e^{0}+q \delta e^{+})\wedge (\alpha e^{0}+q\beta e^{+})] \\ & = \ver^{1,1}(\alpha e^{0}\wedge \delta e^{+}) +  \ver^{1,1}(\beta e^{+}\wedge \gamma e^{0}) - q^{-1}\ver^{1,1}(\gamma e^{0}\wedge \beta e^{+}) - q^{-1}\ver^{1,1}(\delta e^{+}\wedge\alpha e^{0})}\eqa{\nonumber  \qquad \quad \quad \ \  & = \Delta_{A}(\alpha)\ver(e^{0})\Delta_{A}(\delta)\Delta_{\Omega^{1}(A)}(e^{+}) + \Delta_{A}(\beta)\Delta_{\Omega^{1}(A)}(e^{+})\Delta_{A}(\gamma)\ver(e^{0}) \\ & - q^{-1}\Delta_{A}(\gamma)\ver(e^{0})\Delta_{A}(\beta)\Delta_{\Omega^{1}(A)}(e^{+}) - q^{-1}\Delta_{A}(\delta)\Delta_{\Omega^{1}(A)}(e^{+})\Delta_{A}(\alpha)\ver(e^{0}) \\ & = (\alpha\otimes t)(1\otimes t^{-1}\dd t)(\delta\otimes t^{-1})(e^{+}\otimes t^{2}) + (\beta\otimes t^{-1})(e^{+}\otimes t^{2})(\gamma\otimes t)(1\otimes t^{-1}\dd t) \\ & - q^{-1}(\gamma\otimes t)(1\otimes t^{-1}\dd t)(\beta\otimes t^{-1})(e^{+}\otimes t^{2})- q^{-1}(\delta\otimes t^{-1})(e^{+}\otimes t^{2})(\alpha\otimes t)(1\otimes t^{-1}\dd t) \\ & = -\alpha\delta e^{+}\otimes \dd t t + \beta e^{+}\gamma \otimes t \dd t + q^{-1}\gamma \beta e^{+}\otimes \dd t t - q^{-1}\delta e^{+}\alpha\otimes t \dd t \\ & = (-q^{2}\alpha\delta e^{+} +q\beta \gamma e^{+} +q \gamma \beta e^{+} - \delta \alpha e^{+})\otimes t \dd t \\ & = -(q^{2}(\alpha\delta -q^{-1}\gamma\beta) + (\delta\alpha -q\gamma \beta))e^{+}\otimes t \dd t \\ & = - (q^{2}+1)e^{+}\otimes t \dd t,  } which equals \eqa{ \nonumber\ver^{1,1}(-q^{2}[2]_{q^{-2}}e^{+}\wedge e^{0}) & = -(1+q^{2})e^{+}\otimes t \dd t. } Moreover   
\eqa{\nonumber \ver^{1,1}\circ \dd (e^{0}) & = \ver^{1,1}\circ \dd (\delta \dd \alpha-q\beta \dd \gamma) \\ & = \ver^{1,1}(\dd\delta \wedge\dd \alpha - q \dd \beta\wedge \dd \gamma) \\ & = \ver^{1,1}[(\gamma e^{-}-q^{-2}\delta e^{0})\wedge (\alpha e^{0}+q\beta e^{+}) - q (\alpha e^{-}-q^{-2}\beta e^{0})\wedge(\gamma e^{0}+q\delta e^{+})]\\ & = \ver^{1,1}(\gamma e^{-}\wedge \alpha e^{0}-q^{-1}\delta e^{0}\wedge\beta e^{+}-q\alpha e^{-}\wedge \gamma e^{0} +q^{-1}\beta e^{0}\wedge \delta e^{+}) \\ & = \Delta_{A}(\gamma)\Delta_{\Omega^{1}(A)}(e^{-})\Delta_{A}(\alpha)\ver(e^{0})- q^{-1}\Delta_{A}(\delta)\ver(e^{0})\Delta_{A}(\beta)\Delta_{\Omega^{1}(A)}(e^{+}) \\ & - q \Delta_{A}(\alpha)\Delta_{\Omega^{1}(A)}(e^{-})\Delta_{A}(\gamma)\ver(e^{0})+ q^{-1}\Delta_{A}(\beta)\ver(e^{0})\Delta_{A}(\delta)\Delta_{\Omega^{1}(A)}(e^{+}) \\ & = (\gamma\otimes t)(e^{-}\otimes t^{-2})(\alpha\otimes t)(1\otimes t^{-1}\dd t) - q^{-1}(\delta\otimes t^{-1})(1\otimes t^{-1}\dd t)(\beta\otimes t^{-1})(e^{+}\otimes t^{2}) \\ & -q(\alpha\otimes t)(e^{-}\otimes t^{-2})(\gamma\otimes t)(1\otimes t^{-1}\dd t)+ q^{-1}(\beta \otimes t^{-1})(1\otimes t^{-1}\dd t)(\delta \otimes t^{-1})(e^{+}\otimes t^{2})\\ & =  \gamma e^{-}\alpha \otimes t^{-1}\dd t +q\delta \beta e^{+} \otimes t^{-1}\dd t  -q\alpha e^{-}\gamma \otimes t^{-1}\dd t + q\beta \delta e^{+}\otimes t^{-1}\dd t \\ & = (q\gamma \alpha +q\delta \beta -q^2\alpha \gamma +q^2\beta\delta ) e^{+}\otimes t^{-1}\dd t \\ & = 0.  } 

\eqa{\nonumber\ver^{1,1}\circ \dd (e^{-}) & = \ver^{1,1}\circ \dd (\delta \dd \beta - q \beta \dd \delta) \\ & = \ver^{1,1}(\dd \delta \wedge \dd \beta -q\dd \beta \wedge \dd \delta) \\ & = \ver^{1,1}[(\gamma e^{-}-q^{-2}\delta e^{0})\wedge (\alpha e^{-}-q^{-2}\beta e^{0})-q(\alpha e^{-}-q^{-2}\beta e^{0})\wedge (\gamma e^{-}-q^{-2}\delta e^{0})] \\ & = \ver^{1,1}(-q^{-2}\gamma e^{-}\wedge \beta e^{0} -q^{-2}\delta e^{0}\wedge \alpha e^{-} +q^{-1}\alpha e^{-}\wedge \delta e^{0} +q^{-1}\beta e^{0}\wedge \gamma e^{-})\\ & = -q^{-2}\Delta_{A}(\gamma)\Delta_{\Omega^{1}(A)}(e^{-})\Delta_{A}(\beta)\ver(e^{0}) - q^{-2} \Delta_{A}(\delta)\ver(e^{0}) \Delta_{A}(\alpha)\Delta_{\Omega^{1}(A)}(e^{-}) \\ & +q^{-1}\Delta_{A}(\alpha)\Delta_{\Omega^{1}(A)}(e^{-})\Delta_{A}(\delta)\ver(e^{0}) + q^{-1}\Delta_{A}(\beta) \ver(e^{0}) \Delta_{A}(\gamma)\Delta_{\Omega^{1}(A)}(e^{-}) \\ & = -q^{-2}(\gamma\otimes t)(e^{-}\otimes t^{-2})(\beta \otimes t^{-1})(1\otimes t^{-1}\dd t) - q^{-2}(\delta\otimes t^{-1})(1\otimes t^{-1}\dd t)(\alpha\otimes t)(e^{-}\otimes t^{-2})}\eqa{\nonumber \qquad \quad \ \ \ \ \  & + q^{-1}(\alpha\otimes t)(e^{-}\otimes t^{-2})(\delta \otimes t^{-1})(1\otimes t^{-1}\dd t) + q^{-1}(\beta\otimes t^{-1})(1\otimes t^{-1}\dd t)(\gamma \otimes t)(e^{-}\otimes t^{-2}) \\  & = -q^{-2}\gamma e^{-}\beta \otimes t^{-3}\dd t + q^{-2} \delta \alpha e^{-}\otimes t^{-2}\dd t t^{-1} + q^{-1}\alpha e^{-}\delta \otimes t^{-3}\dd t + q^{-1} \beta \gamma e^{-}t^{-2}\dd t t^{-1}    \\ & =  -q^{-3}\gamma \beta e^{-} \otimes t^{-3}\dd t + q^{-4}\delta \alpha e^{-}\otimes t^{-3}\dd t   +q^{-2}\alpha \delta e^{-} \otimes t^{-3}\dd t - q^{-3}\beta \gamma e^{-}t^{-3}\dd t  \\ & = q^{-2}(-q^{-1}\gamma \beta + q^{-2}\delta \alpha + \alpha\delta  - q^{-1}\beta \gamma)e^{-}\otimes t^{-3}\dd t \\ & = q^{-2}(-q^{-1}\gamma \beta + q^{-2}(\alpha\delta + (q-q^{-1})\beta\gamma)+\alpha\delta -q^{-1}\beta\gamma) e^{-}\otimes t^{-3}\dd t \\ & = q^{-2}(1+q^{-2})e^{-}\otimes t^{-3}\dd t, } which equals \eqa{\nonumber \ver^{1,1}(q^{-2}[2]_{q^{-2}}e^{-}\wedge e^{0}) & = q^{-2}(1+q^{-2}) \ver^{1,1}(e^{-}\wedge e^{0}) \\ & = q^{-2}(1+q^{-2}) e^{-}\otimes t^{-3}\dd t. }

Finally \eqa{ \nonumber \ver^{1,1}(-q^{\mp 4}e^{0}\wedge e^{\pm}) & = -q^{\mp 4}(\Delta_{\Omega^{1}(A)}(e^{0})\ver(e^{\pm}) + \ver(e^{0}) \Delta_{\Omega^{1}(A)}(e^{\pm})) \\ & = -q^{\mp 4}\ver(e^{0}) \Delta_{\Omega^{1}(A)}(e^{\pm}) \\ & = - q^{\mp 4}(1\otimes t^{-1}\dd t )( e^{\pm}\otimes t^{\pm 2}) \\ & = -q^{\mp 4}e^{\pm} \otimes t^{-1}\dd t \ t^{\pm 2} \\ & =-q^{\mp 4} e^{\pm}\otimes t^{-1}\dd t \ t^{\pm 2}\\ & =e ^{\pm} \otimes t^{\pm 2-1}\dd t\\ & = \ver^{1,1}(e^{\pm}\wedge e^{0}).  }  

In the various calculations we exploited $$\dd t t^{2} = q^{4} t^{2} \dd t, \quad  \dd t t^{-2} = q^{-4}t^{-2}\dd t $$ beside the various commutation relations between the generators of the algebra and Proposition \eqref{decomposition of vertical maps}.  We extend to $\Omega^{2}(A)$ as \begin{equation}
	 	\begin{aligned}\nonumber
	 		\Delta_{\Omega^{2}(A)}(ae^{+}\wedge e^{-}  & + b e^{+}\wedge e^{0} + ce^{-}\wedge e^{0}) \\  & = \Delta_{A}(a)\Delta_{\Omega^{2}(A)}(e^{+}\wedge e^{-}) + \Delta_{A}(b)\Delta_{\Omega^{2}(A)} (e^{+}\wedge e^{0}) + \Delta_{A}(c)\Delta_{\Omega^{2}(A)}(e^{-}\wedge e^{0}), \\ \mathrm{ver}^{1,1}(ae^{+}\wedge e^{-}  & + b e^{+}\wedge e^{0} + ce^{-}\wedge e^{0}) \\ &  = \Delta_{A}(a)\mathrm{ver}^{1,1}(e^{+}\wedge e^{-}) + \Delta_{A}(b)\mathrm{ver}^{1,1}(e^{+}\wedge e^{0})  \Delta_{A}(c)\mathrm{ver}^{1,1}(e^{-}\wedge e^{0}), \\ \mathrm{ver}^{0,2}(ae^{+}\wedge e^{-} & + b e^{+}\wedge e^{0} + ce^{-}\wedge e^{0})  = 0. 
	 	\end{aligned}
	 \end{equation}
 \qed 
\end{prop}
We define a third order differential calculus $(\Omega^{3}(A),\wedge,\dd)$ over $A$ as the free left $A-$module generated as $\Omega^{3}(A)=\span_{A}\{e^{+}\wedge e^{-}\wedge e^{0}\},$ with commutation relations \eqa{\dd (e^{+}\wedge e^{-})=-q^{-2}\dd (e^{-}\wedge e^{+}), \quad \dd (e^{\pm}\wedge e^{0}) =-q^{\mp 4}\dd(e^{0}\wedge e^{\pm}).}     

\begin{prop}
$\Omega^{3}(A)$ is right $H-$covariant. Moreover the vertical maps \eqa{ \ver^{2,1}& : \Omega^{3}(A)\rightarrow \Omega^{2}(A)\otimes \Omega^{1}(H), \\ \ver^{1,2}& :\Omega^{3}(A)\rightarrow \Omega^{1}(A)\otimes \Omega^{2}(H), \\ \ver^{0,3}& :\Omega^{3}(A)\rightarrow A\otimes \Omega^{3}(H),} are well defined. 
  \proof  We start noticing $\ver^{0,3}$ and $\ver^{1,2}$ are trivial since $\Omega^{k}(H)=0$ for $k>1$.  
  We need to show the coaction $\Delta_{A}:A\rightarrow A\otimes H$ extends to $\Omega^{3}(A)$  in the correct way. On the generator we have: 
  \eqa{\nonumber\Delta_{\Omega^{3}(A)}(e^{+}\wedge e^{-}\wedge e^{0}) & = \Delta_{\Omega^{1}(A)}(e^{+})\Delta_{\Omega^{1}(A)}(e^{-})\Delta_{\Omega^{1}(A)}(e^{0}) \\ & = (e^{+}\otimes t^{2})(e^{-}\otimes t^{-2})(e^{0}\otimes 1) \\ & = e^{+}\wedge e^{-}\wedge e^{0}\otimes 1; \\ \ver^{2,1}((e^{+}\wedge e^{-})\wedge e^{0}) & = \Delta_{\Omega^{2}(A)}(e^{+}\wedge e^{-})\ver(e^{0}) + \ver^{1,1}(e^{+}\wedge e^{-}) \Delta_{\Omega^{1}(A)}(e^{0}) \\ & = (e^{+}\wedge e^{-}\otimes 1)[1 \otimes t^{-1}\dd t] \\ & =  (e^{+}\wedge e^{-}) \otimes t^{-1}\dd t \\ & = \ver^{2,1}(e^{+}\wedge (e^{-}\wedge e^{0})).  }  Accordingly we define \eqa{\Delta_{\Omega^{3}(A)}(e^{+}\wedge e^{-}\wedge e^{0}) & :=( e^{+}\wedge e^{-}\wedge e^{0})\otimes 1, \\  \ver^{2,1}(e^{+}\wedge e^{-}\wedge e^{0})& := (e^{+}\wedge e^{-}\wedge e^{0}) \otimes t^{-1}\dd t.} With these definitions the map $\Delta_{A}^{3}:=\Delta_{\Omega^{3}(A)}+\ver^{2,1}$ supplies the extension of $\Delta_{A}$ to the differential calculus $\Omega^{3}(A)$ provided the commutation relations between elements of the calculus are preserved. We have
  
  \eqa{ \nonumber\Delta_{\Omega^{3}(A)}\circ \dd (e^{\pm}\wedge e^{0}) & = \Delta_{\Omega^{3}(A)}(\dd e^{\pm}\wedge e^{0}- e^{\pm}\wedge \dd e^{0}) \\ & = \Delta_{\Omega^{2}(A)}(\dd e^{\pm})\Delta_{\Omega^{1}(A)}(e^{0}) - \Delta_{\Omega^{1}(A)}(e^{\pm})\Delta_{\Omega^{2}(A)}(\dd e^{0}) \\ & = (\dd e^{\pm}\otimes t^{\pm 2})(e^{0}\otimes 1)- (e^{\pm}\otimes t^{\pm})(\dd e^{0}\otimes 1) \\ & =(\dd e^{\pm}\wedge e^{0}-e^{\pm}\wedge \dd e^{0})\otimes t^{\pm} \\ & = \dd (e^{\pm}\wedge e^{0})\otimes t^{2} \\ & = -q^{\mp 4}\dd (e^{0}\wedge e^{\pm}) \otimes t^{\pm} \\ & = \Delta_{\Omega^{3}(A)}\tonde{-q^{\mp 4}\dd(e^{0}\wedge e^{\pm})},}  
  and \eqa{\nonumber \Delta_{\Omega^{3}(A)}\circ \dd (e^{+}\wedge e^{-}) & = \Delta_{\Omega^{3}(A)}(\dd e^{+}\wedge e^{-}-e^{+}\wedge \dd e^{-}) \\ & =\Delta_{\Omega^{2}(A)}(\dd e^{+})\Delta_{\Omega^{1}(A)}(e^{-}) - \Delta_{\Omega^{1}(A)}(e^{+})\Delta_{\Omega^{2}(A)}(\dd e^{-}) \\ & = (\dd e^{+}\otimes t^{2})(e^{-}\otimes t^{-2}) - (e^{+}\otimes t^{2})(\dd e^{-}\otimes t^{-2}) \\ & = (\dd e^{+}\wedge e^{-}- e^{+}\wedge \dd e^{-})\otimes 1 \\ & = \dd (e^{+}\wedge e^{-})\otimes 1 \\ & =  -q^{-2}\dd (e^{-}\wedge e^{+})\otimes 1 \\ & = \Delta_{\Omega^{3}(A)}(-q^{-2}\dd(e^{-}\wedge e^{+})).  }  Moving to the vertical maps we find

  \eqa{\nonumber\ver^{2,1}\circ \dd(e^{+}\wedge e^{-}) & = \ver^{2,1}(\dd e^{+}\wedge e^{-} - e^{+}\wedge \dd e^{-}) \\ & = \Delta_{\Omega^{2}(A)}(\dd e^{+}) \ver(e^{-})+ \ver^{1,1}(\dd e^{+})\Delta_{\Omega^{1}(A)}(e^{-})\\ &  - \Delta_{\Omega^{1}(A)}(e^{+})\ver^{1,1}(\dd e^{-})- \ver(e^{+})\Delta_{\Omega^{2}(A)}(\dd e^{-}) \\ & =-[(q^{2}+1)e^{+}\otimes t \dd t](e^{-}\otimes t^{-2})- (e^{+}\otimes t^{2})(q^{-2}(1+q^{-2})e^{-}\otimes t^{-3}\dd t)\\ & = q^{-4}(q^{2}+1)e^{+}\wedge e^{-}\otimes t^{-1}\dd t -q^{-2}(1+q^{-2})e^{+}\wedge e^{-} \otimes t^{-1}\dd t \\ & = -(q^{-2}+q^{-4}-q^{-2}-q^{-4})e^{+}\wedge e^{-}\otimes t^{-1}\dd t \\ & = 0, } which equals. 
  \eqa{\nonumber -q^{-2}\ver^{2,1}\circ \dd (e^{-}\wedge e^{+}) & = -q^{-2}\ver^{2,1}(\dd e^{-}\wedge e^{+} - e^{-}\wedge \dd e^{+}) \\ & -q^{-2}\ver^{1,1}(\dd e^{-})\Delta_{\Omega^{1})(A)}(e^{+})- q^{-2}\Delta_{\Omega^{1}(A)}(e^{-}) \ver^{1,1}(\dd e^{+}) \\ & = -q^{-2}(q^{-2}(1+q^{-2}) e^{-}\otimes t^{-3}\dd t)(e^{+}\otimes t^{2}) - q^{-2}(e^{-}\otimes t^{-2})((1+q^{2})e^{+}\otimes t\dd t) \\ & = q^{-4}(1+q^{-2})e^{-}\wedge e^{+}\otimes t^{-3}\dd t t^{2} - q^{-2}(1+q^{2})e^{-}\wedge e^{+} \otimes t^{-1}\dd t \\ & = (1+q^{-2})e^{-}\wedge e^{+}\otimes t^{-1}\dd t - q^{-2}(1+q^{2}) e^{-}\wedge e^{+}\otimes t^{-1}\dd t  \\ & = [1+q^{-2}-q^{-2}-1]e^{-}\wedge e^{+}\otimes t^{-1}\dd t \\ & = -q^{2}[1+q^{-2}-q^{-2}-1]e^{+}\wedge e^{-}\otimes t^{-1}\dd t \\ & = 0.   } Moreover
  \eqa{\nonumber-\ver^{2,1}(e^{0}\wedge e^{+}) & = - \ver^{2,1}(\dd e^{0}\wedge e^{+}-e^{0}\wedge \dd e^{+}) \\ & = - (\ver^{1,1}(\dd e^{0}) \Delta_{\Omega^{1}(A)}(e^{+})- \Delta_{\Omega^{1}(A)}(e^{0})\ver^{1,1}(\dd e^{+})- \ver(e^{0})\Delta_{\Omega^{2}(A)}(\dd e^{+})) \\ & =\Delta_{\Omega^{1}(A)}(e^{0})\ver^{1,1}(\dd e^{+})+  \ver(e^{0})\Delta_{\Omega^{2}(A)}(\dd e^{+}) \\ & = (e^{0}\otimes 1)(-(q^{2}+1)e^{+}\otimes t^{-1}\dd t)+(1\otimes t^{-1}\dd t)(\dd e^{+} \otimes t^{2}) \\ & =  -(q^{2}+1)e^{0}\wedge e^{+}\otimes t^{-1}\dd t - q^{4}\dd e^{+}\otimes t \dd t,   } which equals
  \eqa{\nonumber q^{4}\ver^{2,1}\circ \dd(e^{+}\wedge e^{0})&  = q^{4}\ver^{2,1}(\dd e^{+}\wedge e^{0}- e^{+}\wedge \dd e^{0}) \\ & = q^{4}(\ver^{1,1}(\dd e^{+})\Delta_{\Omega^{1}(A)}(e^{0})- \Delta_{\Omega^{1}(A)}(\dd e^{+})\ver(e^{0})) \\ & = q^{4}[(-(q^{2}+1)e^{+}\otimes t^{-1}\dd t )(e^{0}\otimes 1) - (\dd e^{+}\otimes t^{2})(1\otimes t^{-1}\dd t)] \\ & = q^{4}(q^{2}+1)e^{+}\wedge e^{0} \otimes t^{-1}\dd t - q^{4}\dd e^{+}\otimes t \dd t \\ & = -(q^{2}+1)e^{0}\wedge e^{+} \otimes t^{-1}\dd t - q^{4}\dd e^{+}\otimes t \dd t. }  
  Finally
    \eqa{\nonumber-\ver^{2,1}(e^{0}\wedge e^{-}) & = - \ver^{2,1}(\dd e^{0}\wedge e^{-}-e^{0}\wedge \dd e^{-}) \\ & = - (\ver^{1,1}(\dd e^{0}) \Delta_{\Omega^{1}(A)}(e^{-})- \Delta_{\Omega^{1}(A)}(e^{0})\ver^{1,1}(\dd e^{-})- \ver(e^{0})\Delta_{\Omega^{2}(A)}(\dd e^{-})) \\ & =\Delta_{\Omega^{1}(A)}(e^{0})\ver^{1,1}(\dd e^{-})+  \ver(e^{0})\Delta_{\Omega^{2}(A)}(\dd e^{-}) \\ & = (e^{0}\otimes 1)(q^{-2}(1+q^{-2})e^{-}\otimes t^{-3}\dd t)+(1\otimes t^{-1}\dd t)(\dd e^{-} \otimes t^{-2}) \\ & =  q^{-2}(1+q^{-2})e^{0}\wedge e^{-}\otimes t^{-3}\dd t - q^{-4}\dd e^{-}\otimes t^{-3}\dd t\\ & = (q^{-2}(1+q^{-2})e^{0}\wedge e^{-}  - q^{-4}\dd e^{-} )\otimes t^{-3}\dd t,   } 
    that is
    \eqa{\nonumber q^{-4}\ver^{2,1}\circ \dd(e^{-}\wedge e^{0})&  = q^{-4}\ver^{2,1}(\dd e^{-}\wedge e^{0}- e^{-}\wedge \dd e^{0}) \\ & = q^{-4}(\ver^{1,1}(\dd e^{-})\Delta_{\Omega^{1}(A)}(e^{0})- \Delta_{\Omega^{1}(A)}(\dd e^{-})\ver(e^{0})) \\ & = q^{-4}[(q^{-2}(q^{-2}+1)e^{-}\otimes t^{-3}\dd t )(e^{0}\otimes 1) - (\dd e^{-}\otimes t^{-2})(1\otimes t^{-1}\dd t)] \\ & = -q^{-4}q^{-2}(q^{-2}+1)e^{-}\wedge e^{0} \otimes t^{-3}\dd t - q^{-4}\dd e^{-}\otimes t^{-3} \dd t \\ & = (-q^{-4}q^{-2}(q^{-2}+1)e^{-}\wedge e^{0} -q^{-4}\dd e^{-})\otimes t^{-3}\dd t \\ & = [q^{-2}(1+q^{-2})e^{0}\wedge e^{-}-q^{-4}\dd e^{-}]\otimes t^{-3}\dd t. }  
  We extend to $\Omega^{3}(A)$ as \begin{equation}
	 	\begin{aligned} \Delta_{\Omega^{3}(A)}(ae^{+}\wedge e^{-}\wedge e^{0})&=\Delta_{A}(a)\Delta_{\Omega^{3}(A)}(e^{+}\wedge e^{-}\wedge e^{0}), \\ \mathrm{ver}^{2,1}(ae^{+}\wedge e^{-}\wedge e^{0})&=\Delta_{A}(a)\mathrm{ver}^{2,1}(e^{+}\wedge e^{-}\wedge e^{0}),
	 	\end{aligned}
	 \end{equation}
the other vertical maps vanishing. 
  
   \qed  
\end{prop}
\begin{thm}
$\Omega^{\bullet}(A)$ is a complete differential calculus. \proof 
According to the last propositions we have the existence of vertical maps and right $H-$covariance up to the third order calculus. Moreover $\Omega^k(A)=0$ for $k>3$ since $e^+\wedge e^+=e^{-}\wedge e^{-}=e^{0}\wedge e^{0}=0$.  
  We conclude $\Omega^{\bullet}(A)= A \oplus \Omega^{1}(A) \oplus \Omega^{2}(A) \oplus \Omega^{3}(A) $ is a complete differential calculus over $A$.  \qed 
	
\end{thm}
\begin{rmk}In \cite{beggs-majid} (Proposition 2.35, page 111) it is stated that $\Omega^{1}(B)$ is spanned by \eqa{\dd x & = -1\beta \delta e^{+} - \alpha \gamma e^{-}, \\ \dd z & = \dd (\gamma \delta)= \delta^{2}e^{+}+ \gamma^{2}e^{-}, \\ \dd \bar{z}& =-q\dd (\alpha\beta)=-q\beta^{2}e^{+}-q \alpha^{2}e^{-}.}  Still in \cite{beggs-majid} (page 112) the following relations are provided. \eqa{ \delta e^{+} & = \alpha \dd z + q^{-1} \gamma \dd x, \\ \beta e^{+} & =  q^{-2} \gamma \dd z - q \alpha \dd x, \\ \alpha e^{-} & = q^{2}\beta \dd x - q^{-1} \delta \dd \bar{z}, \\ \gamma e^{-} & = - \delta \dd x - q\beta \dd z.} These are of particular interest in the proof of the next proposition. \qed 
\end{rmk}

\begin{prop} $\Omega^{\bullet}(B)$ is a differential calculus.  
	\proof  let us consider a generic element $\omega\in \Omega^{1}(B)$. We may write \eqa{ \omega= ae^{+}+be^{-}+ce^{0}, \quad a,b,c\in A,} and since we must have $\Delta_{A}^{1}(\omega)=\omega \otimes 1$ we find \eqa{ \Delta_{A}^{1}(ae^{+}+be^{-}+ce^{0}) & = \Delta_{A}(a)\Delta_{A}^{\wedge}(e^{+}) + \Delta_{A}(b)\Delta_{A}^{\wedge}(e^{-})+\Delta_{A}(c)\Delta_{A}^{\wedge}(e^{0}) \\ & = (a\otimes t^{|a|})(e^{+}\otimes t^{2}) + (b\otimes t^{|b|})(e^{-}\otimes t^{-2}) \\ & + (c\otimes t^{|c|})[e^{0}\otimes 1 +1\otimes t^{-1}\dd t] \\ & = ae^{+}\otimes t^{|a|}t^{2} + be^{-}\otimes t^{-2}t^{|b|} + ce^{+}\otimes t^{|c|} + c\otimes t^{|c|}t^{-1}\dd t  \\ & = (a e^{+} + b e^{-} + ce^{0})\otimes 1 } if and only if we have one the following combinations \eqa{
	 |a|=-2,|b|=2,c=0,\quad  a=c=0,|b|=2,\quad  b=c=0,|a|=-2, \quad a=b=c=0. 
}
	 Accordingly every element in $\Omega^{1}(B)$ must be of the form $\omega=a e^{+} + b e^{-}$ with the stated prescriptions. 
	We write down the possible elements of degree $2,-2$ in terms of generators: $$ \text{degree 2}: \alpha^{2},\gamma^{2},\alpha\gamma; \quad \text{degree -2}:\delta^{2},\beta^{2},\beta\gamma.$$  
	
	We now show that every possible element of the form $a e^{+}+b e^{-}$ with $a$ of degree $-2$ and $b$ of degree $2$ gives rise to an element of the form $B\dd B$. For elements of the form $be^{-}$ we find \eqa{ \alpha^{2} e^{-} & = \alpha (q^{2}\beta \dd x - q^{-1} \delta \dd z) \\ & = q^{2}\alpha \beta \dd x - q^{-2} \alpha\delta \dd z; \\  \gamma^{2}e^{-} & = \gamma (-\delta \dd x - q \beta \dd z) \\ & = - \gamma\delta \dd x- q \gamma \beta \dd z;  \\ 
\alpha\gamma e^{-} & = \alpha(-\delta \dd x- q \beta \dd z) \\ & = -\alpha \delta \dd x - q \alpha \beta \dd z.} Each of the above expressions define elements of the form $B\dd B$ as we are always pairing degree $1$ and degree $-1$ elements. By the same reasoning we exploit elements of the form $ae^{+}$, for which \eqa{ \beta^{2}e^{+} & = \beta (q^{-2} \gamma \dd z - q \alpha \dd x) \\ & = q^{-2}\beta \gamma \dd z- q \beta \alpha \dd x; \\ \delta^{2}e^{+} & = \delta (\alpha \dd z + q^{-1}\gamma \dd x) \\ & = \delta \alpha \dd z + q^{-1}\delta \gamma \dd x; \\ \delta \beta e^{+} & = \delta (q^{-2}\gamma \dd z -q\alpha \dd x) \\ & = q^{-2} \delta \gamma \dd z - q \delta \alpha \dd x. } Moreover, in \cite{beggs-majid} (Proposition 2.35, page 113) it is stated that the volume form can be expressed in terms of elements in $B$ and $\Omega^{1}(B)$. Therefore $\Omega^{\bullet}(B)$ is a differential calculus.      
	\qed  
	
\end{prop}
\subsection{Crossed product calculus}
In this section we study the relevant example of crossed product calculi. We recall that Theorem \ref{Doi-Takeuchi} provides a $1$ to $1$ correspondence between crossed product algebras and cleft extensions $B\subseteq B\sharp_{\sigma}H$.  Since every cleft extension is also a faithfully flat Hopf-Galois extension, we conclude every extension $B\subseteq B\sharp_{\sigma}H$ is a quantum principal bundle.  Classically a cleft extension can be thought as a trivialisation of the bundle. 

In section \ref{crossed product algebras} we considered $B$ a $\sigma-$twisted left $H-$module algebra,  and a Hopf algebra $H$. We introduced the crossed product algebra as the tensor product $B\otimes H$ equipped with an associative unital product $\mu_{\sharp_{\sigma}}:(B\otimes H)\otimes (B\otimes H)\rightarrow(B\otimes H)$. For details we remind to Lemma \ref{crossed product algebras are associative unital}. 

\begin{defn}[\cite{sciandra2023noncommutative}, Definition 3.3]\label{cross product calculus}
	Let $(\Omega^{1}(B),\dd_{B})$ be a first order differential calculus on a $\sigma-$twisted $H-$module algebra $B$ with measure $\cdot:H\otimes B\rightarrow B$ and 2-cocycle $\sigma:H\otimes H \rightarrow B$. We say that $(\Omega^{1}(B),\dd_{B})$ is a $\sigma-$twisted $H-$module differential calculus if there exists a linear map $\cdot :H\otimes \Omega^{1}(B)\rightarrow\Omega^{1}(B)$ such that \eqa{ h\cdot (b\dd_{B}b') & = (h_{1}\cdot b)(h_{2}\cdot \dd_{B}b') \\  \dd_{B}(h\cdot b)& =h\cdot \dd_{B}b \\ \dd_{B}\circ \sigma & = 0,} for every $h\in H$ and $b,b'\in B$. 
\end{defn}
\begin{thm}[\cite{sciandra2023noncommutative},Theorem 3.7] Given $(\Omega^{1}(H),\dd_{h})$ a bicovariant first order differential calculus on $H$ and $(\Omega^{1}(B),\dd_{B})$ a $\sigma-$twisted $H-$module first order differential calculus on $B$, we obtain a right $H-$covariant first order differential calculus $(\Omega^{1}(B\sharp_{\sigma}H),\dd_{\sharp_{\sigma}})$ on $B\sharp_{\sigma}H$, where \eqa{& \Omega^{1}(B\sharp_{\sigma}H):=(\Omega^{1}(B)\otimes H)\oplus(B\otimes \Omega^{1}(H)), \\ & \dd_{\sharp_{\sigma}}:B\sharp_{\sigma}H\rightarrow \Omega^{1}(B\sharp_{\sigma}H),\quad b\otimes h \mapsto \dd_{B}b\otimes h + b\otimes \dd_{H}h. }
	\qed 
\end{thm}

\begin{rmk} In the proof of this Theorem an explicit form for the right $H-$coaction on $\Omega^{1}(B\sharp_{\sigma}H)$ is provided. Since $B\sharp_{\sigma}H$ is a right $H-$comodule algebra with \eqa{\Delta_{B\sharp_{\sigma}H}:B\sharp_{\sigma}H\rightarrow B\sharp_{\sigma}H\otimes H,\quad b\otimes h \mapsto b\otimes h_{1}\otimes h_{2},} we define \eqa{\Delta_{B\sharp_{\sigma}H}^{1}:\Omega^{1} (B\sharp_{\sigma}H) &\rightarrow \Omega^{1}( (B\sharp_{\sigma}H) \otimes H), \\  \omega\otimes h + b\otimes \eta& \mapsto \omega\otimes h_{1}\otimes h_{2} + b\otimes \eta_{0}\otimes \eta_{1} + b \otimes \eta_{-1}\otimes \eta_{0}.}  
This map is differentiable, indeed \eqa{ \Delta_{B\sharp_{\sigma}H}^{1}\circ \dd_{\sharp_{\sigma}}(b\otimes h) & = \Delta_{B\sharp_{\sigma}H}^{1}(\dd_{B}b\otimes h +b\otimes \dd_{H}h) \\ & = \Delta_{B\sharp_{\sigma}H}^{1}(\dd_{B}b\otimes h) + \Delta_{B\sharp_{\sigma}H}^{1}(b\otimes \dd_{H}h) \\ & = \dd_{B}b\otimes h_{1}\otimes h_{2} + b\otimes (\dd_{H}h)_{0}\otimes (\dd_{H}h)_{1} + b \otimes (\dd_{H}h)_{-1}\otimes (\dd_{H}h)_{0} \\ & = \dd_{B}b\otimes h_{1}\otimes h_{2} + b\otimes \dd_{H}h_{1}\otimes h_{2} + b\otimes h_{1}\otimes \dd_{H}h_{2} \\ & = \dd_{\sharp_{\sigma}}(b\otimes h_{1})\otimes h_{2} + (b\otimes h_{1})\otimes \dd_{H}h_{2} \\ & = \dd_{\sharp_{\sigma}}\circ \Delta_{B\sharp_{\sigma}H}(b\otimes h),} where we exploited left/right colinearity of differential on $H$. Therefore $\Omega^{1}(B\sharp_{\sigma}H)$ is a right $H-$covariant first order differential calculus over $B\sharp_{\sigma}H$, and the vertical map $\ver=\ver^{0,1}$ is well defined.  \qed  \end{rmk}

We generalise Definition \ref{cross product calculus} to higher order forms on the crossed product algebra $B\sharp_{\sigma}H$. 

\begin{defn}[\cite{sciandra2023noncommutative},Definition 3.13] Let $(\Omega^{\bullet}(B),\dd_{B})$ be a differential calculus on $B$. We say that it is a $\sigma-$twisted $H-$module differential calculus if there exists linear maps $\cdot:H\otimes \Omega^{k}(B)\rightarrow \Omega^{k}(B)$, for all $k\geq 1$, such that\eqa{ h\cdot (b^{0}\dd_{B}b^{1}\wedge \cdots \wedge \dd_{B}b^{k}) & = (h_{1}\cdot b^{0})(h_{2}\cdot \dd_{B}b^{1})\wedge\cdots \wedge(h_{k+1}\cdot \dd_{B}b^{k}), \\ h \cdot \dd_{B}b & = \dd_{B}(h\cdot b), \\ \dd_{B}\circ \sigma & = 0,} for all $b,b^{0},\cdots, b^{k}\in B$ and $h\in H$.  
	
\end{defn}

\begin{thm}\label{DC on CPA are complete}[\cite{sciandra2023noncommutative},Theorem 3.15] Let $(\Omega^{\bullet}(B),\dd_{B})$ be a $\sigma-$twisted $H-$module differential calculus over $B$, and let $(\Omega^{\bullet}(H),\dd_{H})$ be a bicovariant differential calculus over $H$. Let us define \eqa{\Omega^{n}(B\sharp_{\sigma}H):= \bigoplus_{i=0}^{n}\Omega^{n-i}(B)\otimes \Omega^{i}(H),} for all $n\geq 0$, and let \eqa{ (\omega\otimes \eta)\wedge(\omega'\otimes \eta')& := (-1)^{jk}(\omega \wedge(\eta_{-2}\cdot \omega')\sigma(\eta_{-1}\otimes \eta'_{-1}))\otimes (\eta_{0}\wedge \eta'_{0}), \\ \dd_{\sharp_{\sigma}}(\omega\otimes \eta) & := \dd_{B}\omega \otimes \eta + (-1)^{i}\omega\otimes \dd_{H}\eta,} for $\omega\in \Omega^{i}(B)$, $\eta \in \Omega^{j}(H)$ and $\omega'\in \Omega^{k}(B)$, $\eta'\in \Omega^{\ell}(H).$ Then $(\Omega^{\bullet}(B\sharp_{\sigma}H,\dd_{\sharp_{\sigma}})$ is a right $H-$covariant differential calculus on $B\sharp_{\sigma}H$ with respect to which the right $H-$coaction $\Delta_{B\sharp_{\sigma}H}$ is differentiable.  \qed 
	
\end{thm}

 \begin{prop} 
Let $\Omega^{\bullet}(H)$ be a complete differential calculus and let $\Omega^{\bullet}(B\sharp_{\sigma}H)$ be the corresponding (complete) crossed product calculus. Differential forms over the base space are the ones $\Omega^\bullet(B)\otimes 1 \cong \Omega^\bullet(B). $

\proof  Every differential calculus defined as in Theorem \ref{DC on CPA are complete} is complete since the right $H-$coaction extends to higher order forms as a morphism $\Delta_{B\sharp_{\sigma}H}^{\wedge}:\Omega^{\bullet}(B\sharp_{\sigma}H)\rightarrow \Omega^{\bullet}(B\sharp_{\sigma}H)\otimes \Omega^{\bullet}(H).$ Let  $\beta\otimes \gamma$ be a coinvariant form in $ \Omega^{\bullet}(B\sharp_{\sigma}H)$. We have \eqa{ \Delta_{B\sharp_{\sigma}H}^{\bullet}(\beta\otimes \gamma) & = \beta\otimes \gamma_{[1]}\otimes \gamma_{[2]} \\ & = \beta\otimes \gamma\otimes 1, } if and only if $\Delta^{\bullet}(\gamma)=\gamma\otimes 1.$ If $\gamma\in \Omega^{0}(H)=H$, than $\gamma$ must be a scalar multiple of the unit.

On the other hand, given a $k-$form $\gamma=h^0\dd h^{1}\wedge\cdots\wedge \dd h^{k}$ on $\Omega^\bullet(H)$, we have  
\eqa{ h^0_1\dots h^k_1 \otimes h^0_2 \dd h^1_2 \wedge\cdots\wedge \dd h^k_2 = 0.}  

Applying $(\epsilon \otimes \id)$ to the above equation gives \eqa{\gamma = h^0 \dd h^1 \wedge\cdots \wedge \dd h^k = 0.}  \qed

 \end{prop}
 
\cleardoublepage

  \end{document}